\documentclass[11pt]{article}

\usepackage{style}

\usepackage[title]{appendix}

\usepackage{tikz-cd} 
\usepackage{mathrsfs}
\usepackage{dynkin-diagrams}
\usepackage{mathdots}
\usepackage{enumitem}
\usepackage{stmaryrd}
\usepackage{multirow}

\usepackage{array} 
\usepackage{bm}
\usepackage[font=small]{caption}

\newcommand{\CD}{\mathrm{CD}}
\numberwithin{equation}{section}

\newcommand{\pr}{\mathrm{pr}}
\newcommand{\CP}{\mathbb{CP}}
\newcommand{\Vect}{\mathrm{Vec}}
\newcommand{\KO}{\mathrm{KO}}
\newcommand{\solie}{\mathfrak{so}}
\newcommand{\CH}{\mathbb{CH}}
\newcommand{\K}{\mathbb{K}}

\newcommand{\Cl}{\mathrm{Cl}}
\newcommand{\Gtwo}{\mathrm{G}_2}
\newcommand{\Mat}{\mathrm{Mat}}
\newcommand{\image}{\mathrm{image}}

\newcommand{\Gtwosymmetricspace}{\Gr_{(3,0)}^\times(\R^{3,4})}
\newcommand{\Gtwofullflags}{\mathcal{F}_{1,2}^\times}

\newcommand{\dbrack}[1]{\llbracket #1 \rrbracket}
\newcommand{\sig}{\mathrm{sig}}
\newcommand{\Spin}{\mathrm{Spin}}
\newcommand{\Pin}{\mathrm{Pin}}
\newcommand{\Ha}{\mathbb{H}}
\newcommand{\Pho}{\mathrm{Pho}}
\newcommand{\Ein}{\mathrm{Ein}}
\newcommand{\Iso}{\mathrm{Iso}}

\newcommand{\Gtwosplit}{\mathrm{G}_2^{\mathrm{'}}}
\newcommand{\g}{\mathfrak{g}}
\newcommand{\Oct}{\mathbb{O}}
\newcommand{\Ann}{\mathrm{Ann}}
\newcommand{\graph}{\mathrm{graph}}
\newcommand{\fullflag}{\mathcal{F}_{1,2}^\times}

\title{Transverse Spheres in Flag Manifolds}
\author{Parker Evans and J. Maxwell Riestenberg  \thanks{Date: \today.}}
\date{}

\makeatletter
\def\thanks#1{\protected@xdef\@thanks{\@thanks
        \protect\footnotetext{#1}}}
\makeatother

\begin{document}

\maketitle

\abstract{ 
For some partial flag manifolds of semisimple real Lie groups, including many full flag manifolds, transverse circles are known to be locally maximally transverse. 
We complete the classification of all partial flag manifolds of split real Lie groups with this property. 
As a consequence, $\{7\}$-Anosov subgroups of split $E_7$ are virtually free or surface groups.

On the other hand, using spinors, we find transverse spheres of arbitrarily large dimension in certain full flag manifolds of Cartan-Killing types $A,B,D$. 
These transverse spheres are verified to be maximally transverse with tools from topological $K$-theory. 
The aforementioned classification follows from constructions of transverse $m$-spheres, $m \geq 2$, that complement the previously known restrictions as well as the new $E_7$ restriction.  
Additionally, when $G$ is split of type $G_2, B_3,$ or $D_4$, the full flag manifold admits a fibration by maximally transverse 3-spheres. }  
\tableofcontents

\section{Introduction}

Let $G$ be a semisimple real Lie group. 
For every subset  $\Theta \subseteq \Delta$ of simple restricted roots $\Delta$, there is an associated flag manifold $\mathcal{F}_\Theta$; see Section \ref{Sec:FlagPreliminaries}. 
\emph{Transversality} can be defined in each \emph{self-opposite} flag manifold $\mathcal{F}_{\Theta}$ as follows: a subset $S \subset \mathcal{F}_{\Theta}$ is transverse when any distinct pair of points $x,y \in S$ lie in the unique open $G$ orbit in $\mathcal{F}_{\Theta}^2$. 
This paper primarily concerns new examples of transverse subsets of flag manifolds.

Our motivation comes from the study of discrete subgroups of Lie groups. 
Over the last two decades, ($\Theta$-)\emph{Anosov subgroups} of non-compact higher rank (semisimple) real Lie groups $G$ have emerged as the appropriate generalization of convex cocompact subgroups of rank one Lie groups to the higher rank setting \cite{Lab06,GW12}. 
These are hyperbolic groups $\Gamma <G $ whose Gromov boundary $\partial \Gamma$ embeds as a transverse subset $\Lambda$ of the flag manifold $\mathcal{F}_{\Theta}$, characterized by additional dynamical properties of the action of $\Gamma$ on $\Lambda$ \cite{GGKW17a,KLP17}. 
An important fundamental problem is then to first understand transverse subsets of the flag manifolds $\mathcal{F}_\Theta$.

Transverse circles abound in every self-opposite flag manifold $\mathcal{F}_{\Theta}$ of any semisimple real Lie group $G$.
For example, up to finite cover, any such group admits at least one copy of $\SL(2,\R)$ whose limit set is a transverse circle in the associated full flag manifold.
A natural question of interest is when a flag manifold $\mathcal{F}_{\Theta}$ contains transverse subsets larger than a circle. To formalize the question, call a transverse subset $S \subset \mathcal{F}_{\Theta}$ \emph{locally maximally transverse} when it has an open neighborhood $U \subset \mathcal{F}_{\Theta}$ so that $S$ is maximally transverse in $U$. 

\begin{question}\label{q:circles maximal?}
    Fix $G$ and $\Theta$. Is every transverse circle in $\mathcal{F}_\Theta$ locally maximally transverse?
\end{question}

There are now several known cases with a positive answer to Question \ref{q:circles maximal?}, as follows:
Tsouvalas handled $\Gr_{2k+1}(\R^{4k+2})$ in \cite{tsouvalas2020borel}, Dey handled $\Flag(\R^d)$, the $\SL(d,\R)$-full flag manifold, for $d \in \{2,3,4,5,6\} \bmod 8$ \cite{Dey25}, Dey-Greenberg-R. and Pozzetti-Tsouvalas handled every $\Sp(2n,\R)$-flag manifold $\mathcal{F}_{\Theta}$, with positive answer exactly when $\Theta \cap (2\mathbb{N}+1) \neq \emptyset$ \cite{DGR24, PT24}, and Kineider-Troubat handled  $\Iso_{p}(\R^{p,p+1})$ for $p \equiv \{1,2\} \bmod 4$ and $\Iso_{p}^{\pm}(\R^{p,p})$ for $p \equiv 2 \bmod 4$ \cite{KT24}. 
Positive answers are particularly interesting, because they imply restrictions on Anosov subgroups; see the discussion below in Section \ref{Sec:Anosov}. 

\medskip

Building on the prior work mentioned in the previous paragraph,
we resolve Question \ref{q:circles maximal?} for every self-opposite flag manifold $\mathcal{F}_{\Theta}$ of every split real Lie group $G$; see Theorem \ref{Thm:IntroThmEquivalence} for a slightly sharper statement. 
\begin{theorem}
Let $\mathcal{F}_{\Theta}$ be a self-opposite flag manifold of a split real simple Lie group $G$. Transverse circles in $\mathcal{F}_{\Theta}$ are locally maximally transverse precisely for pairs $(G,\Theta)$ in Table \ref{Table:CirclesMaximal?}. 
\end{theorem}

\begin{table}
\small
\centering
\captionsetup{width=0.8\textwidth}
\renewcommand{\arraystretch}{1.2}
\begin{tabular}{|c|cc |}\hline
Cartan-Killing Type & Conditions on $n$ &  $\Theta \subseteq \Delta $ with Positive answer to Question \ref{q:circles maximal?} \\ \hline 
\multirow{3}{*}{$A_{n-1}$ $(n \geq 2)$}
    & $n \in \{-1,0,1\} \bmod 8 $ & none \\\cline{2-3}
    & $n \in \{2,6\} \bmod 8$ & $\Theta \cap \{\frac{n}{2}\} \neq \emptyset $ \\ \cline{2-3} 
    & $n \in \{3,4,5\} \bmod 8$ & $\Theta = \Delta$ \\\hline
\multirow{2}{*}{$B_n$ $(n \geq 1)$} 
    & $n \in \{1,2\} \bmod 4$ & $\Theta \cap \{n\} \neq \emptyset $ \\\cline{2-3}
    & $n \in \{0,3\} \bmod4$ & none  \\\hline
$C_n$ $(n \geq 1)$ & & $\Theta \cap (2\mathbb{N}+1) \neq \emptyset $ \\ \hline
\multirow{2}{*}{$D_n$ $(n \geq 4)$} 
    & $n \in \{0,1,3\} \bmod 4$ & none  \\\cline{2-3}
    & $n \equiv 2 \bmod 4$ & $\Theta \cap \{n^+, n^-\} \neq \emptyset $ \\\hline
$G_2, F_4, E_6, E_8$ &  & none \\ \hline 
$E_7$ & & $\{7\} \cap \Theta \neq \emptyset$ \\ \hline 
\end{tabular}
\caption{\small{The complete answer to Question \ref{q:circles maximal?} for $G$ split simple and $\mathcal{F}_{\Theta}$ a $G$-self-opposite flag manifold.
The statements are in terms of the Dynkin diagrams as in Figures \ref{fig:DynkinAp}, \ref{fig:DynkinBp}, \ref{fig:DynkinDp}, \ref{fig:dynkin diagram E7} in the $A,B,D, E_7$ cases, respectively. The only new restrictions occur for type $E_7$. The known restrictions in the $A,B,D$ cases are verified to be optimal.}}
\label{Table:CirclesMaximal?}
\end{table}

Our approach is based on a construction of transverse spheres via \emph{spinor representations}.
These spheres are then propagated to other flag manifolds with an iterated direct sum construction along with transversality-preserving maps of flag manifolds.
It turns out that the previously known positive answers to Question \ref{q:circles maximal?} indeed cover all cases, except for certain partial flag manifolds of $E_7$; see Theorem \ref{Thm:MainTheoremExceptionalTypes}. 
In every negative case, we find a transverse $m$-sphere in the given flag manifold $\mathcal{F}_{\Theta}$ for $m \geq 2$. 

Before the present work, there was no known example of a transverse subset of a full flag manifold of a split real Lie group properly containing a circle. 
In contrast, we obtain arbitrarily large maximally transverse spheres in (the possible) full flag manifolds of split real groups of Cartan-Killing type $A,B,D$. 
For example, for $d = 2\cdot 16^n$, the full flag manifold $\Flag(\R^d)$ contains a maximally transverse $(8n)$-sphere. 

It turns out that there is a connection between the striking mod $8$ phenomenon controlling transverse spheres in $\Flag(\R^d)$, as seen in \cite{Dey25} and Table \ref{Table:CirclesMaximal?}, and Bott periodicity.
This is explained, at least partly, by the Atiyah-Bott-Shapiro isomorphism \cite{ABS64} linking the representation theory of Clifford algebras to real $\mathrm{K}$-theory of spheres.
We leverage this relationship to prove that many of the transverse spheres we construct in full flag manifolds are, in fact, maximally transverse. 
Direct sum constructions of maximally transverse spheres often remain maximally transverse, 
as witnessed by stable homotopy classes of vector bundles,
and so we are able to produce many different maximally transverse spheres of varying dimensions in the same full flag manifold. 

In the remainder of the introduction, we explain our techniques, main results, and their consequences in more detail. 
Broadly speaking, our goals are twofold: to \emph{construct} and \emph{obstruct} transverse spheres in flag manifolds of split real Lie groups. 
As it turns out, the spheres produced by the most basic version of the construction, as well as some of their cousins in other full flag manifolds, cannot arise as Anosov limit sets, see Section \ref{Sec:Anosov}. 
On the other hand, by deforming the original construction, one obtains other maximally transverse spheres that we cannot currently obstruct from being Anosov limit sets. 

\subsection{Main Results} 
We state our results according to the Cartan-Killing type of the split real Lie group $G$. Due to the connections with Clifford algebras and spinors, the \emph{Radon-Hurwitz} numbers $\rho(d)$ will appear in nearly all of our theorem statements. Precisely, factorize the positive integer $d$ as $d = (2\ell+1)16^m2^j$, for integers $\ell, m \geq 0$ and $0\leq j\leq 3$, and then  
\begin{align}\label{RadonHurwitz}
   \rho(d) \coloneqq 8m +2^j.
\end{align} 
Observe that $k \ge 1$ implies $\rho(4k) \ge 4$. 
\medskip

We begin with the type $A$ case. 
The self-opposite partial flag manifolds associated to $\SL(d,\R)$ are of the form $\Flag_\Theta(\R^d)$, where $\Theta$ is a \emph{symmetric} subset of $\dbrack{d} \coloneqq \{0,1,\dots,d-1,d\}$, i.e.\ $k$ is in $\Theta$ if and only if $d-k$ is in $\Theta$, see Section \ref{Subsec:TypeAFlags}. 
In particular, the full flag manifold is denoted $\Flag(\R^d) = \Flag_{\dbrack{d}}(\R^d)$. 

\begin{theorem}[Type A {[Theorem \ref{Thm:AnTransverseSpheres}]}]
\label{Thm:MainTheoremTypeA}
    Fix an integer $d \ge 2$. 
    Write $d= 8k+\epsilon$, for $k \geq 0$ and $\epsilon \in \{-1,0,1,2, \dots, 6\}$. Suppose: 
    \begin{enumerate}[label=(\alph*)]
        \item $\epsilon \in \{-1,0,1\}$. Write $4k =n2^j$ for some integers $j\geq 2$, $n \geq 1$, such that $j \in \{2,3\} \bmod 4$ or $n$ is odd. Then there is a maximally transverse $(\rho(2^j)-1)$-sphere in $\Flag(\R^d)$. 
        \item $\epsilon=2$. If $\frac{d}{2} \in \Theta$, then every transverse circle in $\Flag_{\Theta}(\R^d)$ is maximally transverse. Otherwise, $\Theta$ is a symmetric subset of $\dbrack{d} \setminus \{\frac{d}{2}\}$ and there exists a transverse $(\rho(4k)-1)$-sphere in $\Flag_{\Theta}(\R^d)$. 
        \item $\epsilon \in \{3,4,5\}$. Every transverse circle in $\Flag(\R^d)$ is locally maximally transverse. Otherwise, $\Theta$ is a proper symmetric subset of $\dbrack{d}$, and there exists a transverse $2$-sphere in $\Flag_\Theta(\R^d)$.
        \item $\epsilon =6$. If $\frac{d}{2} \in \Theta$, then every transverse circle in $\Flag_{\Theta}(\R^d)$ is maximally transverse. Otherwise, $\Theta$ is a symmetric subset of $\dbrack{d} \setminus \{\frac{d}{2}\}$ and there exists a transverse $(\rho(4k+4)-1)$-sphere in $\Flag_{\Theta}(\R^d)$.      
    \end{enumerate}

\end{theorem}

We note that the assertions about maximally transverse circles in (b), (d) of Theorem \ref{Thm:MainTheoremTypeA} follow directly from the proof of \cite[Theorem 1.1]{tsouvalas2020borel} and the assertions about locally maximally transverse circles in (c) are a special case of \cite[Corollary B]{Dey25}. 
The fact that a transverse $3$-sphere in $\Flag_\Theta(\R^7)$ (resp.\ transverse 7-sphere in $\Flag_{\Theta}(\R^{15})$) must be maximally transverse when $\{3,4\}\subset \Theta$ (resp.\ $\{7,8\} \in \Theta$) is a special case of \cite[Lemma 5.4]{TZ24}. \medskip 

Next we consider $G=\SO_0(p,p+1)$. 
Its associated partial flag manifolds $\mathcal{F}_\Theta = \Iso_{\Theta}(\R^{p,p+1})$ consist of chains of isotropic subspaces, as recalled in Subsection \ref{Subsec:TypeBFlags}. 
The full flag manifold is $\mathcal{F}_{\Delta} = \Iso_{\llbracket p \rrbracket}(\R^{p,p+1})$, where $\dbrack{p} = \{1,2,\dots, p\}$. 

\begin{theorem}[Type B {[Theorem \ref{thm:BnTransverseSpheres}]}]\label{Thm:MainTheoremTypeB}
  Let $p \geq 1$.  Suppose: 
  \begin{enumerate}[label = (\alph*)]
      \item $p \equiv 0$ mod 4. 
      For any integer $2^j$, with $j \geq 2$, dividing $p$, there exists a maximally transverse $(\rho(2^j)-1)$-sphere in the full flag manifold $\Iso_{\dbrack{p}}(\R^{p,p+1})$.
      \item $p \equiv 1$ mod 4. If $p \in \Theta$, then every transverse circle in $\Iso_{\Theta}(\R^{p,p+1})$ is locally maximally transverse. Otherwise, $\Theta \subseteq \dbrack{p-1}$, and there exists a transverse $(\rho(p-1)-1)$-sphere in $\Iso_{\Theta}(\R^{p,p+1})$. 
      \item $p \equiv 2$ mod 4. If $p \in \Theta$, then every transverse circle in $\Iso_{\Theta}(\R^{p,p+1})$ is locally maximally transverse. Otherwise, $\Theta \subseteq \dbrack{p-1}$, and there exists a transverse 2-sphere in $\Iso_{\Theta}(\R^{p,p+1})$.
      \item $p \equiv 3$ mod 4. For any integer $2^j$, with $j \geq 2$, dividing $p+1$, there exists a maximally transverse $(\rho(2^j)-1)$-sphere in the full flag manifold $\Iso_{\dbrack{p}}(\R^{p,p+1})$. 
  \end{enumerate} 
  \end{theorem}

The assertions about locally maximally transverse circles in (b), (c) were proven by \cite{KT24}. \medskip 

Question \ref{q:circles maximal?} has already been resolved for type $C$. 
The answer was recently obtained by Dey-Greenberg-R.\ \cite{DGR24}, with overlap in the independent work of Pozzetti-Tsouvalas \cite{PT24}. 
The $\Sp(2n,\R)$-partial flag manifolds $\mathcal{F}_\Theta$ consist of chains of (symplectic-)isotropic subspaces, where $\Theta$ lists the dimensions of the subspaces. 
When $\Theta$ contains only even integers, there exists a transverse $2$-sphere, which is actually the limit set of a rank $1$ sub-symmetric space isometric to $\mathbb{H}^3$.
On the other hand, if $\Theta$ contains an odd integer, then every transverse circle in $\mathcal{F}_\Theta$ is locally maximally transverse; in fact, Pozzetti-Tsouvalas proved such circles are maximally transverse \cite{PT24}.
\medskip

The type $D$ case corresponds to $G=\SO_0(p,p)$.  
We write $\Delta = \{1,\dots,p-2,p^+,p^-\}$, where the labels correspond to the nodes of the Dynkin diagram as in Figure \ref{fig:DynkinDp}. 
There are natural identifications $\mathcal{F}_\Delta =\Iso_{\dbrack{p-1}}(\R^{p,p})$ and  $\mathcal{F}_{\Delta \backslash \{p^+,p^-\}} =\Iso_{\dbrack{p-2}}(\R^{p,p}) $; see Subsection \ref{Subsec:TypeDFlags}. 

\begin{theorem}[Type D {[Theorem \ref{Thm:DnTransverseSpheres}]}]\label{Thm:MainTheoremTypeD}
Let $p \geq 4$. 
Suppose:
    \begin{enumerate}[label =(\alph*)]
    \item $p \equiv 0 \bmod 4$. 
    For any integer $2^j$, with $j\geq 2$, that divides $p$, there exists a maximally transverse $(\rho(2^j)-1)$-sphere in the full flag manifold $\Iso_{\dbrack{p-1}}(\R^{p,p})$.
    \item $p \equiv 1 \bmod 4$. 
    For any integer $2^j$, with $j\geq 2$, that divides $p-1$, there exists a maximally transverse $(\rho(2^j)-1)$-sphere in the full flag manifold $\Iso_{\dbrack{p-1}}(\R^{p,p})$. 
    \item $p \equiv 2 \bmod 4$. If $\Theta \subseteq \Delta$ satisfies $\{p^+,p^-\} \cap \Theta \neq \emptyset$, then every transverse circle in $\mathcal{F}_{\Theta}$ is locally maximally transverse. Otherwise, $\Theta \subseteq \dbrack{p-2}$, and there is a transverse $(\rho(p-2)-1)$-sphere in $\Iso_{\Theta}(\R^{p,p})$.
    \item $p \equiv 3 \bmod 4$. For any integer $2^j$, with $j \geq 2$, that divides $p+1$, as long as $2^j \neq p+1$, there is a maximally transverse $(\rho(2^j)-1)$-sphere in the full flag manifold $\Iso_{\dbrack{p-1}}(\R^{p,p})$.
    \end{enumerate}
  \end{theorem}

The local maximal transversality of transverse circles in case (b) was proven by Kineider-Troubat \cite{KT24}. \medskip 

Finally, we discuss the flag manifolds of split real Lie groups of exceptional type. 
Part of the statement is in terms of the labeling of the $E_7$ Dynkin diagram in Figure \ref{fig:dynkin diagram E7}. 
Note that we have stronger results in the $G_2$ case, as stated in Section \ref{Intro:Fibrations}. 

\begin{theorem}[Exceptional cases]\label{Thm:MainTheoremExceptionalTypes}
    Let $\mathcal{F}(G)$ be the full flag manifold of $G$. \begin{enumerate}[label=(\alph*)]
        \item There is a maximally transverse $3$-sphere in $\mathcal{F}(G)$ for $G \in \{G_2,F_4,E_6\}$.
        \item There is a maximally transverse $7$-sphere in $\mathcal{F}(E_8)$. 
        \item\label{item:E7sphere} Let $G = E_7$.
        There exists a transverse $3$-sphere in $\mathcal{F}_{\Delta \setminus \{7\}}$.
        On the other hand, if $\Theta$ contains $\{7\}$, then every transverse circle in $\mathcal{F}_\Theta$ is maximally transverse. 
    \end{enumerate}
\end{theorem}

In the $E_7$ case, opposed to the other exceptional cases, it is not possible to build a higher-dimensional transverse sphere in the full flag manifold. 
A new restriction on Anosov subgroups follows from Theorem \ref{Thm:MainTheoremExceptionalTypes}\ref{item:E7sphere} and standard techniques; see Section \ref{Sec:Anosov} and Question \ref{q:Sambarino}. 

\begin{corollary}
    If $\Gamma < E_7$ is $\{7\}$-Anosov, then $\Gamma$ is virtually a free or surface group.
    In particular, if $\Gamma < E_7$ is Borel Anosov, then $\Gamma$ is virtually a free or surface group. 
\end{corollary}

\subsubsection{Further Transverse Spheres}

The transverse spheres we construct in flag manifolds of split real Lie groups can be further embedded in flag manifolds of other semisimple Lie groups.
An interesting special case is the inclusion of a split real form into a complex simple Lie group $G' \hookrightarrow G^{\mathbb{C}}$. 
While there are obvious examples of transverse $2$-spheres in $\Flag(G^{\mathbb{C}})$ which arise as limit sets of any principal $(P)\SL(2,\C) \hookrightarrow G^{\mathbb{C}}$, our construction yields many new examples of transverse $m$-spheres, $m \ge 3$, in the full flag manifold of $G^{\mathbb{C}}$.
For example, this applies to the full flag manifold $\Flag(\C^d)$ of $\SL(d,\C)$ for $d \in \{-1,0,1\} \bmod 8$. 

More generally, any semisimple real Lie algebra $\mathfrak{g}$ admits a unique \emph{maximal split real subalgebra $\mathfrak{g'} < \mathfrak{g}$} \cite{KR71} which satisfies $\rank_{\R}(\g) = \rank_{\R}(\g')$.
Combining that inclusion with our construction provides traction for addressing Question \ref{q:circles maximal?} in the general case.
We do not attempt a complete answer to Question \ref{q:circles maximal?} for the general (i.e.\ non-split) case in the present paper.

\subsection{Clifford Algebras, Spinors, and the Atiyah-Bott-Shapiro Isomorphism}
In this portion of the introduction, we provide an overview of our approach for constructing transverse spheres with spinors, proving maximal transversality, and finding transverse deformations of the original construction. 

\subsubsection{Transverse spheres from spinors}\label{Subsec:IntroSpinorSpheres}

Our main theorems rely on an intermediary construction of transverse spheres, Theorem \ref{Thm:TransverseSpinorSphere}, which mildly generalizes the following result.

\begin{theorem}[Arbitrarily Large Transverse  Spheres]\label{Thm:IntroArbitrarilyLargeSpheres}
Let $n \ge 2$ and define $d=d(n)$ as follows: write $n = 8k+r$, for $k \ge 0$ and $1\leq r \leq 8$, and set
\begin{equation}\label{eqn:definition of d(n)}
    d(n) \coloneqq \begin{cases}
    16^k,  & r =1\\
     2\cdot 16^k,  & r = 2 \\
      4\cdot 16^k,  & r \in \{3,4 \}\\
       8\cdot 16^k,  & r \in \{5,6,7,8\}.\\ \end{cases}
\end{equation}
There is a transverse $(n-1)$-sphere in the full flag manifold $\Iso_{\dbrack{d-1}}(\R^{d-1,d})$. 
\end{theorem}

By transversality-preserving embeddings of full flag manifolds, we obtain a transverse $(n-1)$-sphere in the  
$A_{2d-2}$, $A_{2d-1}$, $A_{2d}$, $B_{d}$, $D_{d}$ full flag manifolds as a corollary.
These transverse spheres then propagate to other partial flag manifolds of classical and exceptional type 
via transversality-preserving embeddings of flag manifolds and direct sum constructions for flags. \medskip

We briefly sketch the proof of Theorem \ref{Thm:IntroArbitrarilyLargeSpheres}.

Recall that $\Cl(n) \coloneqq \Cl(V)$, for any Euclidean vector space $(V,q) \cong \R^{n,0}$, is the $\R$-algebra obtained from the tensor algebra $T(V) = \bigoplus_{i=0}^{\infty} V^{\otimes k}$ by imposing the relations $v\otimes v = -q(v)1$. 
Products in $\Cl(n)$ are denoted by juxtaposition. 
Then $\Spin(n) < \Cl(n)$ is the multiplicative subgroup consisting of all elements of the form $v_1v_2\cdots v_{2k} \in \Cl(n)$, where $v_{i} \in Q_+(V,q)$ are each unit elements. 
By definition of $d\coloneqq d(n)$ in Theorem \ref{Thm:IntroArbitrarilyLargeSpheres} and \cite{LM89}, there is an irreducible $d$-dimensional representation $\Spin(n) \rightarrow \SO(d)$; see Section \ref{Sec:SpheresWithSpinors}. 
Let $S_n^+$ be an irreducible $d$-dimensional $\Spin(n)$-module, equipped with a $\Spin(n)$-invariant Euclidean metric.  
\medskip 

We choose a splitting $\R^{d-1,d} = \R^{d-1,0} \oplus \R^{0,d}$ and any complete isotropic full flag $F \in \Iso_{\dbrack{d-1}}(\R^{d-1,d})$. Identify $\R^{0,d}$ with $S_n^+$, the latter with its signature flipped. 
There are orthonormal bases $\{e_i\}$ of $\R^{d-1,0}$ and $\{s_i\}$ of $\R^{0,d}$ so that $F^k = \spann \{e_i+s_i\}_{1\le i \le k}$. 
By definition of $\Spin(n)$ via the Clifford algebra of $\R^{n,0}$, there is a natural injective map $f \colon \mathbb{S}^{n-1} \to \Spin(n)$ given by $f(x)= x_0x $, for any $x_0 \in Q_+(\R^{n,0})$. 
The transverse sphere $\Lambda$ in $\Iso_{\dbrack{d-1}}(\R^{d-1,d})$ is then given by the `orbit' of the sphere $\image(f) \subset \Spin(n)$ on the basepoint flag $F$.
Explicitly: 
\begin{align}\label{SpinorSpheresBasic}
    \Lambda = \left\{ V^{\bullet}(x) \in \Iso_{\dbrack{d-1}}(\R^{d-1,d}) \mid x \in \mathbb{S}^{n-1}, V^i(x) = \spann \{ e_1+f(x)s_1, \, \dots, \,e_i+f(x)s_i \}\right\}
\end{align} 
We emphasize that the above construction works for any chosen basepoint flag $F$. \medskip 

In the special case that $n=3$, the sphere $\mathbb{S}^3$ is a Lie group. 
In the $G_2,B_3,$ and $D_4$ cases, every orbit of the $\mathbb{S}^3\cong\Sp(1)$-action is a maximally transverse 3-sphere.
In fact, these orbits are the fibers of a principal $\Sp(1)$-fibration of the corresponding full flag manifolds, as we describe in Section \ref{Intro:Fibrations}. 

\subsubsection{Maximally Transverse Spheres}

In fact, the spheres in Theorem \ref{Thm:IntroArbitrarilyLargeSpheres} are \emph{maximally} transverse in certain cases, depending on $n \bmod 8$.  
This can be interpreted as a manifestation of Bott periodicity via the Atiyah-Bott-Shapiro isomorphism.

\begin{theorem}[{Theorem \ref{Thm:ABSMaximallyTransverse}}]\label{thm:introthmmaximality}
The transverse $(n-1)$-spheres $\Lambda$ in Theorem \ref{Thm:IntroArbitrarilyLargeSpheres} are maximally transverse if and only if $n \in \{0,1,2,4\} \bmod 8$. 
\end{theorem}

We now briefly sketch the proof of this result. 
Note that when $n \in \{3,5,6,7\} \bmod 8$, the transverse sphere $\Lambda$ sits in a larger transverse sphere in the same full flag manifold by the proof of Theorem \ref{Thm:IntroArbitrarilyLargeSpheres}, so $\Lambda$ is not maximally transverse. \medskip  

Now, suppose $n \in \{0,1,2,4\} \bmod 8$, and we summarize the proof of maximal transversality. 
An observation due to Dey is that a subset $S \subset \mathcal{F}$ that is transverse and homotopically non-trivial in $\mathcal{F}$ must be maximally transverse. 
Indeed, by contrapositive, when $S$ is not maximally transverse there a point $x \in \mathcal{F}$ transverse to $S$. 
Then $S$ is contained in the open Schubert cell $\mathscr{C}_{x} = \{ F \in \mathcal{F} \mid x \pitchfork F\}$ and must be contractible in $\mathcal{F}$.  

To prove the spheres $\Lambda$ in Theorem \ref{Thm:IntroArbitrarilyLargeSpheres} are maximally transverse, it then suffices to show they are homotopically non-trivial in $\Flag(G)$. For $G$ split simple, the maximal compact subgroup $K <G $ is a finite cover of $\Flag(G)$, so it suffices to show the lift of $\Lambda$ to $K$ is homotopically non-trivial. 
We sketch the $B_{d-1}$ case for concreteness. 
In this case, $K = \SO(d-1) \times \SO(d)$ and the lift obtains the form $(\id, f)$. The map $f: \mathbb{S}^{n-1} \rightarrow \SO(d)$ can be verified to be homotopically non-trivial using an isomorphism defined by Atiyah-Bott-Shapiro in \cite{ABS64}. 
The idea is as follows. The map $f$ is associated to an irreducible representation $\eta: \Spin(n) \rightarrow \SO(d)$. 
Associated to $f$ is a rank $d$ vector bundle $E_f$ over $\mathbb{S}^{n}$, using $f$ as clutching data. Then $f$ is homotopically trivial if and only if $E_f$ is a topologically trivial vector bundle. 
The Atiyah-Bott-Shapiro isomorphism allows us to determine whether $E_f$ is non-trivial (even stably non-trivial) in terms of representation-theoretic conditions on $\eta$. 
The latter conditions hold exactly when $n \in \{0,1,2,4\} \bmod 8$. 
In this way, $[f]\neq 0 \in \pi_{n-1}\SO(d)$ is confirmed and the associated sphere $\Lambda \subset \Flag(\R^{d-1,d})$ is maximally transverse.  

\medskip

In low dimensions, the spheres we consider satisfy a stronger form of maximal transversality: the projection to any partial flag manifold remains maximally transverse. 

\begin{theorem}[{Theorem \ref{Cor:CompAlgMaxTransverlity}}]\label{thm:introMaxlTransversalityEveryFlag}
    Let $a \in \{2,4,8\}$ and $\epsilon \in \{-1,0\}$ and $\mathrm{pr}_{\Theta} \colon \Flag(\R^{a+\epsilon,a}) \rightarrow \Iso_{\Theta}(\R^{a+\epsilon,a})$ be the natural projection to any $\SO_0(a+\epsilon, a)$-flag manifold.
    If $\Lambda \subset \Flag(\R^{a+\epsilon,a})$ is a transverse $(a-1)$-sphere from Corollary \ref{Thm:IntroArbitrarilyLargeSpheres}, then $\mathrm{pr}_{\Theta}(\Lambda)$ is maximally transverse.
\end{theorem}

This property is related to the fact that $\mathbb{S}^{n-1}$ is an $H$-space exactly when $n \in \{1,2,4,8\}$ \cite{Ada60}.  
More precisely, we use that the map $f\colon \mathbb{S}^{n-1} \rightarrow \SO(d)$ by $f(x) = \eta(x_0x)$, has some/every orbit map $\mathbb{S}^{n-1} \rightarrow \mathbb{S}^{d-1}$, via $x \mapsto f(x)\cdot y_0$, a homeomorphism exactly when $n=d$, in which case $\mathbb{S}^{n-1}$ is equipped with a kind of multiplication map. 
The stronger maximal transversality result in Theorem \ref{thm:introMaxlTransversalityEveryFlag} does not hold when $n < d-1$. 
Indeed, the projection of an $(n-1)$-sphere $\Lambda \subset \Flag(\R^{d-\epsilon, d})$ to $\Ein(\R^{d+\epsilon,d})$ cannot be maximally transverse when $n < d-1$.  

Atiyah-Bott-Shapiro also proved a rather amazing fact in \cite{ABS64}: every vector bundle over $\mathbb{S}^{k+1}$ is stably isomorphic to a `spinor bundle' constructed via clutching function $\eta: \mathbb{S}^k \rightarrow \mathrm{O}(d)$, where $\eta$ is the restriction of a Clifford algebra representation $\eta: \Cl(\R^{k+1,0}) \rightarrow \End(\R^d)$ to the unit sphere $Q_+(\R^{k+1,0})\cong \mathbb{S}^k$. In other words, using the isomorphism $\widetilde{\KO}(\mathbb{S}^{k+1})\cong \pi_k(\mathbf{O})$, they prove every non-trivial stable homotopy class $\omega \in \pi_k(\mathbf{O})$ can be represented by a `spinor sphere' $\eta: \mathbb{S}^k \rightarrow \mathrm{O}(d)$ as above. We find the following shadow cast by the deep work of \cite{ABS64}.   

\begin{theorem}
Let $0 \neq \omega \in \pi_k(\mathbf{O})$. Then $\omega =[f] $ for a representative $f: \mathbb{S}^k \rightarrow \SO(d)$ such that the orbit map $\xi: \mathbb{S}^k \rightarrow \Flag(\R^d)$ by $\xi(x) = f(x) \cdot F$, for some flag $F \in \Flag(\R^d)$, is a maximally transverse sphere in $\Flag(\R^d)$. 
\end{theorem}

\subsubsection{Deformations}\label{Subsec:Deform+Obstruct}

The transverse spheres built from spinors in Section \ref{Subsec:IntroSpinorSpheres} are algebraic, but admit `generic' deformations.

Associated to a strictly 1-Lipschitz map $\phi \colon \mathbb{S}^{n-1} \rightarrow \mathbb{S}^{n-1}$ we obtain a modified transverse sphere in $\Flag(\R^{d,d})$ by now allowing $\Spin(n)$ to act on both factors of $\R^{d,d}$. 
Recall the notation $f \colon \mathbb{S}^{n-1} \rightarrow \SO(d)$ with $f(x) =\eta(
x_0x)$ defined via an irreducible representation $\eta \colon \Spin(n) \rightarrow \SO(d)$. 
We consider the sphere $\Lambda' \subset \Flag(\R^{d,d})$ given by 
\begin{align}\label{SpinorSpheresDeformed}
    \Lambda' = \{ V^{\bullet} \in \Flag(\R^{d,d}) \mid x \in \mathbb{S}^{n-1}, V^i = \spann \{ f(\phi(x))e_1+f(x)s_1, \, \dots, \,f(\phi(x))e_i+f(x)s_i \}
\}.
\end{align} 
In this light, the sphere $\Lambda$ in \eqref{SpinorSpheresBasic} corresponds to the constant map $\phi(x) = -x_0$. 

\begin{lemma}\label{thm:IntroDeformedSphere}
If $\phi \colon \mathbb{S}^{n-1} \rightarrow \mathbb{S}^{n-1}$ is a strictly 1-Lipschitz map, then $\Lambda'$ in \eqref{SpinorSpheresDeformed} is a maximally transverse $(n-1)$-sphere in $\Flag(\R^{d,d})$.  
\end{lemma}
Using well-chosen contracting maps $\phi$ in Lemma \ref{thm:IntroDeformedSphere}, along with some further geometric arguments, we are able to produce maximally transverse $(n-1)$-spheres, for $n \in \{4,8\}$ in any full flag manifold $\Flag(\R^{d,d})$, where $d$ is divisible by $4$ or $8$, whose projection to the Einstein universe is linearly full.
\begin{theorem}\label{thm:introthmfullsphere}
    Let $n \in \{4,8\}$ and $d \equiv 0 \bmod n$. 
    There is a maximally transverse $(n-1)$-sphere $\Lambda' \subset \Flag(\R^{d,d})$ such that $\spann(\pr_1(\Lambda') ) = \R^{d,d}$, for $\pr_1\colon \Flag(\R^{d,d}) \rightarrow \Ein^{d-1,d-1}$ the natural projection. 
    Moreover, there is a maximally transverse $(n-1)$-sphere $\Lambda' \subset \Flag(\R^{2d})$ such that $\spann(\pr_1(\Lambda') ) = \R^{2d}$, for $\pr_1\colon \Flag(\R^{2d}) \rightarrow \mathbb{RP}^{2d-1}$ the natural projection. 
\end{theorem}

\subsection{Fibrations by Maximally Transverse Spheres}\label{Intro:Fibrations}

In this section, we discuss the full flag manifolds in the $G_2, B_3, D_4$ cases. These full flag manifolds admit principal $\mathrm{Sp}(1)$-bundle fibrations such that every fiber is a maximally transverse 3-sphere. 

\subsubsection{The Case of \texorpdfstring{$\Gtwosplit$}{G2'}} 

The group $\Gtwosplit$ can be viewed as $\Gtwosplit = \Aut(\R^{3,4}, \times_{3,4})$, where $\times_{3,4}: \R^{3,4} \times \R^{3,4} \rightarrow \R^{3,4}$ is the \emph{cross-product} on $\R^{3,4}$.

Denote a choice of simple roots for $\Gtwosplit$ by $\{\alpha_1, \alpha_2\}$, with $\alpha_1$ the short root. The corresponding partial flag manifolds are identified with the following model spaces: 
\begin{align}
    \Gtwosplit/P_{\alpha_1} \cong \Ein^{2,3} &= \{ \,[x] \in \mathbb{P}(\R^{3,4}) \; | \; q_{3,4}(x) = 0 \},\\ 
   \Gtwosplit/P_{\alpha_2} \cong \Pho^\times  &= \{\, \omega \in \Pho(\R^{3,4}) \; |\; \omega \times_{3,4} \omega = 0\}, \\ 
    \Gtwosplit/P_{\{\alpha_1, \alpha_2\}} \cong \mathcal{F}_{1,2}^\times &= \{  \,(\ell, \omega) \in \Ein^{2,3} \times \Pho^\times | \; \ell \subset \omega \}.
\end{align}
We gather further background on $\Gtwosplit$ in Subsection \ref{Sec:G2Prelims}.

Now, let $\X$ denote the $\Gtwosplit$-Riemannian symmetric space. There is a natural identification $\X \cong \Gr_{(3,0)}^\times(\R^{3,4})$, where $\Gr_{(3,0)}^\times(\R^{3,4}) \coloneqq \{ P \in \Gr_{(3,0)}(\R^{3,4}) \; |\; P \times_{3,4} P = P \}. $
For $P \in \X$, let $H_P$ denote the \emph{pointwise} stabilizer of $P$; we have $H_P \cong_{\mathbf{Lie}} \Sp(1)$. For each choice of $P$, one gets a fibration of $\Gtwofullflags$ by $H_P$-orbits. Each such orbit is a particular interesting submanifold. 
\begin{theorem}[{Theorem \ref{Thm:G2Fibration}}]
    Let $P \in \X$. There is a principal $H_P$-bundle $\mathcal{F}_{1,2}^\times \to \Flag(P) \cong \Flag(\R^3)$ with maximally transverse fibers.
\end{theorem}

\begin{figure}[ht]
\centering
\includegraphics[width=.6\textwidth]{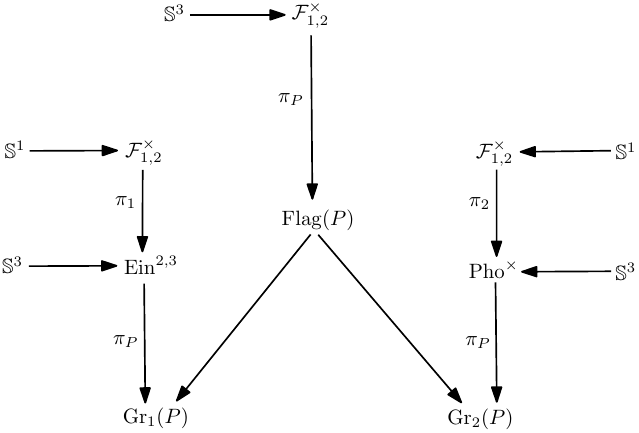}
\caption{\emph{Principal bundle fibrations of $\Gtwosplit$-flag manifolds relative to choice of $P \in \Gtwosymmetricspace$.}}
\label{Fig:GtwoFlagFibrationsSimple}
\end{figure}

The maximal transversality is easily seen from transversality. Indeed, the projection of these $\mathbb{S}^3$-fibers to $\Ein^{2,3}$ are each copies of $\Ein^{0,3}$, which are maximally transverse. 
The diagram in Figure \ref{Fig:GtwoFlagFibrationsSimple} encodes the fibrations of $\Ein^{2,3}, \Pho^\times,$ and $\Gtwofullflags$ relative to the choice of $P \in \Gtwosymmetricspace$. The precise details of the fibrations can be found in Theorem \ref{Thm:G2Fibration}. 

\vspace{-2ex}
\subsubsection{The case of \texorpdfstring{$\SO_0(3,4)$}{SO0(3,4)}}

We first view $\R^{3,4}$ as $\R^{3,0} \oplus \Ha^{0,1}$, where $\Ha^{k,l}$ denotes $(\Ha^{k+l},g_{k,l})$, where $g_{k,l}$ is a quaternion-hermitian form of signature $(k,l)$. 
We let the group $\Sp(1) \cong \SU(2) \cong \mathbb{S}^3$ include into $\SO_0(3,4)$ by declaring $ x\in Q_+(\Ha) \cong \Sp(1) $ to act trivially on $\R^{3,0}$ and to act on $\Ha^{0,1}$ by left multiplication $L_{x}$. 
Recall also that the full flag manifold $\mathcal{F}_{\Delta}$ for $G=\SO_0(3,4)$ is $\Iso_{\dbrack{3}}(\R^{3,4})$. 
We prove that $\Sp(1)$ acts freely on $\Iso_{\dbrack{3}}(\R^{3,4})$, with each orbit a maximally transverse 3-sphere.
\begin{theorem}\label{Thm:IntroB3Fibration}
Let $F \in \Iso_{\dbrack{3}}(\R^{3,4})$. 
Then the orbit $\Sp(1) \cdot F$ is a maximally transverse 3-sphere. 
\end{theorem}
Assembling these orbits, we obtain a principal bundle $\Iso_{\dbrack{3}}(\R^{3,4}) \rightarrow \Iso_{\dbrack{3}}(\R^{3,4})/\Sp(1)$ with fibers that are each maximally transverse 3-spheres. 
See Section \ref{Sec:B3Fibration} for more details. 

\subsubsection{The case of \texorpdfstring{$\SO_0(4,4)$}{SO0(4,4)}}

The full flag manifold $\Iso_{\dbrack{3}}(\R^{4,4})$ is very nearly the direct product of four maximally transverse 3-spheres.

To prove the claim, we recall the natural action of $H \coloneqq \Sp(1)^4$ on $\Iso_{\dbrack{3}}(\R^{4,4})$. 
We consider four $\Sp(1)$-subgroups of $\SO_0(4,4)$, denoted $H_L^{\pm}, H_{R}^{\pm}$, so that $H := H_L^+ \times H_L^-\times H_R^{+} \times H_R^-$. First, decompose $\R^{4,4} = \Ha^{1,0} \oplus \Ha^{0,1}$. The subgroup $K <\SO_0(4,4)$ stabilizing this splitting is $K \cong \SO(4) \times \SO(4)$, a maximal compact subgroup of $G$. Recall the standard isomorphism $Q_+(\Ha) \cong \Sp(1)$ by $x \mapsto L_x$. 
Finally, define the map $\varphi: Q_+(\Ha)^4 \rightarrow K<\SO_0(4,4)$ by 
\[ \varphi(u,v,w,x) = L_u\circ R_{v^{-1}} \oplus  L_w\circ R_{x^{-1}}.\]
The map $\varphi$ is a 4-1 cover of $\Sp(1)^4$ onto $K$. Indeed, the 2-1 homomorphism $Q_+(\Ha) \times Q_+(\Ha) \rightarrow \SO(4)$ by $(x,y) \mapsto L_x \circ  R_{y^{-1}}$ exhibits the isomorphism $\Spin(4) \cong \Sp(1) \times \Sp(1)$. The subgroup $H_{L}^+$ corresponds to $Q_+(\Ha) \times \{\mathrm{id}\} \times \{\mathrm{id}\}\times \{\mathrm{id}\} <H$, and similarly for the others. 
\begin{theorem}\label{Thm:IntroD4Fibration}
    Via the map $\varphi$, the group $\Sp(1)^4$ acts transitively on $\Iso_{\dbrack{3}}(\R^{4,4})$, with finite stabilizer. 
    Moreover, each $\Sp(1)$-subgroup $H_{L/R}^{\pm}$ has maximally transverse orbits. 
\end{theorem}

See Section \ref{Subsec:D4Fibration} for more details.

\subsection{Connection with Anosov Subgroups}\label{Sec:Anosov}

Let $G$ be a semisimple Lie group and $\Theta \subset \Delta$ a subset of simple restricted roots with associated flag manifold $\mathcal{F}_\Theta$. 
A $\Theta$-\emph{Anosov} subgroup $\Gamma < G$ is a discrete subgroup $\Gamma$ that is Gromov hyperbolic, with a continuous, equivariant, transverse embedding $\xi \colon \partial \Gamma \rightarrow \mathcal{F}_{\Theta}$ satisfying additional dynamical conditions; see \cite{GW12,GGKW17a} for a formal definition. 
The image of $\xi$ is called the (flag) limit set of $\Gamma$, see \cite{KLP17}. 
It is a difficult and interesting problem to characterize the subsets of the flag manifold $\mathcal{F}_\Theta$ that can arise as limit sets of Anosov subgroups.

\subsubsection{Sambarino's Question, the Inversion Involution \& Property (I)}

At the moment, it is generally unclear which hyperbolic groups can arise as $\Theta$-Anosov subgroups of $G$ for a fixed pair $(G,\Theta)$. 
When $\Theta' \subset \Theta$, every $\Theta$-Anosov subgroup is $\Theta'$-Anosov, so the most restricted condition is the $\Delta$-Anosov condition, which is also called \emph{Borel Anosov}.
The present situation is highlighted by the following question of Sambarino, which remains open for $d \in \{-1,0,1\} \bmod 8$. 

\begin{question}[Sambarino]\label{q:Sambarino}
Is every Borel Anosov subgroup of $\SL(d,\R)$ virtually isomorphic to a free or surface group?
\end{question}

Prior positive results for Question \ref{q:Sambarino} have been obtained by Canary-Tsouvalas when $d \in \{3,4\}$ \cite{CT20}, Tsouvalas when $d \equiv 2 \bmod 4$ \cite{tsouvalas2020borel}, and by Dey when $d \in \{2,3,4,5,6\} \bmod 8$ \cite{Dey25}. 

Sambarino's question \ref{q:Sambarino} can be asked more generally for any fixed pair $(G,\Theta)$. 
Precisely: suppose $G$ is a noncompact semisimple real Lie group. Is every $\Theta$-Anosov subgroup of $G$ virtually isomorphic to a free or surface group?
For some pairs $(G,\Theta)$ this general form is uninteresting due to easy negative examples, but in general positive answers are interesting. 
Any positive answer to Question \ref{q:circles maximal?} implies a positive answer to this ``generalized Sambarino's question," and such results are obtained in \cite{tsouvalas2020borel,DGR24,PT24,KT24}.

Dey's approach \cite{Dey25} to Sambarino's question elegantly relates to Question \ref{q:circles maximal?}. 
To make this connection precise, we recall the \emph{inversion involution}, an important ingredient in Dey's argument.
Fix a pair of transverse flags $x,z \in \mathcal{F}_\Theta$. 
Define $\mathscr{C}_x =\{y \in \mathcal{F}_{\Theta} \mid y \pitchfork x\}$ as the open Schubert cell of points transverse to $x$. 
For any point $z \in \mathscr{C}_x$, the orbit map $u \mapsto u\cdot z$ identifies the unipotent radical $U_x$ of the stabilizer subgroup $\Stab_{G}(x)$ with $\mathscr{C}_x$. 
Set $\Omega\coloneqq \mathscr{C}_x \cap \mathscr{C}_z$ as the open subset of flags $y$ transverse to both $x$ and $z$. 
It is straightforward to see that the inversion map $\iota \colon U_x \to U_x$, given by $\iota(u) = u^{-1}$, preserves $\Omega$. 

\begin{definition}[{\cite[Definition 2.3]{DGR24}}]
    A self-opposite flag manifold $\mathcal{F}$ has Property (I) if the inversion map $\iota \colon \Omega \rightarrow \Omega$ fixes no connected component. 
\end{definition}

The relationship between Property (I) and Question \ref{q:circles maximal?} is explained by Theorem \ref{Thm:IntroThmEquivalence} below, which 
answers the question raised in \cite[Remark 1.1]{Dey25} as well as \cite[Question 2.4]{DGR24} for split real groups.
Our results demonstrate that Property (I) fails to provide any further restrictions in the case of split real groups $G$, except for $E_7$. 
\begin{theorem}\label{Thm:IntroThmEquivalence}
    Let $\mathcal{F}_\Theta$ be a self-opposite flag manifold of a split real Lie group $G$.
    The following are equivalent:
    \begin{enumerate}[label=(\alph*)]
        \item $\mathcal{F}_\Theta$ has Property (I).
        \item Every transverse circle in $\mathcal{F}_\Theta$ is locally maximally transverse.
        \item $\mathcal{F}_\Theta$ does not contain a transverse $2$-sphere.
    \end{enumerate}
\end{theorem}

\begin{proof}
    The implication (a) $\implies$ (b) is \cite[Lemma 2.7]{DGR24}, based on arguments of \cite{Dey25}.  
    The implication (b) $\implies$ (c)  is immediate. 
    For the implication (c) $\implies$ (a), first suppose that $G$ is simple. 
    In our main theorems we determine, case-by-case, precisely when $\mathcal{F}_\Theta$ admits a transverse $2$-sphere. 
    Here we have applied results of \cite{tsouvalas2020borel, Dey25, DGR24, KT24} to verify that Property (I) holds in the cases complementing our construction(s). 
    This resolves the implication (c) $\implies$ (a) when $G$ is simple. 
    When $G$ is semi-simple, its associated flag manifold $\mathcal{F}_\Theta(G)$ is a product of flag manifolds $\mathcal{F}_{\Theta_i}(G_i)$ with each $G_i$ simple split.
    If any one of $\mathcal{F}_{\Theta_i}(G_i)$ has Property (I), then so does $\mathcal{F}_\Theta(G)$ \cite[Proposition 2.5]{DGR24}. 
    Otherwise, each factor admits a transverse $2$-sphere $f_i \colon \mathbb{S}^2 \to \mathcal{F}_{\Theta_i}(G_i)$ and $f = (f_i) \colon \mathbb{S}^2 \to \mathcal{F}_\Theta(G)$ is a transverse $2$-sphere. 
\end{proof}

Property (I) had not been determined for many partial flag manifolds of $\SL(d,\R)$. 
For  example, we show that Property (I) cannot hold for $\Flag(\R^{8k})$, which was left open by Dey, see \cite[Remark 2.5]{Dey25}.
That Property (I) must fail for $\Gtwosplit$ is clear from Marsh-Rietsch \cite{MR02}: they prove that the intersection of big open Schubert cells $ \Omega = \mathscr{C}_F \cap \mathscr{C}_{F'}$ has 11 connected components for any given pair of transverse flags $F \pitchfork F' \in \Flag(\Gtwosplit)$. 
Hence any involution on the components $\pi_0(\Omega)$ has a fixed point. 
For the remaining (exceptional) cases, the number of components of the corresponding double Schubert cell is even \cite{GSV03,Zel00} and Property (I) was open.
Property (I) was determined for partial flag manifolds of $\Sp(2n,\R)$ in \cite{DGR24} and for partial flag manifolds of $\SO(p,q)$ in \cite{KT24}.

\subsubsection{Obstructions} 

It is natural to ask which closed transverse subsets of flag manifolds may be realized as limit sets of Anosov subgroups.
We show in Section \ref{Sec:Obstructions} that the transverse spheres arising from the basic form of our construction do not arise as Anosov limit sets. 
Actually, we prove something stronger: in Corollary \ref{Cor:ObstructCircles}, we show that any subset $\Lambda'$ of the transverse spheres $\Lambda$ of the form \eqref{SpinorSpheresBasic} cannot be a limit set of a Borel Anosov subgroup $\Gamma$ when $\Lambda'$ contains a circle.
The key ingredient is a naturally associated filtration of the underlying vector space $\R^d$ which would necessarily be preserved by the subgroup $\Gamma$. 
The semisimplification of $\Gamma$ then preserves a nontrivial decomposition of $\R^d$ which can be ruled out by reducing to results of Canary-Tsouvalas \cite{CT20}.

These arguments do not apply to the deformed spheres $\Lambda'$ from Theorem \ref{thm:introthmfullsphere}. 
Indeed, any hypothetical Borel Anosov subgroup $\Gamma < \SL(8k,\R)$ with limit set a $3$-sphere $\Lambda'$ from Theorem \ref{thm:introthmfullsphere} would have to be irreducible.
We currently do not see an obstruction to such a potential representation, so we ask:

\begin{question}\label{q:borel anosov in G2,SO34,SO44}
    Let $G \in \{\Gtwosplit,\SO(3,4),\SO(4,4)\}$.
    Does there exist a Borel Anosov subgroup $\Gamma < G$ with $\partial \Gamma \cong \mathbb{S}^3$?
\end{question}

It is intriguing to consider the possibility of a positive answer to Question \ref{q:borel anosov in G2,SO34,SO44}. 
For concreteness, take $G=\SO(4,4)$.
The limit set of a Borel Anosov subgroup $\Gamma<\SO(4,4)$ with $3$-sphere boundary would naturally project to a transverse $3$-sphere in $\Ein^{3,3}$, and so would be $\mathbb{H}^{4,3}$-convex cocompact in the sense of Danciger-Gueritaud-Kassel \cite{DGK18}.
Beyrer-Kassel \cite{BK25} have recently shown that such an inclusion $\Gamma \hookrightarrow G$ lies in a component of the character variety consisting entirely of projective Anosov representations (in their terminology, a higher higher Teichmuller space). 

On the other hand, it seems natural to conjecture that the answer to Question \ref{q:borel anosov in G2,SO34,SO44} is no, since it would provide an incredibly strong answer to (the generalized form of) Sambarino's question.
The same discussion applies to the case of hypothetical Borel Anosov subgroups of $\SO(7,8)$ or $\SO(8,8)$ with $7$-sphere boundary.

\subsection{Structure of the paper} 
 
In Section \ref{Sec:Prelims}, we recall preliminaries on parabolic subgroups, flag manifolds, and transversality-preserving maps between flag manifolds. 

Across Sections \ref{Sec:Spinors}, \ref{Sec:G2Transversality}, \ref{Sec:TransverseSpheresDivisionAlgebras}, we build transverse spheres with spinors and explore the consequences. 
In Section \ref{Sec:Spinors}, we prove Theorem \ref{Thm:IntroArbitrarilyLargeSpheres}, which constructs arbitrarily large transverse spheres in certain type $B, D$ flag manifolds using spinors and Clifford algebras. 
We then construct deformations of the spinor spheres that remain transverse. 
In Section \ref{Sec:G2Transversality}, we consider the exceptional split real Lie group $\Gtwosplit$ of type $G_2$ and produce the fibrations of $\Flag(\Gtwosplit)$ by maximally transverse 3-spheres and discuss the corollaries.
The transverse 3-spheres in $\Flag(\Gtwosplit)$ are, in fact, built by the same recipe as in Section \ref{Sec:Spinors}. 
In Section \ref{Sec:TransverseSpheresDivisionAlgebras}, we construct transverse spheres in some low rank type $B$, $D$, cases using \emph{normed division algebras} ($\R, \C,\Ha, \Oct)$ and their split counterparts $(\C', \Ha', \Oct'$), generalizing the construction from Section \ref{Sec:G2Transversality}. 
This construction is also a special case of the construction from Section \ref{Sec:Spinors}.  
In the $B_3, D_4$ cases, we introduce fibrations of the full flag manifolds by maximally transverse 3-spheres. 
We then prove Theorem  \ref{thm:introMaxlTransversalityEveryFlag}.  

In Section \ref{Sec:MaximalTransversality}, we prove Theorem \ref{thm:introthmmaximality}, that the higher-dimensional transverse spheres from Section \ref{Sec:Spinors}, as well as their deformations, are  maximally transverse in the $B_{d-1}$ and $D_d$ cases when $d=2^k$ is a power of two, depending on $k \bmod 8$, via the Atiyah-Bott-Shapiro isomorphism \cite{ABS64}. 
We also verify maximal transversality of inclusions of these spheres to other full flag manifolds. 

The next portion of the paper consists of bootstrapping off the transverse spheres constructed in Sections \ref{Sec:Spinors}, \ref{Sec:G2Transversality},  \ref{Sec:TransverseSpheresDivisionAlgebras}. 
In Section \ref{Sec:IsotropicFlags}, we prove Theorems \ref{Thm:MainTheoremTypeB} and \ref{Thm:MainTheoremTypeD} by using a direct sum construction to build transverse subsets of flag manifolds, specifically tailored to the case of isotropic flags. 
In Section \ref{Sec:AnTransversality}, we prove Theorem \ref{Thm:MainTheoremTypeA} on the heels of a more general direct sum construction to build transverse subsets of self-opposite $\SL(d,\R)$-flag manifolds. 
In Section \ref{sec:exceptional}, we handle the remaining exceptional cases, namely of type $F_4, E_6, E_7, E_8$. 

In Section \ref{Sec:Obstructions}, we demonstrate that some of the transverse spheres we consider cannot be realized as limit sets of Borel Anosov subgroups. 
These include the basic construction of transverse spheres in full flag manifolds of classical split groups. 
We do the same for any subset of these spheres that contains a circle.

\subsection*{Acknowledgments} 
We are grateful to Dick Canary, Subhadip Dey, Zack Greenberg, Fanny Kassel, Clarence Kineider, Yosuke Morita, Alex Nolte, Jesus Sanchez, Kostas Tsouvalas, and Anna Wienhard for their encouragement and interesting conversations related to this paper. We additionally thank Alex Nolte for useful comments on a preliminary draft. 
The first author thanks Jesus Sanchez especially for many insightful discussions on spinors, comments on Lemma \ref{Lem:PinAction} (b), and for explaining the work of \cite{ABS64}. 

Part of this research was completed during the research program `Higher rank geometric structures' at the Institut Henri Poincar\'e. The authors thank the organizers for the invitation and acknowledge support of the Institut (UAR 839 CNRS-Sorbonne Universit\'e), and LabEx CARMIN (ANR-10-LABX-59-01).

\section{Preliminaries} \label{Sec:Prelims}

In this section, we review definitions of parabolic subgroups/subalgebras, flag manifolds, and transversality in (self-opposite) flag manifolds. We then review some standard geometric models for the flag manifolds of classical groups of type $A,B,D$ that we shall use throughout the paper.  

\subsection{Flag Manifolds and Transversality}\label{Sec:FlagPreliminaries}

Let $G$ be a connected real Lie group with finite center and semisimple Lie algebra $\mathfrak{g}$.
We fix a Cartan decomposition 
    \[ \mathfrak{g} = \mathfrak{k} \oplus \mathfrak{p} \]
and a maximal abelian subspace $\mathfrak{a} \subset \mathfrak{p}$.
Letting $\Sigma$ denote the set of (restricted) roots of $(\g, \mathfrak{a})$ we have the \emph{(restricted) root space decomposition}
    \[ \mathfrak{g} = \mathfrak{g}_0 \oplus \bigoplus_{\alpha \in \Sigma} \mathfrak{g}_\alpha .\]
We choose a set of positive roots $\Sigma^+$ and denote the corresponding simple roots by $\Delta$.
Then the positive Weyl chamber is given by
    \[ \mathfrak{a}^+ = \{ A \in \mathfrak{a} : \alpha(A)>0, \forall \alpha \in \Delta \}.\]
The (restricted) Weyl group $W$ is generated by reflections in the walls $\ker \alpha$ with $\alpha \in \Delta$.
There is a unique element $w_0$ in $W$ taking $\mathfrak{a}^+$ to $-\mathfrak{a}^+$. 
The map $\iota=-w_0: \mathfrak{a} \rightarrow \mathfrak{a}$ is called the \emph{opposition involution}. 
It is nontrivial precisely when the root system is of type $A_n, D_{2n+1}$ or $E_6$. The opposition involution $\iota$ induces an involution of $\Delta$ denoted by the same symbol. 

We fix a non-empty subset $\Theta \subseteq \Delta$ of simple roots. 
Define the $\Theta$-height of a positive root $\sigma \in \Sigma^+$ by writing $\sigma = \sum_{\delta \in \Delta} c^{\delta} \delta$ and setting $\mathrm{ht}_{\Theta}(\sigma) \coloneqq \sum_{\delta \in \Theta} c^{\delta}$.
Associated to $\Theta$ are the following nilpotent subalgebras: 

\[ \mathfrak{u}_\Theta = \bigoplus_{\mathrm{ht}_{\Theta}(\sigma) >0} \mathfrak{g}_{\sigma}, \quad \mathfrak{u}_\Theta^{opp} = \bigoplus_{\mathrm{ht}_{\Theta}(\sigma)< 0} \mathfrak{g}_{\sigma}.\]
The normalizer of $\mathfrak{u}_\Theta$ (resp.\ $\mathfrak{u}^{opp}_\Theta$) with respect to the adjoint action is the \emph{parabolic subgroup} $P_\Theta$ (resp.\ $P_\Theta^{opp}$).
We note that $P_{\iota \Theta}$ is conjugate to $P_\Theta^{opp}$.
The \emph{$\Theta$-flag manifold} is the homogeneous space 
    \[ \mathcal{F}_\Theta = G/P_\Theta.\]
We may write $\Flag(G)$ to refer to the full flag manifold $\mathcal{F}_{\Delta}$, especially when the group $G$ is changing frequently. \medskip 

Pairs in $G/P_\Theta \times G/P_\Theta^{opp}$ of the form $(gP_\Theta,gP_\Theta^{opp})$ are called \emph{transverse}. 
The subset of $G/P_\Theta \times G/P_\Theta^{opp}$ consisting of transverse pairs is the unique open and dense $G$-orbit.
A flag manifold $\mathcal{F}_{\Theta}$ is called \emph{self-opposite} if the subgroups $P_{\Theta}$ and $P_{\Theta}^{opp}$ are conjugate under $G$ to one another, or, equivalently, when $\Theta$ is preserved by the opposition involution. 
If $\mathcal{F}_{\Theta}$ is self-opposite, then $G/P_{\Theta}$ and $G/P_{\Theta}^{opp}$ are canonically identified, so we may speak of whether $(x,y) \in \mathcal{F}_{\Theta}^2$ are transverse. 
A subset $S \subset \mathcal{F}_{\Theta}$ of a self-opposite flag manifold is \emph{transverse} if every pair of distinct points $x \neq y \in S$ is transverse, \emph{maximally transverse} if $S$ is not properly contained in another transverse subset of $\mathcal{F}_{\Theta}$, and \emph{locally maximally transverse} if there exists a neighborhood $U$ of $S$ so that $S$ is not properly contained in a transverse subset contained in $U$.  

An element $A$ in $\overline{\mathfrak{a}^+}$ is called \emph{$\Theta$-regular} if $\alpha(A) >0$ for all $\alpha \in \Theta$. 
For such an element $A$, there is an associated parabolic subgroup $P_A$ satisfying   
    \[ P_A \coloneqq \left\{ g \in G : \lim_{t\to\infty} \exp(-tA)g\exp(tA) \text{ exists} \right\} \subseteq P_\Theta,\]
with equality exactly when each $\beta \in \Delta  \setminus \Theta$ is zero on $A$, see e.g.\ \cite{eberlein96book}.

Now let $G'$ be another semisimple Lie group, with similar choices, and a subset $\Theta'$ of $\Delta'$.
A map of flag manifolds is called \emph{transversality-preserving} if it takes transverse pairs to transverse pairs. For any set $X$, a map $f:X \rightarrow \mathcal{F}_{\Theta}$ is called \emph{transverse} if it takes distinct pairs to transverse pairs. Hence,  $S\subset \mathcal{F}_{\Theta}$ is transverse exactly when the inclusion $ S \hookrightarrow \mathcal{F}_{\Theta}$ is a transverse map. 

The following criterion will be frequently used to produce transversality-preserving maps. 
\begin{lemma}[Transversality-Preserving Maps of Flag Manifolds]\label{lem:transverse map criterion}
    Suppose $G \to G'$ induces a map $\mathfrak{a} \to \mathfrak{a}'$ taking $A \in \mathfrak{a}$ to $A' \in \mathfrak{a}'$.
    If $A'$ is $\Theta'$-regular, then there is an induced map $\mathcal{F}_A \to \mathcal{F}_{\Theta'}$.
    Let $\iota$ (resp.\ $\iota'$) be the opposition involution of a chamber containing $A$ (resp.\ $A'$).
    If moreover $A=\iota(A)$, $\Theta'$ is self-opposite, and $\iota'(A')=A'$, then the map $\mathcal{F}_A \to \mathcal{F}_{\Theta'}$ is transversality-preserving.
\end{lemma}

\begin{proof}
    Since $G \to G'$ sends $P_A$ to $P_{A'}$ the map $gP_A \mapsto gP_{\Theta'}$ is well-defined.    
    Transversality follows when $w_0(A)=-\iota(A) \mapsto -\iota'(A)=w_0'(A)$.
\end{proof}

We note some non-examples of transversality-preserving maps. 

\begin{example}\label{Ex:D3B3}
    Consider $\SO_0(3,3) \hookrightarrow \SO_0(3,4)$.
    There is an equivariant map of full flag manifolds, but no such transversality-preserving map. More generally, this holds for $D_{2n-1}$ and $D_{2n}$, for $n \geq 2$.
\end{example}

\begin{example}\label{Ex:AevenAodd}
Consider $\SL(3,\R) \hookrightarrow \SL(4,\R)$. 
As in the previous example, there is an equivariant map $\Flag(\R^{3}) \rightarrow \Flag(\R^4)$, but no such transversality-preserving map. The same holds for $A_{2n}$ and $A_{2n+1}$ for $n \ge 1$. 
\end{example}

We shall often make use of the following fact.

\begin{fact}[Projections of Flag Manifolds]\label{fact:TransverseProjection}
Suppose that $\Theta' \subset \Theta$. 
The natural projection $\mathcal{F}_{\Theta}\rightarrow \mathcal{F}_{\Theta'}$ is transversality-preserving.
\end{fact}

It is an observation due to Dey that a transverse subset $S \subset \mathcal{F}_{\Theta}$ that is non-nullhomotopic must be maximally transverse. 
We shall use this fact often. 
\begin{fact}\label{Fact:MaximallyTransverse}
Suppose that $S \subset \mathcal{F}_{\Theta}$ is transverse and non-nullhomotopic in $\mathcal{F}_{\Theta}$. 
Then $S$ is maximally transverse in $\mathcal{F}_{\Theta}$. 
\end{fact}
Fact \ref{Fact:MaximallyTransverse} holds due to the following: if $S \subset \mathcal{F}_{\Theta}$ is not maximally transverse, then we may select a flag $x \in \mathcal{F}_{\Theta}$ transverse to all of $S$. 
Hence, $S$ is contained in the open Schubert cell $\mathscr{C}_x = \{F \in \mathcal{F}_{\Theta} \mid x \pitchfork F\}$ and thus $S$ is contractible in $\mathcal{F}_{\Theta}$. 

\subsection{Geometric Models for Flag Manifolds of Type \texorpdfstring{$A,B,D$}{A,B,D}}\label{Subsec:ConcreteModelsABD}

We recall well-known geometric models for the partial flag manifolds of $G$ of type $A,B,D$. 
We then discuss transverse embeddings between some type $A,B,D$ full and partial flag manifolds. 

\subsubsection{Type \texorpdfstring{$A$}{A} Case.}\label{Subsec:TypeAFlags}

\begin{figure}[h]
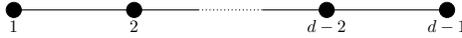

		\centering
		\begin{equation*}
			\begin{dynkinDiagram}[text style/.style={scale=.6},
				edge length=.8cm,
				scale=2,
				labels={1,2,d-2,d-1},
				label macro/.code={{{#1}}}
				]A{}
			\end{dynkinDiagram}
		\end{equation*}
		\caption{\emph{\small{Dynkin diagram of type $A_{d-1}$.}}}
		\label{fig:DynkinAp}
	\end{figure}

We consider type $A_{d-1}$.
We start by recalling standard Lie theory for $\g = \mathfrak{sl}_d(\R)$. 
Choose any basis $(v_i)_{i=1}^d$ for $\R^d$. 
A choice of maximal split torus $\mathfrak{a} < \mathfrak{sl}_d(\R)$ is then given by
\[ \mathfrak{a} = \left \{ \mathrm{diag}(\lambda_1, \dots, \lambda_d) \in \mathfrak{gl}_d(\R) \; \bigg| \; \sum_{i=1}^d \lambda_ i= 0 \right \}.\]
A choice of simple roots for $\Sigma(\g, \mathfrak{a})$ is given by $\alpha_i \coloneqq \lambda_i^*-\lambda_{i+1}^{*}$ for $1 \leq i \leq n-1$. 
For ease, we may refer to the root $\alpha_i$ by its index $i$. 
Under this labeling, the associated partial flag manifold $G/P_{\alpha_i}$ is naturally identified with  $\Gr_i(\R^d)$. 
This correspondence aligns with the $A_{d-1}$-Dynkin diagram, as in Figure \ref{fig:DynkinAp}. 

We now set a notational convention going forwards. Define 
\[ \dbrack{k} \coloneqq \{ 0,1,2,\dots, k\}. \]
The simple roots $\Delta$ of $\Sigma(\g,\mathfrak{a})$ naturally identify with $\{1,2, \dots, d-1\}$. 
On the other hand, it will be convenient for us to imagine $\Theta \subset \dbrack{d}$ instead, where $\Theta = \{0,i_1, \dots, i_k, d\}$ with $i_{j} < i_{j+1}$. 
With our indexing, $\mathcal{F}_{\Theta} \cong \Flag_{\Theta}(\R^d)$ is the collection of nested $\Theta$-index subspaces: 
\[ \Flag_{\Theta}(\R^d) \coloneqq \{ \;(0,V^{i_1},\, \dots, \,V^{i_k},\R^d) \; |\; V^{i_j} \in \Gr_{i_j}(\R^n), \,  V^{i_j} \subset V^{i_{j+1}}\}.\]
Under this convention, $\Flag_{\dbrack{d}}(\R^d) = \Flag(\R^d)$ is the space of full flags in $\R^d$. 
When convenient, we may omit the trivial 0-subspace $\{0\}$ and the $d$-subspace $\R^d$, but in Section \ref{Sec:AnTransversality}, we shall use these two subspaces.

A flag manifold $\mathcal{F}_{\Theta}$ is self-opposite if and only if $\Theta$ is symmetric in the sense that $\iota(\Theta) = \Theta$, for $\iota \colon \dbrack{d} \rightarrow \dbrack{d}$ given by $\iota(k) = d-k$. 
In the case $\Theta$ is symmetric, then a pair of flags $(V^{\bullet}, W^{\bullet}) \in \Flag_{\Theta}(\R^d)^2$ are \emph{transverse} if and only if  $V^k + W^{d-k} = \R^d$ for all $k \in \Theta$. 

\subsubsection{Type \texorpdfstring{$B$}{B} Case.}\label{Subsec:TypeBFlags}

We now consider type $B_p$ and begin by recalling standard Lie theory of $\mathfrak{so}(p,p+1)$.
Let us consider a basis $(x_i)_{i=1}^{2p+1}$ for $\R^{p,p+1}$ such that the quadratic form $q_{p,p+1}$ is represented by the following matrix $[q_{p,p+1}]$:
\[ [q_{p,p+1} ] = \begin{pmatrix} 
                    & & & & & &1 \\ 
                    & & & & &\reflectbox{$\ddots$}&  \\
                    & & & & 1& &  \\
                    & & &-1 & & & \\     
                    & & 1& & & &\\ 
                    &\reflectbox{$\ddots$}& & & & & \\   
                    1 & & & & & &\\
\end{pmatrix} .\] 
In such a basis, the diagonal transformations in $\g=\mathfrak{so}(p,p+1)$ form a maximal split torus $\mathfrak{a}$ given by 
\[ \mathfrak{a} = \{ \mathrm{diag}(\lambda_1, \lambda_2 , \dots, \lambda_p,0,-\lambda_p, \,\dots , -\lambda_2,-\lambda_1) \in \mathfrak{gl}_{2p+1}(\R)\; | \; \lambda_i \in \R   \}.\]
A choice of $p$ simple roots for $\Sigma(\g,\mathfrak{a})$ is given by $\alpha_i = \lambda_i^*-\lambda_{i+1}^*$ for $1 \leq i \leq p-1$ and $\alpha_p =\lambda_p^*$. 
With this labeling, one finds the associated flag manifolds are $\mathcal{F}_{\Theta}$ are given by $\mathcal{F}_{\Theta} \cong \Iso_{\Theta}(\R^{p,p+1})$, where $\Iso_{\Theta}(\R^{p,p+1})$ denotes totally isotropic $\Theta$-index flags that are nested. 
Write $\Theta = \{i_1, \dots, i_k\}$ in ascending order and 
\[ \Iso_{\Theta}(\R^{p,p+1}) \coloneqq \{(F^{i_1}, \dots, F^{i_k}) \; | \; F^{i_j} \in \Iso_{i_j}(\R^{p,p+1}), \, F^{i_j} \subset F^{i_{j+1}}\},\] 
where $\Iso_{k}(\R^{p,p+1})$ denotes the Grassmannian of totally isotropic $k$-planes in $\R^{p,p+1}$. 
A basic fact is that $\Iso_k(\R^{p,q}) = \emptyset$ for $k > \max \{p,q\}$. 

\begin{figure}[h]
		\centering
		\begin{equation*}
			\begin{dynkinDiagram}[text style/.style={scale=.6},
				edge length=.8cm,
				scale=2,
				labels={1,2,p-2,p-1,p},
				label macro/.code={{{#1}}}
				]B{}
			\end{dynkinDiagram}
		\end{equation*}
		\caption{\emph{\small{Dynkin diagram of type $B_p$.}}}
		\label{fig:DynkinBp}
	\end{figure}

In the $B_p$ case, the opposition involution is trivial and every flag manifold $\Iso_{\Theta}(\R^{p,p+1})$ is self-opposite. 
Two flags $V^{\bullet}, W^\bullet \in \Iso_{\Theta}(\R^{p,p+1})$ are transverse if and only if $V^k+(W^{k})^{\bot} = \R^{p,p+1}$ for all $k \in \Theta$. \medskip 

\begin{remark}
The first two isotropic Grassmannians enjoy standard alternate names and notations, in the case of any mixed signature. 
These are the ``Einstein universe" $\Ein^{p-1,q-1} \coloneqq\Iso_{1}(\R^{p,q})$, and the ``space of photons" $\Pho(\R^{p,q})\coloneqq \Iso_2(\R^{p,q})$. 
\end{remark}

\subsubsection{Type \texorpdfstring{$D$}{D} Case.}\label{Subsec:TypeDFlags}

\begin{figure}[h]
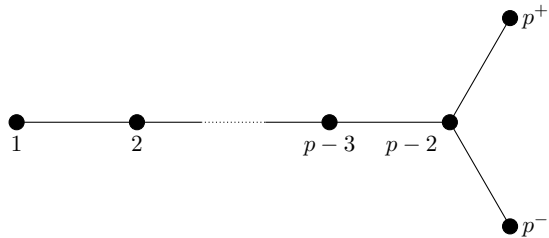

		\centering
		\begin{equation*}
			\begin{dynkinDiagram}[text style/.style={scale=.8},
				edge length=.8cm,
				scale=2,
				labels={1,2,p-3,p-2,p^+, p^-},
				label macro/.code={{{#1}}}
				]D{}
			\end{dynkinDiagram}
		\end{equation*}
		\caption{\emph{\small{Dynkin diagram of type $D_p$.}}}
		\label{fig:DynkinDp}
	\end{figure}

For the case of $D_p$-geometry, it is useful to consider a basis $(x_i)_{i=1}^{2p}$ for $\R^{p,p}$ such that the quadratic form $q_{p,p}$ obtains the anti-diagonal form 
\[ [q_{p,p} ] = \begin{pmatrix} & & 1 \\ 
                                & \reflectbox{$\ddots$} &  \\
                               1 & &  \\
\end{pmatrix} .\] 
This choice of basis is convenient as the diagonal transformations in $\g =\solie(p,p)$ form a maximal split torus $\mathfrak{a} < \solie(p,p)$ in this basis.
Indeed, one finds $\mathfrak{a}$ to be the following:
\[ \mathfrak{a} = \{  \mathrm{diag}(\lambda_1, \lambda_2 , \dots, \lambda_p,-\lambda_p, \,\dots , -\lambda_2,-\lambda_1) \in \mathfrak{gl}_{2p}(\R)\; | \; \lambda_i \in \R \}.\]
Next, label the simple roots of $\Sigma(\g,\mathfrak{a})$  by $\Delta = \{ 1, 2, \dots, p-2\} \cup \{p^+, p^-\}$, according to the Dynkin diagram, see Figure \ref{fig:DynkinDp}.
For $1 \leq i \leq p-2$, the $i^{th}$ root is $\lambda_i^*-\lambda_{i+1}^*$. 
The final two exceptional roots are as follows. 
The root we call $p^+$ is $\lambda_{p-1}^*-\lambda_p^*$ and the root we call $p^-$ is $\lambda_{p-1}^*+\lambda_{p}^*$. 
Under this labeling, there are $\SO_0(p,p)$-equivariant identifications $\mathcal{F}_{\{i\}} \cong \Iso_{i}(\R^{p,p})$ for $1 \leq i \leq p-2$. 
This corresponds to the identification, for $1\leq i \leq p-2$, given by $P_{\alpha_i} = \Stab_{\SO_0(p,p)}(V^i)$, where $V^i = \spann_{\R} \{ x_1,\, x_2 , \,\dots, \,x_i \}$.

For the final two roots, things are a bit more subtle, as the group $\SO_0(p,p)$ does not act transitively on $\Iso_p(\R^{p,p})$. 
One can distinguish the two $\SO_0(p,p)$ orbits in $\Iso_p(\R^{p,p})$ as follows. Fix a maximal spacelike $p$-plane $P \subset \R^{p,p}$ to write $\R^{p,p} = P\oplus P^\bot$. Then every isotropic $p$-plane $T$ obtains the form $T = \graph(\varphi_{T,P})$ for some anti-isometry $\varphi_{T,P}: P \rightarrow P^\bot$. Recall that $\SO_0(p,p)$ preserves both space and time-orientation. For any choice of oriented bases of $P, P^\bot$, the orbit invariant of $T$ is $\det(A) \in \{-1,+1\}$, with $A \in \Mat_p(\R)$ representing $\varphi_{T,P}$ in the chosen bases. Equivalently, we ask if $\varphi_{T,P}$ takes a space-oriented basis of $P$ to a time-oriented basis of $P^\bot$.\footnote{Of course, this discussion requires a background choice of space and time-orientation on $\R^{p,p}$, and changing one (space/time) orientation but not the other will change $\det(A)$.}

We now introduce some more notation for these orbits. Let $ \Iso_{p}^{\pm}(\R^{p,p})$ denote a labeling of the two $\SO_0(p,p)$-orbits such that  
$\mathcal{F}_{\{p^{\pm}\}} \cong \Iso_{p}^{\pm}(\R^{p,p})$. With this notation, 
\begin{align}
    \begin{cases}
        P_{p^+} = \Stab_{\SO_0(p,p)}(T_+), \mathrm{where} \; T_+ = T_{p-1} \oplus \R \{x_{p} \}  \;\\
        P_{p^-} = \Stab_{\SO_0(p,p)}(T_-),\mathrm{where} \; T_- = T_{p-1} \oplus \R \{ x_{p+1} \} .
    \end{cases}
\end{align}
The two flag manifolds $\Iso_{p}^{\pm}(\R^{p,p})$ are self-opposite if and only if $p$ is even. Indeed, when $p$ is odd, then $\iota(p^+) = p^-$. Note that $\iota(j) = j$ for $1\leq j \leq p-2$ regardless of parity of $p$. 

When convenient, we may avoid studying $\Iso_{p}^{\pm}(\R^{p,p})$ altogether when considering $\SO_0(p,p)$-full flags, by using the identification $\mathcal{F}_{\{p^+, p^-\} } \cong \Iso_{p-1}(\R^{p,p})$. 
Geometrically, the model for $\mathcal{F}_{\{p^+, p^-\}}$ is natural for the following reason. 
\begin{proposition}\label{Prop:TwoParents}
Given any sub-maximal isotropic plane $T \in \Iso_{p-1}(\R^{p,p})$, the space $T$ is contained in exactly two isotropic $p$-planes $T_{\pm}$, where $T_{\pm} \in \Iso_p^{\pm}(\R^{p,p})$. Denote $\nu_{\pm}:\Iso_{p-1}(\R^{p,p}) \rightarrow \Iso_{p}^{\pm}(\R^{p,p})$ as the $\SO_0(p,p)$-equivariant maps $\nu_{\pm}(T) = T_{\pm}$. 
\end{proposition}
\begin{proof}
Observe that $T^\bot $ has signature $(1,1, p-1)$, with the signature indicating the positive, negative, and isotropic parts, respectively. Thus, $\Ein(T^\bot) \cong \mathbb{S}^0$. By an appropriate choice of labels  $\ell_{\pm} \in \Ein(T^{\bot})$, one finds $T_{\pm} = T \oplus \ell_{\pm}$. 
\end{proof} 
In this way, one can imagine $\mathcal{F}_{\{p^+,p^-\}} \hookrightarrow \mathcal{F}_{p^+} \times \mathcal{F}_{p^-}$. As in the type $A,B$ cases, the type $D$ flag manifolds mirror the shape of the corresponding Dynkin diagram; the wishbone shape shown in Figure \ref{fig:DynkinDp} is the cause of the exceptions with the simple roots $p^{\pm}$.

With these remarks in place, we introduce some notation for $\SO_0(p,p)$-full flags. We maintain the notation $\dbrack{k} \coloneqq \{0,1,2,\dots, k\}$ from before. 

We will often write $\Flag(\R^{p,p})$ for the first naive geometric model of $\SO_0(p,p)$-full flags $(F^1, \dots, F^{p-2}, F^p_+, F^p_-)$ that are nested appropriately. 
As seen in Proposition \ref{Prop:TwoParents}, this model is redundant. 
Next, we denote
\[\Flag_+(\R^{p,p}) \coloneqq \left\{ F = (F^i)_{i=1}^{p} \in \Iso_{\dbrack{p}}(\R^{p, p}) \mid 1 \leq i \leq p-1:  F^i \in \Iso_i(\R^{p,p}), F^{i} \subset F^{i+1}, \, F^{p} \in \Iso^+_{p}(\R^{p,p}) \right\}. \] 
Finally, for $\Theta \subset \dbrack{p-1}$, we write $\Theta = \{i_1, \dots, i_j\}$ in order and then use $\Iso_{\Theta}(\R^{p,p})$ to denote  
\[\Iso_{\Theta} (\R^{p,p}) \coloneqq \{ F =(F^1, \dots, F^{j} ) \mid F^{i_\ell} \in \Iso_{i_{\ell}}(\R^{p,p}), F^{i} \subset F^{i+1}\} .\]
There are $\SO_0(p,p)$-equivariant diffeomorphisms 
\begin{align}
    \begin{cases}
        \Iso_{\dbrack{p-1}}(\R^{p,p}) \rightarrow \Flag(\R^{p,p}) \; &\mathrm{by} \;  (F^1, \dots, F^{p-1}) \mapsto \big(F^1, \dots, F^{p-2}, \nu_{+}(F^{p-1}), \nu_{-}(F^{p-1}) \big),\\
        \Iso_{\dbrack{p-1}}(\R^{p,p}) \rightarrow \Flag_+(\R^{p,p}) \; &\mathrm{by} \;  (F^1, \dots, F^{p-1}) \mapsto \big(F^1, \dots, F^{p-2}, \nu_+(F^{p-1})\, \big) .
    \end{cases}
\end{align}
When $p$ is even, we may be more liberal with the notation $\Iso_{\Theta}$. 
For $\Theta = \dbrack{p-1} \cup \{p^+\}$ we write $\Iso_{\Theta}(\R^{p,p}) =\Flag_+(\R^{p,p})$, and note that each component subspace is self-opposite.

When studying $\SO_0(p,p)$-full flags by themselves, the most convenient model is $\Iso_{\dbrack{p-1}}(\R^{p,p})$. However, in direct sum constructions of isotropic flags in Section \ref{Sec:IsotropicFlags}, we will need the top-dimensional isotropic subspace, so we shall use the model $\Flag_+(\R^{2p,2p})$ for the $\SO_0(p,p)$-full flag manifold. \medskip 

We now briefly elaborate on the embedding $\Iso_{p-1}(\R^{p,p})\hookrightarrow \Iso_{p}^+(\R^{p,p}) \times \Iso_p^-(\R^{p,p})$ and how transversality looks in these models depending on parity of $p$. \medskip 

Let $p$ be odd. Take $F_i \in \Iso_{p-1}(\R^{p,p})$ for $i \in \{1,2\}$. Then we observe the following: 
\[ F_1 \pitchfork F_2 \iff F_1 + F_2^\bot = \R^{p,p} \iff \begin{cases} F_1^+ +F_2^- = \R^{p,p}\\
                 F_1^- +F_2^+ = \R^{p,p} \end{cases}  \iff (F_1^+,F_1^-)\pitchfork (F_2^+,F_2^-). \] 

When $p$ is even, $\Iso_{p}^{\pm}(\R^{p,p})$ are self-opposite. In this case, for $F_i \in \Iso_{p-1}(\R^{p,p})$ with $i \in \{1,2\}$, 
\[F_1 \pitchfork F_2 \iff  F_1 +F_2^\bot = \R^{p,p} \iff 
 \begin{cases} F_1^+ +F_2^+ = \R^{p,p}\\
                 F_1^- +F_2^- = \R^{p,p} \end{cases} \iff (F_1^+,F_1^-) \pitchfork (F_2^+,F_2^-). \] 
In summary, regardless of parity of $p$, for $F_i \in \Iso_{p-1}(\R^{p,p})$, we have $F_1 \pitchfork F_2$ if and only if $F_1+ F_2^\bot = \R^{p,p}$, and we may ignore considerations of $F_i^{\pm}$. Moreover, using $\Iso_{\dbrack{p-1}}(\R^{p,p})$ as the model for the $\SO_0(p,p)$ full flag manifold, for $F_i = (F_i^1,\dots, F_{i}^{p-1}) \in \Iso_{\dbrack{p-1}}(\R^{p,p})$, we have $F_1 \pitchfork F_2$ if and only if $\R^{p,p} = F_{1}^{j} + (F_{2}^{j})^{\bot}$ for all indices $ 1\leq j \leq p-1$.

\subsubsection{Transversality-Preserving Maps via Geometric Models}\label{Subsec:TranverseEmbeddingsABD}

Now, we give straightforward descriptions of (equivariant) transversality-preserving maps of flag manifolds of type $A,B,D$ in the geometric models presented.  
\begin{example}\label{Ex:D2nA4n-1}
    ($D_{2n},A_{4n-1})$. There is a transversality-preserving map $\Flag_+(\R^{2n,2n}) {\hookrightarrow} \Flag(\R^{4n})$
    \[ (F^1, \dots, \,F^{2n-1},\, F^{2n}_+) \mapsto (F^1, \dots, \, F^{2n-1}, \,F^{2n}_+,\, (F^{2n-1})^\bot, \dots, \,(F^1)^\bot).\] 
\end{example}

\begin{example}\label{Ex:D2B2n}
    $(D_{2n},B_{2n})$. There is a transversality-preserving map $\Flag_+(\R^{2n, 2n}) {\hookrightarrow} \Iso_{\dbrack{2n}}(\R^{2n, 2n+1})$ 
    \[ (F^1, \dots, F^{n-1}, F^n_+) \mapsto (F^1, \dots, F^{n-1}, F^{n}_+ ). \]
    The map appears trivial, but the flags have changed homes. 
\end{example}

\begin{example}\label{Ex:BnDn+1}
    $(B_{n},D_{n+1})$. As in the previous example, the transversality-preserving map $\Iso_{\dbrack{n}}(\R^{n,n+1}){\hookrightarrow} \Iso_{\dbrack{n}}(\R^{n+1, n+1})$ looks  like the identity map:
    \[ (F^1, \dots, F^{n}) \mapsto (F^1, \dots, F^{n}) .\]
\end{example}

\begin{example}\label{Ex:BnA2n}
    $(B_{n},A_{2n})$. A transversality-preserving map $\Flag(\R^{n,n+1}) {\hookrightarrow} \Flag(\R^{2n+1})$ is given by 
    \[ (F^1, \dots, \,F^{n-1}, \, F^{n}) \mapsto  (F_1, \dots, \,F^{n-1}, \, F^n, \, (F^{n})^\bot, \,(F^{n-1})^\bot, \dots, \,(F^1)^\bot).\] 
\end{example}

\begin{example}\label{Ex:GeneralChains}
    Combining all the previous results, one finds the following chains of transversality-preserving embeddings: 
    \begin{itemize}
        \item $B_{2n-1} \hookrightarrow D_{2n} \hookrightarrow B_{2n} \hookrightarrow A_{4n}$.  
        \item  $B_{2n-1} \hookrightarrow D_{2n} \hookrightarrow A_{4n-1}.$
        \item $B_{2n-1} \hookrightarrow D_{2n} \hookrightarrow B_{2n} \hookrightarrow D_{2n+1}$.
    \end{itemize}
\end{example} 
\medskip 

\begin{example}
    $(A_{2n-1},A_{2n})$. Write $\R^{2n+1} = \R^{2n} \oplus L$. One such transversality-preserving map  
    $\Flag(\R^{2n}) {\hookrightarrow} \Flag(\R^{2n+1})$ is given by  
    \[ (F^1, \,F^2, \dots, \,F^{2n-1}) \mapsto \big(F^1, \dots, \,F^{n}, \,F^n\oplus L, \,F^{n+1} \oplus L, \dots, \,F^{2n-1} \oplus L \big).\]
\end{example}

\begin{example}
    Non-example: $(D_{2n+1},B_{2n+1})$, as an addendum to Example \ref{Ex:D3B3}. There are two  $\SO_0(2p+1,2p+1)$-equivariant maps $\psi_{\pm}: \Flag(\R^{2p+1,2p+1}) \rightarrow \Flag(\R^{2p+1, 2p+2})$ given by 
    \[ \psi_{\pm} := \bigg ( (F^1, \dots, F^{p-2}, F^{p}_+,F^p_-) \mapsto (F^1, \dots, F^{p-2}, F^p_+ \cap F^{p}_-, F^{p}_{\pm})\; \bigg). \]
    The maps $\psi_{\pm}$ are not transversality-preserving, essentially because $\Iso_{2p+1}^{\pm}(\R^{2p+1,2p+1})$ is not self-opposite. 
\end{example}

\begin{example}\label{Ex:DnA2nMinusn}
    $(D_{n},A_{2n})$. There is a transversality-preserving embedding $\Iso_{\dbrack{n-1}}(\R^{n,n}) \hookrightarrow \Flag_{\dbrack{2n}\setminus \{n\}}
        (\R^{2n})$ given by
    \[ (F_1, \dots, \,F_{n-1}) \mapsto  (F_1, \dots, \,F_{n-1}, \,F_{n-1}^\bot, \dots, \,F_1^\bot).\] 
    When $n$ is odd, we can do no better; there is no  equivariant transversality-preserving embedding $\Iso_{\dbrack{n-1}}(\R^{n,n}) \rightarrow \Flag(\R^{2n})$ in this case.
\end{example}

\section{Transverse Spheres via Spinors}\label{Sec:Spinors}

In this section, for any dimension $n$, we construct a transverse $(n-1)$-sphere in $\Flag(\R^{d(n)+\epsilon, d(n)})$ for an explicit integer $d(n) \sim 2^{n/2}$ depending on $n$ and any $\epsilon \in \{-1,0,1\}$. 
The construction is done with Clifford algebras and spinors. 
We then discuss deformations of the resulting transverse spheres. 
Maximal transversality of the spheres is discussed in Section \ref{Sec:MaximalTransversality}.

\subsection{Background on Clifford Algebras \& Spinors}

We recall the details on (real) Clifford algebras and spinors necessary for the construction of transverse spheres of arbitrarily large dimension. 
For a comprehensive account of spinors, we refer the reader to the excellent book of Lawson \& Michelsohn \cite{LM89}.
For recent work which builds \emph{explicit} models for the irreducible representations of the real algebras $\Cl(n)$, see \cite{San25}.  

\subsubsection{Definitions and Properties of Clifford Algebras} 

We recall the standard construction of real (Euclidean) Clifford algebras and some of their essential properties.

\begin{definition}\label{Defn:Clifford}
Let $(V,q) \cong \R^{n,0}$ be a Euclidean vector space over $\R$. Define the tensor algebra $T(V) = \bigoplus_{i=0}^{\infty} V^{\otimes i}$, where $V^0 = \R$. 
The \textbf{Clifford algebra} $\Cl(V)$ is the $\R$-algebra defined by $\Cl(V) \coloneqq  T(V)/ I $, where $I$ is the two-sided ideal generated by elements of the form $v\otimes v +q(v)1$, for $v \in V$.  
\end{definition} 

Alternatively, the algebra $\Cl(V)$ is reached from $T(V)$ by declaring relations of the form
\begin{align}\label{CliffordCondition}
    u \otimes v+ v\otimes u = -2q(u,v) \,1,
\end{align}
which follows by polarizing the relation $v \otimes v = -q(v) 1$. 

Going forward, we may denote the algebra product $xy$ of $x, y \in \Cl(V)$ by juxtaposition. Next, we recall some important facts about the structure of $\Cl(V)$. \medskip 

There is a canonical \emph{vector space} isomorphism $\Cl(V) \cong \Lambda^\bullet V$, where $\Lambda^\bullet V$ is the exterior algebra of $V$. This map is induced by the map of ``simple tensors''
\[ e_{i_1}e_{i_2}\cdots e_{i_n} \mapsto e_{i_1} \wedge \dots \wedge e_{i_n} .\]
In particular, $\dim_{\R}(\Cl(V))= 2^m$, where $\dim_{\R} (V)= m$. This isomorphism induces a canonical $\Z$-grading on $\Cl(V)$. We write $\Cl(V) = \bigoplus_{i\in \Z}\Cl_{i}(V)$ for this $\Z$-grading. We will seldom use this grading. 

The Clifford algebra also carries a coarse $\Z_2$-grading, which is actually an \emph{algebra} grading. Here, the even and odd graded subspaces $\Cl^0(V)$, $\Cl^1(V)$, respectively, are given by 
\begin{align}
    \begin{cases}
        \Cl^0(V) = \bigoplus_{i \in \Z_{\geq0}} \Cl_{2i}(V)\\
        \Cl^1(V) = \bigoplus_{i \in \Z_{\geq0}}\Cl_{2i+1}(V).
    \end{cases}
\end{align}
We write $\Cl(V) = \Cl^0(V) \oplus \Cl^1(V)$ for this algebra grading. In particular, $\Cl^i(V) \,\Cl^j(V) = \Cl^{i+j}(V)$ for $i, j \in \Z_2$, so that $\Cl^0(V)$ is a subalgebra. We will frequently use this grading.\medskip 

One of the most basic and useful properties $\Cl(V)$ is its \emph{universal property}.

\begin{proposition}[Universal Property of $\Cl(V)$]\label{Prop:CliffordUniversal}
If $\rho: V \rightarrow \mathcal{A}$ is a linear map to an associative unital $\R$-algebra $\mathcal{A}$ such that $\rho(v)^2 = - q(v) 1_{\mathcal{A}}$, then there is a unique homomorphism of $\R$-algebras 
$\hat{\rho}: \Cl(V) \rightarrow \mathcal{A}$ such that $\hat{\rho} \circ \iota = \rho$, for $\iota: V \hookrightarrow \Cl(V)$. The map $\rho$ is called the \textbf{defining map} of $\hat{\rho}$. 
\end{proposition} 

The proof of Proposition \ref{Prop:CliffordUniversal} is straightforward. The map $\rho$ extends uniquely to a map $T(V) \rightarrow \mathcal{A}$ in the obvious way, and this map descends to $\Cl(V)$ due to the condition $\rho(v)^2 = -q(v)1_{\mathcal{A}}$. We remark that one can \emph{define} the Clifford algebra via Proposition \ref{Prop:CliffordUniversal}; there is a unique $\R$-algebra, up to isomorphism, satisfying the universal property, and it is $\Cl(V)$ as defined in Definition \ref{Defn:Clifford}.

As a corollary to the proposition, we see the functorial nature of $\Cl(\cdot)$. 
\begin{corollary}
Let $V, V'$ be Euclidean vector spaces of the same dimension. Then any isometry $\varphi: V\rightarrow V'$ induces an isomorphism of real algebras $\Cl(V) \rightarrow \Cl(V')$. 
\end{corollary}
We now write $\Cl(n)$ to unambiguously denote $\Cl(V)$ for any Euclidean vector space $V$ of dimension $n$. \medskip 

Finally, we define $\Pin(V)$ and $\Spin(V)$. Recall the notation $Q_{\epsilon}(V,q) = \{x \in V \mid q(v)  =\epsilon 1\}$, for $\epsilon \in \{+,-\}$.
\begin{definition}
Let $V$ be a Euclidean vector space. Denote by $\Cl^\times (V)$ the multiplicative subgroup of invertible elements of $\Cl(V)$. The \textbf{pin group} $\Pin(V) < \Cl^\times(V)$ is the subgroup of $\Cl^\times(V)$ generated by $Q_+(V) \subset \Cl^\times(V)$. The \textbf{spin group} $\Spin(V)$ is given by $\Spin(V) \coloneqq \Pin(V) \cap \Cl^0(V)$.
\end{definition}

Following the definition, one finds that $\Spin(V)$ admits the following direct characterization: 
\[ \Spin(V) = \left\{ \; \prod_{i=1}^{2k} v_i \in \Cl(V) \mid v_i \in Q_+(V), k \in \Z_+\right \}.\]

The group $\Spin(n)$ is actually the universal cover of $\SO(n)$ when $n \ge 3$.
Indeed there is a (natural) nontrivial 2-1 Lie group covering map $\Spin(n) \rightarrow \SO(n)$ and $\Spin(n)$ is simply-connected when $n \geq 3$. 
See \cite[Ch 1, Theorem 2.9, 2.10]{LM89}.

\subsubsection{Spinor Representations}

We now recall some theory of representations of Clifford algebras and spin groups and some of their features. 
Almost all the material comes from \cite{LM89} and is well-known. 
This expository section is just for the reader's convenience, before we present the construction of transverse spheres using spinors.

The following fact from \cite[Proposition 5.16, Corollary 5.17]{LM89} will be essential in our construction of transverse spheres. 

\begin{proposition}\label{prop:SpinMetric}
Let $c: V \rightarrow \End_{\R}(S)$ be the defining map of a representation of $\Cl(V)$. There exists an inner product $g: S \times S \rightarrow \R$ such that for all $x \in V$, 
\begin{enumerate}[label=(\alph*)]
    \item $c(x) \in O(S, g)$ if $x \in Q_+(V)$.
    \item $c(x)^* =-c(x)$, where the adjoint is with respect to $g$. 
\end{enumerate}
We call such a metric $g$ a \textbf{spin metric}. 
 \end{proposition}
Observe that Proposition \ref{prop:SpinMetric} (a) says that a representation $\Cl(V) \rightarrow \End_{\R}(S)$ and a spin metric $g$ on $S$ induces a representation $\Pin(V) \rightarrow O(S,g)$. 

In fact, these spin metrics carry further properties related to whether the Clifford representation is real, complex, or quaternionic, but we shall not need the additional structure. We note that while \cite{LM89} gives a straightforward (abstract) proof of existence by an averaging argument, these spin metrics also appear naturally, and can be explicitly written down, in the models of \cite{San25}.  

We recall the definition of a \emph{spinor representation} from \cite[Definition 5.11]{LM89}. 
\begin{definition}
A \textbf{spinor representation} $\rho\colon \Spin(n) \rightarrow \GL(W)$ is a representation induced by restriction of an irreducible real representation $\hat{\rho}\colon \Cl^0(n) \rightarrow \End_{\R}(W)$ to $\Spin(n) < \Cl^0(n)$. 
\end{definition}

There are explicit descriptions of $\Cl(n)$ as a matrix algebra $\Mat_m(\K)$, or direct sum $\Mat_m(\K) \oplus \Mat_m(\K)$, for $\K \in \{\R,\C, \Ha\}$, for all integers $n$, where $m$ depends on $n$. 
The key result to this end is the following periodicity property \cite[Theorem 4.3]{LM89}. 
\begin{theorem}[Periodicity of $\Cl(n)$]\label{thm:CliffordPeriodicity}
Let $n \in \Z_+$. Then $\Cl(n+8) \cong_{\R-\mathrm{alg}} \Cl(n) \otimes_{\R} \Cl(8)$. 
\end{theorem}

Hence, by Table \ref{Table:Clifford} and Theorem \ref{thm:CliffordPeriodicity}, one knows the structure of $\Cl(n)$ as an $\R$-algebra for all $n \in \Z_+$. 
We include also in Table \ref{Table:Clifford} the spinor modules $S_n^+$ in the low-dimensional cases, as well as all the low-dimensional concrete models for the abstract group $\Spin(n)$. 
While there is no low-dimensional isomorphism for $\Spin(7), \,\Spin(8)$, these groups are generated by certain multiplications of unit octonions, across $\mathrm{Im} \,\Oct$ and $\Oct$. \medskip 

There is a natural and relevant question of when a representation $\hat{\rho}:\Cl(n) \rightarrow \End_{\R}(W)$ remains irreducible when restricted to $\Spin(n)$. The answer depends on the parity of $n \bmod 8$  \cite[Proposition 5.12]{LM89}. 

\begin{proposition}\label{prop:SpinIrreducible?}
Let $\hat{\rho}\colon \Cl(n) \rightarrow \End_{\R}(S)$ be an irreducible representation and set $\rho \coloneqq \hat{\rho}|_{\Spin(n)}$. Then $\rho$ is irreducible when $n \in \{3,5,6,7\} \bmod 8$ and is a direct sum of two irreducible representations otherwise. 
\end{proposition}

A representation $\hat{\rho}:\Cl(n) \rightarrow \End_{\R}(S)$ has $\hat{\rho}|_{\Spin(n)}$ also irreducible if and only if $\hat{\rho}|_{\Cl^0(n)}$ is irreducible, since $\Spin(n)$ spans the subspace $\Cl^0(n)$. Using this observation, one obtains the proof of Proposition \ref{prop:SpinIrreducible?} via the explicit matrix algebra descriptions of $\Cl(n)$.

By definition, any spinor representation $\eta: \Spin(n) \rightarrow \End_{\R}(W)$ is induced by an irreducible $\Cl^0(n)$-representation. Note that by the explicit nature of $\Cl(n)$ as a matrix algebra, any two irreducible $\Cl(n)$-modules have the same dimension. Since $\Cl^0(n) \cong \Cl(n-1)$, we obtain the following corollary.

\begin{corollary}\label{Cor:SpinIrrep}
Let $S_{n-1}$ be an irreducible $\Cl(n-1)$-module and $S_n^+$ a $\Spin(n)$ spinor representation. Then we have $ \dim_{\R} S_n^+ = \dim_{\R}(S_{n-1})$.
\end{corollary}

\begin{table}[ht] 
\renewcommand{\arraystretch}{1.2}
\centering
\resizebox{0.9\textwidth}{!}{
\begin{tabular}{ |c|c|c|c|c| } 
\hline
Dimension $n$ & Model for $\Cl(n)$ & Model for $\Spin(n)$ & Clifford Module $S_n$ & Spinor Module(s) $S_n^{\pm}$ \\
\hline
 1 & $\Mat_1(\C)$ & $\SO(1)$ & $\C$ & $\R$ \\ 
\hline 
2 & $\Mat_1(\Ha)$ & $\mathrm{U}(1)$ & $\Ha$ & $\C$ \\ 
\hline 
3 &  $\Mat_1(\Ha) \oplus \Mat_1(\Ha)$ & $\Sp(1)$ & $\Ha $ & $\Ha$ \\ 
\hline
4 & $\Mat_2(\Ha)$ & $\Sp(1) \times \Sp(1)$ & $\Ha^2$ & $\Ha^+$, $\Ha^-$ \\
\hline
5 & $\Mat_4(\C)$ & $\Sp(2)$ & $\C^4$ & $\Ha^2$ \\
\hline
6 & $\Mat_8(\R)$ & $\SU(4)$ & $\R^8$ & $\C^4$\\
\hline
7 & $\Mat_8(\R) \oplus \Mat_8(\R)$& $\langle L_{x} \in \SO(\Oct) | \; x \in Q_+(\mathrm{Im} \, \Oct) \rangle $ & $ \R^8$ & $\Oct$ \\
\hline
8 & $\Mat_{16}(\R)$ & $\langle \,L_{x} \oplus R_x \in \SO(\Oct \oplus \Oct) | \; x \in Q_+( \Oct) \rangle $ & $\R^{16}$ & $\Oct^+$, $\Oct^-$\\
\hline
\end{tabular} }
\caption{\emph{\small{An explicit description of $\Cl(n)$ as an $\R$-algebra, irreducible $\Cl(n)$ and $\Spin(n)$-modules, as well as low-dimensional descriptions of $\Spin(n)$ for $1 \leq n \leq 8$ See \cite[Table I \& Table II, Theorem 8.1]{LM89} and \cite{Bry20} or \cite{Bae02} for $n \in \{7,8\}$ and further relations between spinors and octonions.}}}
\label{Table:Clifford}
\end{table}
\medskip

\subsection{Constructing Transverse Spheres}\label{Sec:SpheresWithSpinors}

We now describe some essential properties of the $\Pin(V)$-action on a Clifford module $S$ equipped with a spin metric. 
Property (a) is well-known.

\begin{lemma}[Pin Action on Clifford Module]\label{Lem:PinAction}
Let $c\colon V \rightarrow \End_{\R}(S)$ satisfy the Clifford condition and $\langle \cdot , \cdot \rangle_S$ be a spin metric on $S$. 
Then 
\begin{enumerate}[label=(\alph*)]
    \item (interaction of norms) For any $x,y \in V$ and $w \in S$, 
    \[ \langle c(x) \cdot w, c(y) \cdot w\rangle_{S} = \langle x, y\rangle_{V} \langle w, w\rangle_{S} . \] 
    \item (anti-symmetry) For any $x,y \in V$ and  $v, w \in S$ such that $v\, \bot_{S} \,w$, 
    \[\langle c(x)\cdot v, c(y) \cdot w\rangle_{S} = -\langle c(y)\cdot v, c(x) \cdot w\rangle_{S}. \] 
\end{enumerate}
\end{lemma}

\begin{proof}
(a) Fix any $x \in V, w \in S$. Observe that by the Clifford condition and Proposition \ref{prop:SpinMetric} (b), we have 
\[ \langle c(x) \cdot w, c(x) \cdot w\rangle_{S} = \langle-c(x) \cdot (c(x) \cdot w), w\rangle = \langle q_{V}(x) w, w\rangle = \langle x,x\rangle_{V} \langle w,w\rangle_{S}. \]
Polarizing the identity $||c(x)\cdot w||^2 = ||x||^2\cdot ||w||^2$ leads to 
\[ \langle c(x)\cdot w, c(y) \cdot w\rangle_{S} = \langle x,y\rangle_{V} \langle w,w\rangle_{S} . \]

(b) Using Proposition \ref{prop:SpinMetric} (b) again, observe that 
\begin{align*}
\langle c(x)\cdot v, c(y) \cdot w \rangle &= -\langle c(y) \cdot (c(x) \cdot v), w\rangle = -\langle (c(y)c(x) ) \cdot v,w\rangle  \\
    & =_{(2)} \langle (c(x)c(y))\cdot v, w \rangle + 2 \langle x,y \rangle_{V} \langle v, w\rangle \\
    &= \langle (c(x)c(y))\cdot v, w \rangle = - \langle c(y) \cdot v, c(x) \cdot w\rangle .
\end{align*}
In the step (2), we use the defining Clifford property \eqref{CliffordCondition}. 
\end{proof} 

We now turn to our primary construction of transverse spheres.
In Theorem \ref{Thm:TransverseSpinorSphere} below, we will leverage the spin metric of spinor representations to verify transversality of certain spheres. 
In particular, it is vital to use representations of $\Spin(n)$ induced from $\Cl^0(n)$-representations, rather than, say, non-faithful representations that descend to $\SO(n)$. 
In order to realize this construction in the smallest possible dimension, we want to use the smallest possible (irreducible) representations of the spin group inherited from $\Cl(n)$, which are precisely the spinor representations. 
On the other hand, it is convenient in the proof to record the original irreducible $\Cl(n)$-representation $S_n$ and not just the action of $\Spin(n)$ on its irreducible sub-module $S_n^+$. 

With these remarks in place, we proceed for the setup of the construction of transverse spheres with spinors. 
Let $ \Cl(n) \rightarrow \End_{\R}(S_n)$ be an irreducible representation, and $S_n^+ \subset S_n$ an irreducible $\Spin(n)$-module. 
We shall need the dimension of $S_n^+$. 
By Corollary \ref{Cor:SpinIrrep}, $ \dim_{\R} S_n^+ = \dim_{\R}(S_{n-1})$. Using that $\Cl(8k) \cong \Mat_{16^k}(\R)$, by Theorem \ref{thm:CliffordPeriodicity} and Table \ref{Table:Clifford}, we have the following for $n = 8k+r$ with $0\leq k$ and $1 \leq r\leq 8$: 
\begin{align}\label{SpinModuleDimensions} d(8k+r) \coloneqq\dim_{\R}(S_{8k+r}^+) =\dim_{\R}(S_{8k+r-1}) = \begin{cases}
    16^k,  & r=1\\
     2\cdot 16^k , & r= 2\\
      4\cdot 16^k,  & r \in \{3,4\}\\
       8\cdot 16^k,  & r \in \{5,6,7,8\}.\\
\end{cases} \end{align} 

We now introduce the essential definitions for the theorem. 
Let $c_n \colon V_n \rightarrow S_n$ be the defining map of an irreducible representation $\eta \colon \Cl(n) \rightarrow \End_{\R}(S_n)$. 
Take any spin metric $g_n$ on $S_n$ and let $S_n^+ \subset S_n$ be an irreducible $\Spin(n)$-submodule. 
Choose any spinor $x_0 \in Q_+(V_n)$ and define the map $f \colon Q_+(V_n) \rightarrow \Spin(n)$ by $f(x) = x_0x$. 
Now, for any positive integer $p$, define 
\[ W_n \coloneqq \R^{p,0} \oplus (S_n^+, -g_{S_n}) \cong \R^{p, d}. \] 
Consider the map 
$\iota \colon \Spin(V_n) \hookrightarrow \SO(S_n^+) \hookrightarrow \SO(W_n)$, where the last inclusion is the map $g \mapsto \id_{\R^p} \oplus g$. 
Combining $\iota$ and $f$, we build a map $\beta \colon Q_+(V_n) \rightarrow \SO(W_n)$ by $\beta \coloneqq \iota \circ f.$

\begin{theorem}[Transverse Spheres from Spinors]\label{Thm:TransverseSpinorSphere}
Let $\mathcal{F} = \Flag(\R^{p,d})$ for any positive integer $p$, where $d= d(n)$ as in \eqref{SpinModuleDimensions}. 
Fix any full flag $F \in \mathcal{F}$ and $\Lambda \coloneqq \beta(Q_+(V_n)) \cdot F$ is a transverse $(n-1)$-sphere in $\mathcal{F}$.
\end{theorem}

\begin{proof}

Choose any orthonormal basis $(s_i)_{i=1}^{d}$ of $S_n^+$. 
Next, let $(e_i)_{i=1}^{p}$ be an orthonormal basis of $\R^{p,0}$. 
The pair $ (e_i),(s_i),$ yields a flag $F=(F^i)_{i=1}^{\min(p,d)} \in \Flag(\R^{p,d})$. 
That is, define $z_i = e_i + s_i$ and set $F^i \coloneqq \spann_{\R} \{ z_j \}_{j=1}^i$. When $p=d$, if $(e_i)$ is chosen well, we obtain $F \in \Flag_+(\R^{d,d})$, and the following arguments holds without alteration. 
By construction, each $F^i$ is an isotropic $i$-plane and $F \in \Flag(\R^{p,d})$. 
We will show that $\image(\beta) \cdot F$ is a transverse $(n-1)$-sphere.

Choose any $x \neq y \in Q_+(V_n)$. 
Then define $Z_i(x) = \beta(x) \cdot z_i$ and $Z_i(y) = \beta(y) \cdot z_i$. 
To examine the transversality conditions between $\beta(x) \cdot F$ and $\beta(y) \cdot F$, we form the matrix of dot products $A = A(x,y) = (a_{ij})_{i,j=1}^{\min(p,d)}$ with entries $a_{ij}$, where $a_{ij} = \langle Z_i(x) , Z_j(y) \rangle_{W_n}$. 
Since $c_n(x_0) \in O(S_n, g_n)$, observe the following: 
\[ \begin{cases}
    a_{ii} &= 1 -\langle c_n(x) \cdot s_i, c_n(y) \cdot s_i \rangle_{S_n} \\ 
    a_{ij} &= \langle c_n(x)\cdot s_i, c_n(y) \cdot s_j \rangle_{S_n}
\end{cases} \] 
By Lemma \ref{Lem:PinAction} (a), we see that 
$a_{ii} =1- \langle x,y\rangle_{V_n} > 0$ since $x \neq y$. 
Also, Lemma \ref{Lem:PinAction} (b) gives $a_{ij} = -a_{ji}.$ 
We conclude that $\det(A) > 0 $ by Lemma \ref{Lem:ElementaryLinAlg}. 
Moreover, the same arguments imply $\det(A_m) > 0 $ for any of the $(m\times m)$-minors $A_m = (a_{ij})_{i,j =1}^m$, where $1 \leq m \leq \min(p,d)$. 
Now, $\det(A_m) > 0$ verifies the transversality condition between the isotropic $m$-planes $\beta(x) \cdot F^m$ and $\beta(y) \cdot F^m$. 
Hence, the full flags $\beta(x) \cdot F$ and $\beta(y) \cdot F$ are transverse.  
\end{proof}

Here is the elementary linear algebra lemma used in the proof of Theorem \ref{Thm:TransverseSpinorSphere}. 

\begin{lemma}\label{Lem:ElementaryLinAlg}
    Let $A=X+cI \in \Mat_n(\R)$ with $ c > 0 \in \R$ and $X$ skew-symmetric. Then $\det(A) > 0$. 
\end{lemma}

\begin{proof}
    The eigenvalues of $X$ are purely imaginary. Now, let $c_{M}(\lambda) \coloneqq \det(M- \lambda I)$ be the characteristic polynomial of a matrix $M$. Then
    \[ c_A(\lambda) = \det ( (c-\lambda)I+X) = c_{X}(\lambda-c).\]
    Thus, the eigenvalues of $A$ are $c$-translates of eigenvalues of $X$ and hence are all nonzero.
    Now, we prove $\det(A) > 0$. 
    Suppose the roots of $X$ are $\alpha_1 i, -\alpha_1i, \alpha_2i, -\alpha_2 i, \dots, -\alpha_ki, \alpha_ki$, with $0$ appearing $j$ times. Then one finds $\det(A) = c^j \prod_{j=1}^k (c^2+\alpha_j^2) >0$. 
\end{proof}

for us, the most important cases of Theorem \ref{Thm:TransverseSpinorSphere} occur when $p \in \{d-1,d\}$, so that the associated group $G= \SO_0(p,d)$ is split. 
We will bootstrap off these cases to produce many further transverse spheres. 

As it stands, we describe $d=d(n)$ such that the $B_{d-1}$ or $D_d$ full flag manifold carries a transverse $(n-1)$-sphere. It is useful for later purposes to invert this relationship. The Radon-Hurwitz numbers from \eqref{RadonHurwitz} encode the inversion as follows. In particular, for $d$ a power of two, we have:
\begin{align}\label{RadonHurwitzInversion}
    \begin{cases}
        \rho(16^k)&= 8k+1\\
        \rho(2\cdot 16^k)&= 8k+2\\
        \rho(4 \cdot 16^k)&= 8k+4\\
        \rho(8 \cdot 16^k)&= 8k+8.
    \end{cases}
\end{align}
Thus, a repackaging of Theorem \ref{Thm:TransverseSpinorSphere} is as follows.
\begin{corollary}\label{Cor:IsotropicTransverseSphere}
Let $n = 2^j$ be a power of two, for $j \geq 2$. Then there is a transverse $(\rho(n)-1)$-sphere in $\Flag(\R^{n,n+\epsilon})$ for $\epsilon \in \{-1,0,1\}$.
\end{corollary}

\begin{proof}
The cases $\epsilon \in \{0,-1\}$ are explicitly addressed by Theorem \ref{Thm:TransverseSpinorSphere}. 
The transverses inclusion $\Flag_+(\R^{2n,2n}) \stackrel{\pitchfork}{\hookrightarrow} \Flag(\R^{2n,2n+1})$ from Example \ref{Ex:D2B2n} handles the $\epsilon = 1$ case. 
\end{proof}

\subsection{Deformations of Transverse Spheres}\label{Sec:Deformations}

The construction of transverse spheres in the previous section is algebraic. 
In this section, we show the construction is rather flexible by demonstrating some ``generic" deformations. 
In Section \ref{Sec:Obstructions} we will show the original transverse spheres in Theorem \ref{Thm:TransverseSpinorSphere} (or their descendants constructed in later sections) cannot arise as Anosov limit sets, but those obstructions fail on the generic deformations we obtain here. 

We now sketch the main idea of the deformation. 
For concreteness, we restrict attention to the case of transverse spheres in the full flag manifold of $\R^{d,d}$.
The key is to now allow the lifted sphere in the maximal compact $K =\SO(d) \times \SO(d)$ to move the spacelike vectors, but only a small amount in comparison to the timelike vectors. First, we split $\R^{d,d} = \R^{d,0} \oplus \R^{0,d}$ and identify each factor with the same irreducible $\Spin(n)$-module $S_n^+$. 
Afterwards, select a suitably contracting map $\phi: \mathbb{S}^{n-1} \rightarrow \mathbb{S}^{n-1}$ and then define a map $\psi \colon \mathbb{S}^{n-1} \rightarrow \SO(d) \times \SO(d)$ by taking the graph of $\phi$, combined with the Clifford action. 
Again, let $f \colon \mathbb{S}^{n-1} \rightarrow \SO(d)$ be the map $f(x) = \eta(x_0x)$ for $\eta \colon \Spin(n) \rightarrow \SO(d)$ an irreducible spinor representation. Finally, the map $\psi$ is given by 
\begin{align}\label{PsiDeformation}
    \psi(x) = (f(\phi(x)), f(x)). 
\end{align}
We will show the the orbit action of $\psi$ on a single flag yields a transverse $(n-1)$-sphere. 

\begin{lemma}[Contraction implies Transverse]\label{Lem:Deformation}
Suppose that $\phi: \mathbb{S}^{n-1} \rightarrow \mathbb{S}^{n-1}$ satisfies 
\begin{align}\label{Contraction}
\langle \phi(x), \phi(y) \rangle_{\R^{n,0}} > \langle x, y \rangle_{\R^{n,0}} \; \;(x \neq y). 
\end{align}
Take any flag $F_0 \in \Flag(\R^{d,d})$. Then the map $\xi: \mathbb{S}^{n-1} \rightarrow \Flag(\R^{d,d})$ by $\xi(x) = \psi(x) \cdot F_0$, for $\psi$ defined by \eqref{PsiDeformation}, is a transverse $(n-1)$-sphere. 
\end{lemma}

\begin{proof}
We follow the proof strategy of Theorem \ref{Thm:TransverseSpinorSphere}. Here, let $\eta: \Cl(n) \rightarrow \End_{\R}(S_n)$ be an irreducible representation that restricts to the spinor representation $\eta:\Spin(n) \rightarrow \End_{\R}(S_n^+)$, abusively denoted by the same symbol. To simplify notation, for $x \in \Cl(n), w \in S_n$, we write $x \cdot w$ to denote $\eta(x)(w)$. Let $g$ be a spin metric on $S_n$ and recall $d = \dim_{\R}(S_n^+)$. Define also $\R^{2d,2d} \cong (S_n, g) \oplus (S_n, -g)$. 
Write $F_0 = (F_0^i)_{i=1}^d$, with $F_0^i = \spann \{e_i + s_i\}_{j=1}^i$, with $(e_i), (s_i)$ orthonormal bases of $S_n^+$. Define $w_i = e_i+s_i$ and then set $a_{ij}(x,y) \coloneqq \langle \psi(x) \cdot w_i, \psi(y) \cdot w_j \rangle_{\R^{d,d}}$. Observe that since $\eta(x_0) \in O(S_n)$ by Proposition \ref{prop:SpinMetric}, we have 
\begin{align}\label{SymmetryOfDeformation}
    a_{ii} &= \langle \psi(x)\cdot w_i , \psi(y) \cdot w_i \rangle_{\R^{d,d}} = \langle \phi(x) \cdot e_i, \phi(y) \cdot e_i\rangle_{g} -\langle x\cdot s_i, y\cdot s_i\rangle_{g} = \langle \phi(x),\phi(y)\rangle_{g}-\langle x, y\rangle_{g} > 0
\end{align}
by the hypothesis \eqref{Contraction}. Then, similar to the anti-symmetry in Theorem \ref{Thm:TransverseSpinorSphere}, one finds 
\begin{align*}
    a_{ij} &= \langle \psi(x) \cdot w_i, \psi(y) \cdot w_j\rangle_{\R^{d,d}} = \langle \phi(x) \cdot e_i, \phi(y)\cdot e_j \rangle_{g} +\langle x \cdot s_i, y \cdot s_j \rangle_{-g} \\
    &= -\langle \phi(y)\cdot e_i, \phi(x) \cdot e_j\rangle_{g} -\langle y \cdot s_i , x \cdot s_j\rangle_{-g} = -\langle \psi(y) \cdot w_i, \psi(x) \cdot w_j \rangle_{\R^{d,d}} = -a_{ji}.
\end{align*}
Set $A = (a_{ij})_{i,j=1}^d$ and let $A_m = (a_{ij})_{i,j=1}^m$. Then we conclude that $A_1 > 0$ by \eqref{SymmetryOfDeformation} and we have $\det(A_i) > 0$ by Lemma \ref{Lem:ElementaryLinAlg}. Once again, we conclude $\xi(x)^i$ and $\xi(y)^i$ are transverse for all $x \neq y \in \mathbb{S}^{n-1}$ and $1\leq i \leq d$, so that the map $\xi$ is transverse. 
\end{proof}

In light of the lemma, Theorem \ref{Thm:TransverseSpinorSphere} is just the case $\phi$ is the constant map $\phi(x)=-x_0$.

\begin{remark}
    Condition \eqref{Contraction} on $\phi$ in Corollary \ref{Lem:Deformation} is equivalent to the condition that $\phi$ is strictly 1-Lipschitz with respect to the round metric on the sphere. 
    This is reminiscent of the well-known fact in pseudo-hyperbolic geometry that a connected transverse subset of the Einstein universe obtains the form of a graph of a strictly $1$-Lipschitz map between spheres over any splitting $\R^{p,q} = \R^{p,0} + \R^{0,q}$, see for instance \cite{Mes07,bonsante2010maximal,DGK18}. 
    Note that the image of any such map is contained in a hemisphere \cite[Lemma 2.8]{Hei05}.
\end{remark}

Restricting attention to the case of $\R^{d,d}$ with $d \equiv 0 \bmod 4$, we produce transverse $3,7$-spheres whose span in the Einstein universe is linearly full. 

\begin{theorem}[Linearly Full Deformations]\label{thm:Full37Spheres}
    Let $d \geq 4$ and $\mathrm{pr}_1 \colon \Flag(\R^{d,d}) \rightarrow \Ein^{d-1,d-1}$ be the natural projection. 
    For $n \in \{4,8\}$, if $d \equiv 0 \bmod n$, then there is a transverse $(n-1)$-sphere $\Lambda \subset \Flag(\R^{d,d})$ such that $\spann \pr_{1}(\Lambda) = \R^{d,d}$.
\end{theorem}

\begin{proof}
Suppose $n \in \{4,8\}$. 
The key property is that in these cases, the map $f \colon \mathbb{S}^{n-1} \cong Q_+(\mathbb{K}) \rightarrow \SO(\mathbb{K})$ given by $x \mapsto L_x$, for $\K \in \{\Ha, \Oct\}$, composes with the orbit map to yield a homeomorphism $\mathbb{S}^{n-1} \rightarrow \mathbb{S}^{n-1}$. \medskip 

\textbf{Step 1: Folding}. 

Write $d=nk$. 
We start by producing, for $1 \le i \le k$, explicit maps $F_i:\mathbb{S}^{n-1} \rightarrow \mathbb{S}^{n-1}$ that are strictly $1$-Lipschitz and $2^i$-to-$1$ on some open subset. The maps $F_i$ will be compositions of `elementary folds' $f_j$. 

Choose a $2$-plane $U = \spann\{e_1,e_2\}$ in $\R^n$, with an isometric action of the dihedral group $W$ of order $2^k$. 
It is convenient to fix an orientation on $U$.
Extend the action of $W$ trivially on $U^\perp \subset \R^n$. 
Let $D$ denote the closure of a fundamental domain for the action of $W$ on $\R^n$. 
Write the translates of $D$ as $D_1, \dots, D_{2^k}$, with $D= D_1$, in the cyclic order compatible with the orientation of $U$.
Define $H_i $ as the hyperplane spanned by $D_{2^{i-1}} \cap D_{2^{i-1}+1}$. 
Let $\alpha_i: \R^n \rightarrow \R$ denote a linear map with $\ker(\alpha_i) = H_i$, which is positive on the interior of $D$. 
Let $R_i:\R^n \rightarrow \R^n$ denote the reflection in $H_i$ and define the `elementary fold' $f_i: \mathbb{S}^{n-1} \to \mathbb{S}^{n-1}$ by 
\[ f_i(x) = \begin{cases}
     x  & \text{if} \; \alpha_i(x) \ge 0 \\
     R_i(x)  & \text{otherwise.} \end{cases} \]
Define $\hat{F}_i: \mathbb{S}^{n-1} \rightarrow \mathbb{S}^{n-1}$ by $\hat{F}_i \coloneqq f_{k} \circ f_{k-1} \circ \cdots \circ f_{k-i+1}$. For $1 \le i \le k$, let $\Omega_i$ be given by $\Omega_{i} = \bigcup_{j=1}^{2^{i-1}} (D_{j} \cap \mathbb{S}^{n-1})$. 
Observe that $\Omega_{k+1-i}$ is the image of $\hat{F}_i$, and for $i \le j$, we have $\hat{F}_i \circ R_j = \hat{F}_i$. 
Additionally, for all indices $1\leq i \leq k$, the map $\hat{F}_i$ is $2^{i}$-to-$1$ on $\Omega \coloneqq \bigcup_{j=1}^{2^k} (\mathbb{S}^{n-1} \cap \mathrm{int}(D_j))$. 

The folds $F_j$ will be post-compositions of $\hat{F}_j$ by a strictly contracting map. Let $c: \Omega_k \rightarrow \Omega_k$ be a strictly 1-Lipschitz map given by contracting the geodesic $\frac{\pi}{2}$-ball $\Omega_k$ to a geodesic $\frac{\pi}{4}$-ball with the same center. Define $F_i \coloneqq c \circ \hat{F}_i$. 
We emphasize for later use that $F_i \circ R_j = F_i$ exactly when $i \le j$. More precisely, if $1\leq j < i$, then $F_i \circ R_j(x) \neq F_i(x)$ for $x \in \mathrm{int}(\Omega_{k+1-i}) \backslash H_j$. 
\medskip

\textbf{Step 2: Define the sphere}. 

We build the transverse sphere in $\Flag(\R^{d,d})$ now. Recall that $d= nk$. Let $\psi_l: \mathbb{S}^{n-1} \rightarrow \mathbb{S}^{n-1}$ be homeomorphisms, for $1 \leq l \leq k$, to be prescribed later. Let $f:\mathbb{S}^{n-1} \rightarrow \SO(d)$ now be the usual map $f(x) = \eta(x_0x)$. 
Split $\R^{d,d} = (\R^{n,0})^k \oplus (\R^{0,n})^k$ and correspondingly we define a map
$\hat{\xi}: \mathbb{S}^{n-1} \rightarrow \SO(n,0)^k\times \SO(0,n)^k < \SO(d) \times \SO(d)$ by 
\begin{align}\label{DeformedSphereLift}
 \hat{\xi}(x) = (F_k, \,F_{k-1},\,\dots, \,F_{1},\,f\circ \psi_k,\dots, \,f\circ \psi_1). 
\end{align}
We build a sphere in $\Flag(\R^{d,d})$ by the orbit map $\xi(x) = \hat{\xi}(x)(F^{\bullet})$ for a flag $F^{\bullet} \in \Flag(\R^{d,d})$. 
As usual, we write 
$F^m = \spann \{e_j+s_j\}_{j=1}^m$ for some orthonormal bases $(e_i)_{i=1}^d$ and $(s_i)_{i=1}^d$ of $\R^{d,0}$ and $\R^{0,d}$. 
This time, we do make a prescription on $F^{\bullet}$. 
Under the splittings $\R^{d,0} = (\R^{n,0})^k$ and $\R^{0,d} = (\R^{0,n})^k$, we demand $e_1$ and $s_1$, have nonzero projection on every $\R^{n,0}$ and $\R^{0,n}$ factor, respectively. 

We cannot invoke the direct sum construction, namely Corollary \ref{Cor:IsotropicDirectSumMaps}, to see that $\hat{\xi}\cdot F^{\bullet}$ gives a transverse sphere, because $F^{\bullet}$ is not itself a direct sum of flags coming from the splitting $\R^{d,d} = (\R^{n,n})^k$. 
Instead, transversality follows from a modification of the proof of Theorem \ref{Thm:TransverseSpinorSphere}, as we now explain. 
First, we show the non-deformed map $\hat{\xi}_{0}\colon \mathbb{S}^{n-1}\rightarrow \SO(n)^k\times \SO(n)^k$
\[ \hat{\xi}_{0}(x) = (\id, \id, \dots, \id, f\circ \psi_k, \cdots, f\circ \psi_1), \]
yields a transverse sphere when acting on $F^{\bullet}$.
To this end, we once again form the matrix of dot products $A(x,y) = (a_{ij}(x,y))$, for $x \neq y \in \mathbb{S}^{n-1}$, where, for $1\le i,j \le d$: 
\[ a_{ij}(x,y) = \langle \hat{\xi}_{0}(x) (e_i+s_i), \hat{\xi}_{0}(y) (e_j+s_j) \rangle_{\R^{d,d}}. \]
The matrix $A$ naturally decomposes into pieces. 
We write $e_i = \sum_{l=1}^k b_le^l_i$ and $s_i = \sum_{l=1}^kc_l s^l_i$, where $e_i^l, s_i^l$ live in the $l^{th}$ $\R^{n,0}$ or $\R^{0,n}$ block of $\R^{d,0}$ or $\R^{0,d}$, respectively. 
Then $A = \sum_{l=1}^k A_l$, where $A_l = (a_{ij}^l)_{i,j=1}^d$ and 
\[a_{ij}^l \coloneqq a_{ij}^l(x,y) = \langle \hat{\xi}_{0}(x)  (e_{i}^l+s_i^l), \hat{\xi}_{0}(y) (e_j^l+s_j^l) \,\rangle_{\R^{d,d}}.\] 
Each sub-matrix $A_{ij}^l$ obtains the form $A_{ij}^l = \alpha_l I+X_l$ for some scalar $\alpha_l \geq 0$ and $X_l$ skew-symmetric as in the proof of Theorem \ref{Thm:TransverseSpinorSphere}. Hence, $A = \sum_{l=1}^k A_l$ is also a sum of a scalar and skew-symmetric matrix, with diagonal entries given by  
\begin{align}
    a_{11} = \sum_{l=1}^k \alpha_l&=\sum_{l=1}^k \langle \hat{\xi}_{0}(x)(e_1^l+s^l_1), \hat{\xi}_{0}(y)  (e_1^l+s_1^l)\,\rangle_{\R^{d,d}} = 1-\sum_{l=1}^k \langle \psi_l(x)\cdot s^l_1, \psi_l(y)\cdot s^l_1\rangle_{\R^{d,0}} \\
    &= 1-\sum_{l=1}^k \langle \psi_l(x), \psi_l(y)\rangle_{\R^{d,0}}||s_1^l||^2 >1-\sum_{l=1}^k ||s_1^l||^2 = 0.
\end{align}
Invoking Lemma \ref{Lem:ElementaryLinAlg} again, we find that $\det(A) > 0$, and $\det(A_m) > 0$ for $1\leq m \leq d$, for any of the $(m \times m)$-minors $A_m = (a_{ij})_{i,j=1}^m$. As in the proof of Theorem \ref{Thm:TransverseSpinorSphere}, we conclude that the orbit of $\hat{\xi}_{0}$ defines a transverse $(n-1)$-sphere in $\Flag(\R^{d,d})$. Following the proof of Lemma \ref{Lem:Deformation}, one sees the deformation 
$\hat{\xi}$ of $\hat{\xi}_{0}$ also yields a transverse $(n-1)$-sphere in $\Flag(\R^{d,d})$. 

\medskip 

\textbf{Step 3: Check the linear fullness}. Define $\Lambda \coloneqq \image(\xi)$ and $V_1\coloneqq \spann ( \pr_1(\Lambda))$. We will show $V_1 = \R^{d,d}$. 

We first show that $\R^{0,d} \subset V_1$.
We record the decomposition $\R^{0,d} = (\R^{0,n})^k$ as $\R^{0,d} = \bigoplus_{i=1}^k U_i$. 
Choose a $d$-tuple of pairwise distinct points $\{x_1,\dots,x_d\}$ in the interior of $\Omega_1$ and note that $\{x_1,\dots,x_d,R_k(x_1),\dots,R_k(x_d)\}$ consists of $2d$ pairwise distinct points of $\mathbb{S}^{n-1}$. 

We now choose the suitable self-homeomorphisms $(\psi_l)_{l=1}^k$ of $\mathbb{S}^{n-1}$ in \eqref{DeformedSphereLift}, used to define $\hat{\xi}$.
We may choose $\psi_l$ isotopic to the identity such that $\{ \psi_l(x_j), \psi_l(R_k(x_j) \}_{j=1}^d$ is any arbitrary distinct $2d$ points on $\mathbb{S}^{n-1}$.
In particular, for each $1 \leq l \leq k$, we can choose $\psi_l$ wisely enough that the differences $\psi_l(x_j)-\psi_l(R_k(x_j))$, for $1 \leq j \leq d$, span $\R^n$. 
For any $y \in \mathbb{S}^{n-1}$, recall that the orbit map $\mathbb{S}^{n-1} \to \mathbb{S}^{n-1}$ by $x\mapsto xy$ is a homeomorphism, identifying again $Q_+(\K) \cong \mathbb{S}^{n-1}$. 
Hence, we can arrange that for fixed $l$, $Y_{lj}=(\psi_l(x_j)-\psi_l(R_k(x_j)))s_1^l$ span $U_l$. 

Additionally, we can can choose $\psi_l$ so that the elements $Y_j \coloneqq \sum_{l=1}^k Y_{lj}$ form a basis $(Y_j)_{j=1}^d$ of $\R^{d}$. 
Recall from Step 1 that $F_l \circ R_k = F_l$ for all $1\leq l \leq k$. 
Now, consider suitable differences of elements of the image of $\hat{\xi}\cdot (e_1+s_1)$: 
\begin{align}\label{NegativeSpan}
    (\hat{\xi}(x_j)-\hat{\xi}(R_k(x_j)))\cdot (e_1 +s_1)= \sum_{l=1}^k \big(\psi_l(x_j)-\psi_l(R_k(x_j))\big)s_1^l = \sum_{l=1}^k Y_{lj} = Y_{j}.
\end{align}
Hence, $Y_j \in V_1$ for all $1 \leq j \leq d$ and then $\R^{0,d} \subset V_1$. \medskip

To conclude, it suffices to show that $\R^{d,0} \subset V_1$. Write $\R^{d,0} = \bigoplus_{i=1}^k U_i'$, with $U_i' \cong \R^{n,0}$ the direct sum factors used in defining $\hat{\xi}$. 
Recall that $F_l \circ R_{k-1} = F_l$ for $l \leq k-1$, while $F_k \circ R_{k-1}(x) \ne F_{k}(x)$ for any $x \in \mathrm{int}(\Omega_{1})$.
Then 
    \[ (F_k(x) - F_k(R_{k-1}(x)))e_1^1 = (\hat{\xi}(x)-\hat{\xi}(R_{k-1}(x)))e_1 \in V_1 \]
and since $\image(F_k - F_k \circ R_{k-1})$ spans $\R^n$, the orbit of $e_1^1$ under $F_k-F_k\circ R_{k-1}$ spans $U_1'$.
In general, we find that for any $x \in \mathrm{int}(\Omega_{1})$
    \[ \sum_{j=k-l+1}^{k} (F_j(x) - F_j(R_{k-l}(x)))e_1^j = (\hat{\xi}(x)-\hat{\xi}(R_{k-l}(x)))e_1 \in V_1. \]
Hence, $U_1' \oplus \dots \oplus U_l' \subset V_1$ by induction. 
We conclude that $U_l' \subset  V_1$ for all $1\leq l \leq k$, meaning $\R^{d,0} \subset V_1$.
\end{proof}

We include one more special case of deformations of interest. 
\begin{corollary}\label{Cor:G2deformation}
There exists a transverse 3-sphere $\Lambda \subset \Flag(\Gtwosplit)$ such that $\spann(\pr_1(\Lambda))= \R^{3,4}$, where $\pr_1: \Flag(\Gtwosplit) \rightarrow \Ein^{2,3}$ is the natural projection. 
\end{corollary}

\begin{proof}
The reasoning is similar to that of Theorem \ref{thm:Full37Spheres}. 
Fix a splitting $\R^{3,4} = \R^{1,0}\oplus \R^{2,0} \oplus \R^{0,4}$ and take a full flag $F \coloneqq (\ell, \omega) \in \mathcal{F}_{1,2}^\times \cong \Flag(\Gtwosplit)$ such that $\ell = [e_1+s_1]$ with $e_1 = e_1^1+e_1^2$ a unit vector in $\R^{3,0}$ such that $e_1^i \neq 0$ for $i \in \{1,2\}$, with $s_1 \in Q_-(\R^{0,4})$ arbitrary. Let $\SO(4) \cong K < \Gtwosplit $ be the maximal compact subgroup preserving the splitting $\R^{3,4}= \R^{3,0} \oplus \R^{0,4}$. Recall that topologically $\SO(4) = \SO(3) \times \Sp(1)$. We obtain a diffeomorphism $\varphi:\SO(4) \rightarrow \SO(3) \times \Sp(1)$ of $K$ via the map $\varphi(g) = (g|_{\R^{3,0}}, g|_{\R^{0,4}})$. We will define a continuous map $\hat{\xi}: \mathbb{S}^3 \rightarrow \SO(4)$ such that the 3-sphere $\Lambda$ of is given by $\hat{\xi} \cdot F$. Once more, define $V_1 = \spann(\pr_1(\Lambda))$, so that we must show $V_1 = \R^{3,4}$.  

Next, take a strictly 1-Lipschitz map $F_1: \mathbb{S}^3 \rightarrow \SO(2)$, where $\SO(2)$ is identified with $\mathbb{S}^1$ with its usual distance function. One way to construct the map $F_1$ is to define $F_1(y_1, y_2, y_3, y_4) = \frac{1}{2} y_4$, where we identify $\mathbb{S}^1 = [0,1]/(0\sim 1)$. We also imagine $\SO(2) \hookrightarrow \SO(3)$ induced by the splitting $\R^{3,0} = \R^{1,0} \oplus \R^{2,0}$. 
Finally, the map $\hat{\xi}$ is given by $\hat{\xi}(x) = (F_1(x), \eta(x_0\psi_1(x))) \in \SO(3) \times \Sp(1)$, where again $\eta: Q_+(\Ha) \rightarrow \Sp(1)$ is the map $x \mapsto L_x$ and $\psi_1:\mathbb{S}^3 \rightarrow \mathbb{S}^3$ is a homeomorphism to be chosen later. The subset $\Lambda = \hat{\xi} \cdot F$ is a transverse 3-sphere by Theorem \ref{Thm:G2Fibration} and Lemma \ref{Lem:Deformation}. 

Observe that for all $x,y \in H_4^+ = \{ x \in \mathbb{S}^3 \mid x_4> 0 \}$ with $x_4=  y_4$, we have $F_1(x) = F_1(y)$. Hence, as in the proof of Theorem \ref{thm:Full37Spheres}, see that for any $x\neq y \in H_4^+$ with $x_4=y_4$, we have  
\begin{align}\label{timelikevectors}
(\psi_1(x)-\psi_1(y))s_1 =\big(\hat{\xi}(x) -\hat{\xi}(y)\big) \cdot (e_1+s_1)  \in V_1. 
\end{align}
Just as in the proof of Theorem \ref{thm:Full37Spheres}, prescribing finitely many points $(\psi_1(x_j), \psi_1(y_j))$ appropriately, for $x_j, y_j \in H_4^+$, can force the above elements in \eqref{timelikevectors} to span $\R^{0,4}$. 
Thus, we conclude $\R^{0,4} \subset V_1$ for $\psi_1$ well-chosen. 
Then across $x, y \in H_4^+$, the elements 
\[ \big(\hat{\xi}(x)-\hat{\xi}(y))\big) \cdot e_1=\big(F_1(x)-F_1(y)\big) e_1^2\]
will span $\R^{2,0}$. Hence, $\R^{1,0} \subset V_1$ as well and we conclude that $V_1 = \R^{3,4}$. 
\end{proof}

Observe that for $\Lambda$ one of the spheres from Theorem \ref{thm:Full37Spheres} or Corollary \ref{Cor:G2deformation}, if $\Lambda$ were the flag limit set of a Borel Anosov representation $\eta: \Gamma \rightarrow G$, then $\eta$ would have to be irreducible. This irreducibility stands in direct contrast to the spheres in Theorem \ref{Thm:TransverseSpinorSphere}, which we can easily obstruct in Section \ref{Sec:Obstructions}, by arguing that if $\Lambda$ were the flag limit set of a representation $\eta$, then $\eta$ would have to be highly reducible. 

We expect there are many further constructions of transverse spheres with sufficiently `generic' properties to force groups preserving them to be irreducible.
Since at this time it is not even clear whether one can either obstruct or realize the 3-spheres and 7-spheres from Theorem \ref{thm:Full37Spheres} as Anosov flag limit sets, we do not pursue any further deformations in the present work.

\section{Transversality in \texorpdfstring{$\Gtwosplit$}{G2'}-Flag Manifolds} \label{Sec:G2Transversality}

In this section, we describe a family of principal $\Sp(1)$-bundle fibrations for the the full flag manifold $\Flag(\Gtwosplit) \cong \Gtwofullflags$, indexed by points $P$ in the $\Gtwosplit$-symmetric space. 
The fibers of these fibrations are topological 3-spheres, which turn out to be maximally transverse 3-spheres in $\Gtwofullflags$. 
This is no accident: these transverse spheres are an instance of the construction in Theorem \ref{Thm:TransverseSpinorSphere} when $n=4$. 
We then describe the associated fibrations of the partial flag manifolds $\Ein^{2,3}$ and $\Pho^\times$, encoded in  Figure \ref{Fig:GtwoFlagFibrationsSimple}.

\subsection{Reminders on the Exceptional Group \texorpdfstring{$\Gtwosplit$}{G2'}}\label{Sec:G2Prelims}

The group $\Gtwosplit$ admits a variety of useful characterizations. 
In this section, we study $\Gtwosplit$ through the cross-product on $\R^{3,4}$ and suppress the relationship with the split-octonions $\Oct'$, though the latter is discussed briefly in Section \ref{Subsec:CompositionAlgebras}. 
A comprehensive introduction to $\mathrm{G}_2$, including the complex group $\mathrm{G}_2^{\C}$ and its split and compact real forms $\Gtwosplit$ and $\mathrm{G}_2^c$, can be found in \cite{Fon18}. 

Let us now discuss the \emph{cross-product} $\times_{3,4}\colon \R^{3,4} \times \R^{3,4} \rightarrow \R^{3,4}$ on pseudo-Euclidean space $\R^{3,4} = (\R^{7}, q)$, where $q \coloneqq q_{3,4}$ is a signature $(3,4)$ bilinear form on $\R^7$. 
The map $\times_{3,4}$ is (i) bilinear, (ii) alternating, (iii) $q$-orthogonal map, (iv) normalized by $q(u\times_{3,4} v) = q(u)q(v)-q(u,v)^2$. 
The group $\Gtwosplit$ may be defined as the group of cross-product preserving transformations of $(\R^{3,4}, \times_{3,4})$: 
\begin{align}
    \Gtwosplit \coloneqq \Aut(\R^{3,4}, \times_{3,4})=  \{ g \in \GL(\R^{3,4}) \; | \; g(u \times_{3,4} v) = g(u) \times_{3,4} g(v) \}. 
\end{align}
The group $\Gtwosplit$ above is indeed the split real (adjoint) form of $\Gtwo^{\C}$. We will frequently denote by $\mathcal{C}_x: \R^{3,4} \rightarrow \R^{3,4}$ the map $\mathcal{C}_x(y)= x\times y$, the cross-product endomorphism of $x$. 

The above definition of $\Gtwosplit$ seems more natural given the following: if $\wedge_{3,4} \colon \R^{3,4} \times \R^{3,4} \rightarrow \R^{3,4}$ is any other map satisfying (i)-(iv), then there is a linear isomorphism $(\R^{3,4}, \times_{3,4}) \cong (\R^{3,4}, \wedge_{3,4})$. This uniqueness of cross-products on $\R^{3,4}$ is a striking corollary of Hurwitz's theorem that classifies  composition algebras over $\R$; see Section \ref{Subsec:CompositionAlgebras} for more details. Moreover, the exceptional nature of $\Gtwosplit$ can be traced back to the scarcity of cross-products in the aforementioned sense. Indeed, the only pseudo-Euclidean spaces $\R^{p,q}$ admitting a cross-product (excluding the trivial one-dimensional cases) satisfying axioms (i)-(iv) are those of signature $(3,0),(7,0), (1,2), (3,4)$; see Section \ref{Subsec:CompositionAlgebras}, specifically Proposition \ref{Prop:CrossProductCompAlg}, for more details. 

The group $\Gtwosplit$ is a connected (cf. Proposition \ref{Prop:GtwosplitFirstStiefel}), simple Lie group with trivial center and fundamental group $\pi_1(\Gtwosplit) = \mathbb{Z}_2$ \cite[Proposition 5.17]{Fon18}. 
The Lie algebra $\g_2'$ consists of derivations of $\times_{3,4}$ and is found to be the split real form of the unique complex simple Lie algebra $\g_2^{\C}$ of type $G_2$. 
Often in the literature, this group $\Gtwosplit$ is referred to as $G_{2,2}$, because the rank of $G_2$ is two and the Killing form on $\g_2'$ is of signature $(8,6)$. 

Axiomatically, under the given definition, $\Gtwosplit$ is not assumed to preserve $q_{3,4}$. However, the condition $g \in \Stab_{\GL(7,\R)}(\times_{3,4})$ does imply $g \in \Stab_{\GL(7,\R)}(q_{3,4})\cong O(3,4)$. One way to see that $\times_{3,4}$ secretly encodes $q_{3,4}$ is the following identity: 
$\mathrm{trace}(\mathcal{C}_x|_{x^\bot}) = -\frac{1}{6} q(x) $ (cf. \cite[(13)]{Bae02}), which follows from the double cross-product identity \eqref{DCP}.
Moreover, since $\Gtwosplit$ is connected, one finds the inclusion $\Gtwosplit < \SO_0(3,4)$. 

It is often convenient to describe the $\Gtwosplit$ action on $\R^{3,4}$ via certain Stiefel manifolds. 
The simplest such model is the following Stiefel triplet manifold:  
\[ V_{(+,+,-)}(\R^{3,4}) = \left\{ (u,v,w) \in (\R^{3,4})^3 \mid q(u) = q(v) = -q(w) = +1, \, u\cdot v= u\cdot w = v \cdot w = 0, (u\times v) \cdot w =0. \right\}\] 
The following very well-known result yields a double-fiber bundle description of $\Gtwosplit$ that verifies its connectedness. Additionally, this Stiefel model is a natural analogue of the corresponding model in the compact case (cf. \cite[Remark 5.13]{Fon18} and \cite[Lemma A.8]{HL82}). 
\begin{proposition}[$\Gtwosplit$ Multiplication Stiefel Model]\label{Prop:GtwosplitFirstStiefel}
The group $\Gtwosplit$ acts simply transitively on $V_{(+,+,-)}(\R^{3,4})$. 
\end{proposition}

A proof can be found in \cite[Proposition 2.3]{Eva24}, but the idea is well-known, and is as follows. 
Given a tuple $p = (u,v,w) \in V_{(+,+,-)}(\R^{3,4})$, one constructs an ordered basis $B_p$ for $\R^{3,4}$ by
\[ B_{p} \coloneqq (u,v,u\times v, w, w\times u, w\times v, w\times(u\times v)).\]
Then given $p, p' \in V_{(+,+,-)}(\R^{3,4})$, there is a unique transformation $\psi \in \Gtwosplit$ such that $\psi \cdot p = p'$ and this transformation also is the unique linear map satisfying $\psi \cdot B_p = B_{p'}$. 

Let us now describe the maximal compact subgroup $K \cong_{\mathrm{Lie}} \SO(4)$ of $\Gtwosplit$ as well as a geometric model for the $\Gtwosplit$-symmetric space $\mathbb{X}_{\Gtwosplit}$.
Write 
    \[ \Gr_{(3,0)}^\times(\R^{3,4}) = \{ P \in \Gr_{(3,0)}(\R^{3,4}) \mid P \times_{3,4} P = P \} .\]
We include the following elementary lemma to illustrate the utility of Proposition \ref{Prop:GtwosplitFirstStiefel}.
\begin{lemma}
    Let $P \in \Gr_{(3,0)}^\times(\R^{3,4})$. 
    Then $K_P \coloneqq\Stab_{\Gtwosplit}(P)$ is a maximal compact subgroup of $\Gtwosplit$ isomorphic to $\SO(4)$. 
    Moreover, $\Gtwosplit$ acts transitively on $\Gr_{(3,0)}^\times(\R^{3,4})$.  
\end{lemma}

\begin{proof} 
Let $P \in \Gr_{(3,0)}^\times(\R^{3,4})$. 
We first show $K_P \cong \SO(4)$. 
To this end, consider the map $F: K_P \rightarrow \SO(P^\bot)$ by $\psi \mapsto \psi|_{P^\bot}$ and define the Stiefel manifold
\[ V_{(+,+,-)}(P)\coloneqq \{ (u,v,w) \in V_{(+,+,-)}(\R^{3,4}) \; |\; u,v \in P, w \in P^\bot \}.\] 
Since any $\varphi \in K_P$ preserves the orthogonal splitting $\R^{3,4} = P \oplus P^\bot$, one sees that $K_P$ preserves $V_{(+,+,-)}(P)$. 
In fact, $K_P$ acts simply transitively on $V_{(+,+,-)}(P)$ by Proposition \ref{Prop:GtwosplitFirstStiefel}. 
Indeed, $F$ is injective because $(\psi(u), \psi(v), \psi(w))$ determines $\psi$. Conversely, $F$ is surjective because if if $p=(u,v,w), p'=(x,y,z) \in V_{(+,+,-)}(P)$, there is a unique transformation $\psi \in \Gtwosplit$ such that $\psi \cdot p = p'$ and this map $\psi$ must fix $ \spann \{ u,\,v, \,u\times v\}= P= \spann \{x, \,y, \,x \times y\}$, so $\psi \in K_P$.

Recall the notation $Q_{\epsilon}(V) = \{ x \in V \mid q(x) =\epsilon 1\}$, for $\epsilon \in \{+,-, 0\}$. 
From a slightly different perspective, $K_P$ also acts simply-transitively on the Stiefel manifold 
\[ V_3(P^\bot) = \{ (w_1,w_2,w_3) \in (Q_-(P^\bot))^3 \; | \; w_i \cdot w_j = 0, i \neq j\}, \] 
which just corresponds to looking at the action of $\psi \in K_P$ on $(w, w\times u, w\times v)$, for a fixed tuple $(u,v,w) \in V_{(+,+,-)}(P)$. On the other hand, $\SO(P^\bot) \cong \SO(4)$ acts simply transitively on $V_3(P^\bot)$. We conclude that the restriction map $K_P \rightarrow \SO(P^\bot)$ by $g \mapsto g|_{P^\bot}$ is a Lie group isomorphism. 

Take any points $P_1, P_2 \in \Gr_{(3,0)}^\times(\R^{3,4})$. 
We may write $P_i = \spann_{\R} \{ u_i,\, v_i, \,u_i\times v_i \}$ for some orthonormal elements $u_i, v_i \in Q_+(P_i)$. 
By Proposition \ref{Prop:GtwosplitFirstStiefel}, there is an element $\psi$ of $\Gtwosplit$ such that $\psi \cdot (u_1, v_1) = (u_2,v_2)$ and hence $\psi\cdot P_1 =P_2$. Thus, $K_P$ acts transitively on $\Gr_{(3,0)}^\times(\R^{3,4})$.
\end{proof}

Next, to discuss the flag manifolds of $\Gtwosplit$, it is necessary to discuss annihilators. 
\begin{definition}
Let $u \in \R^{3,4}$. Then the \textbf{annihilator} $\Ann(u)$ of $u$ is given by 
\[ \Ann(u) \coloneqq \{ v \in \R^{3,4} \; | \; u\times_{3,4} v = 0\}.\]
\end{definition}
We take this name ``annihilator'' from the wonderful work \cite{BH14}.\medskip  

We have a complete classification of annihilators, given shortly. First, we recall the \emph{double cross-product} identity \cite[Corollary 3.49]{Kar20} 
\begin{align}\label{DCP}
 u\times (u\times v) = -q(u)v+(u\cdot v)u. 
\end{align}
One consequence of \eqref{DCP} is the following: if $x \in Q_+(\R^{3,4})$, then $J_x\coloneqq \mathcal{C}_x|_{x^\bot}: x^\bot \rightarrow x^\bot$ is an almost-complex structure on $x^\bot$. Making these identifications pointwise, the pseudosphere $\hat{\mathbb{S}}^{2,4} =Q_+(\R^{3,4})$ inherits a non-integrable almost-complex structure, in direct analogue with the non-integrable almost-complex structure on $\mathbb{S}^6$ coming from $\times_{7,0}$ \cite{Leb87}. For study of almost-complex curves in $\hat{\mathbb{S}}^{2,4}$, see \cite{Nie24, Eva24, CT24}. \medskip 

When $u$ is not isotropic, $q(u) \neq 0 $, then $\Ann(u) = \R\{u\}$, as seen by \eqref{DCP}. On the other hand, if  $q(u)=0$, annihilators are non-trivial. 
\begin{proposition}[{\cite[Proposition 7]{BH14}}]\label{Prop:Annihilators}
Let $x \in Q_0(\R^{3,4})$. Then $\Ann(x)$ is a 3-dimensional, isotropic subspace of $\R^{3,4}$. For any spacelike three-plane $P \in \Gr_{(3,0)}^\times(\R^{3,4})$, the subspace $\Ann(u)$ is uniquely the graph of a linear map $\phi_{u,P}: P \rightarrow P^\bot$. Moreover, if $x =u+z$ for $u \in Q_+(P), z\in Q_-(P^\bot)$, then 
\[ \phi_{u,P}(v) = -z \times (u \times v). \] 
\end{proposition}
We remark that any isotropic 3-plane $T \in \Iso_3(\R^{3,4})$ is the graph over any spacelike three-plane $P \in \Gr_{(3,0)}(\R^{3,4})$; the essence of the proposition is to describe the shape annihilator 3-planes take in such a description. 
While the formula obtained in \cite{BH14} for $\phi_{u,P}$ looks slightly different from ours, one can manually verify $(u+z) \times (v-z\times (u\times v)) =0$. Since `annihilator' is a linear condition, we may write $\phi_{l, P}$ for $l \in \Ein^{2,3}$, which is the same as $\phi_{u,P}$ if $l =[u]$. \medskip 

There is one more relevant concept related to annihilators, introduced in \cite{Eva24}.  
\begin{definition}
Let $X = (x_i)_{i=3}^{-3}$ be a basis for $\R^{3,4}$ such that $x_i \times x_j = c_{ij} x_{i+j}$ for some constants $c_{ij} \in \R$. Then we call $X$ an \textbf{$\R$-cross-product basis} for $\R^{3,4}$. 
\end{definition}

A convenient feature of such a basis is that the quadratic form $q_{3,4}$ is automatically ``anti-diagonal'' \cite[Proposition 2.3.6]{Eva24}, so $x_3$ is isotropic, and moreover $\Ann(x_3) = \spann_{\R} \{ x_3, \,x_2,\, x_1 \}$. 
The subalgebra $\mathfrak{a} < \g_2'$ of diagonal transformations in an $\R$-cross-product basis is a maximal split torus of $\g_2'$ of the form
\[ \mathfrak{a} = \{ \mathrm{diag}(r+s, r, s, 0, -s,-r,-r-s) \in \mathfrak{gl}_7\R \; | \; r,s \in \R \}.\]
The associated root system $\Lambda$ of $(\g_2', \mathfrak{a})$ admits a set of simple roots given by $\Delta = \{\alpha,\beta\}$, where $\alpha = r^*-s^*$ and $\beta = s^*$ are choices for the primitive long and short roots, respectively.
Relative to these choices, one may construct the parabolic subgroups $P_{\Delta}, P_{\beta}, P_{\alpha}$ as in Section \ref{Sec:FlagPreliminaries}. 

We can now introduce the three $\Gtwosplit$-flag manifolds. 
The model spaces are as follows: 
\begin{align}
    \Ein^{2,3} &= \mathbb{P}Q_0(\R^{3,4}) = \{ [x] \in \mathbb{P}(\R^{3,4}) \; | \; q(x) = 0 \},\\ 
    \Pho^\times  &= \{ \omega \in \Iso_2(\R^{3,4}) \; |\; \omega \times_{3,4} \omega = 0\}, \\ 
    \mathcal{F}_{1,2}^\times &= \{  \,(\ell, \omega) \in \Ein^{2,3} \times \Pho^\times | \; \ell \subset \omega \}.
\end{align}
We call an element $\omega \in \Pho^\times$ an \emph{annihilator photon} and a pair $(\ell, \omega) \in \mathcal{F}_{1,2}^\times$ a \emph{pointed annihilator photon}. 

Keeping the choices of $X, \alpha, \beta$ fixed as above, define $p_{\beta} \coloneqq \R\{ x_3\} \in \Ein^{2,3}$, $p_{\alpha} = \spann \{ x_3, x_2\} \in \Pho^\times$, $p_{\Delta} \coloneqq (p_{\beta}, p_{\alpha}) \in \mathcal{F}_{1,2}^\times.$ 
In this running notation, we have the following well-known result. 
\begin{proposition}[$\Gtwosplit$-Flag Manifolds]
Let $X$ be an $\R$-cross-product basis for $\R^{3,4}$. Then $P_{\alpha} = \Stab_{\Gtwosplit}(p_{\alpha})$, $P_{\beta} = \Stab_{\Gtwosplit}(p_{\beta})$, $P_{\Delta} = \Stab_{\Gtwosplit}(p_{\Delta})$. The group $\Gtwosplit$ acts transitively on $\mathcal{F}_{1,2}^\times$, $\Pho^\times$, $\Ein^{2,3}$, and there are $\Gtwosplit$-equivariant diffeomorphisms
\[ \begin{cases}
    \Gtwosplit/P_{\beta} \cong \Ein^{2,3} \\
    \Gtwosplit/P_{\alpha} \cong \Pho^\times \\
    \Gtwosplit/P_{\Delta} \cong \mathcal{F}_{1,2}^\times,
\end{cases} \] 
given by $g P_{\sigma} \mapsto g\cdot p_{\sigma}$ for $\sigma \in \{ \alpha, \beta, \Delta \}$. 
\end{proposition}

These models are well-known $\Gtwosplit$ experts \cite{Bry20, MNS21}. 
The proof in the case $\sigma = \Delta$ is contained in \cite[Proposition 2.5.3]{Eva24Thesis}. 
The other cases are similarly handled. 
Here, for the transitivity of the $\Gtwosplit$-action on these flag manifolds, one should not use the Stiefel model from Proposition \ref{Prop:GtwosplitFirstStiefel}, but instead the \emph{null Stiefel triplet model} from \cite[Theorem 12]{BH14}.

Transversality in the $\Gtwosplit$ partial flag manifolds involves only familiar notions from Subsection \ref{Subsec:ConcreteModelsABD}. More precisely, we can explicitly describe ($\Gtwosplit$-)transversality in $\Ein^{2,3}$ and $\Pho^\times$ as follows. Two points $x, x' \in \Ein^{2,3}$ are $\Gtwosplit$-transverse if and only if $q$ defines a non-degenerate pairing between $x$ and $x'$, just as for $G = \SO_0(3,4)$. Similarly, $\omega, \omega'\in \Pho^\times$ are transverse if and only if $q$ defines a non-degenerate pairing between $\omega$ and $\omega'$. 

We now recall the natural embedding of $\Flag(\Gtwosplit)$ into $\Flag(\SO_0(3,4))$. 
\begin{example}\label{Ex:G2B3}
There is an equivariant transversality-preserving map $\iota: \mathcal{F}_{1,2}^\times \hookrightarrow \Iso_{\dbrack{3}}(\R^{3,4})$ by 
\[ (\ell, \omega) \mapsto (\ell, \omega, \Ann(\ell)).\] 
Now, $\Gtwosplit$-transversality and $\SO_0(3,4)$-transversality for $\Ein^{2,3}$ agree. Also, the flag manifold $\Pho^\times $ includes transversely in $\Pho(\R^{3,4})$. Thus, the only remaining transversality consideration involves the 3-planes in the flag. However, here we note 
\[ \ell \pitchfork_{\Ein^{2,3}} \ell' \Longleftrightarrow \Ann(\ell) \pitchfork_{\Iso_3(\R^{3,4})} \Ann(\ell').\] 
In particular, $q|_{\ell + \ell'}$ is non-degenerate if and only if $q|_{\Ann(\ell) + \Ann(\ell')}$ is non-degenerate. This means $\Ann: \Ein^{2,3} \rightarrow \Iso_3(\R^{3,4})$ by $\ell \mapsto \Ann(\ell)$ is transversality-preserving. Hence, $\iota$ is transversality-preserving as well.  
\end{example}

\subsubsection{Pointing towards \texorpdfstring{$\Gtwosplit$}{G2'}-flag manifolds}\label{Sec:PointingTowardsG2'Flags}

We now briefly recall when tangent vectors to a symmetric space point towards certain flag manifolds. We shall only need this notion for $G = \Gtwosplit$, so we focus the discussion entirely on this case. Hence, we set $\mathbb{X} \coloneqq \Gr_{(3,0)}^\times(\R^{3,4})$ for this subsection.\medskip 

Let us fix a basepoint $o \in \X$. Then $K:=\Stab_{\Gtwosplit}(o)$ is a copy of the maximal compact subgroup. The orthogonal Lie algebra  decomposition $\g_2'= \mathfrak{k} \oplus \mathfrak{p}$ is called a Cartan decomposition. Then we fix a maximal abelian subalgebra $\mathfrak{a} \subset \mathfrak{p}$, which is a two-dimensional maximal split torus of $\g_2'$. Let $\overline{\mathfrak{a}}^+$ be a closed Weyl chamber in $\mathfrak{a}$. There is a basis $(x_i)_{i=3}^{-3}$ of a $\R^{3,4}$ such that 
\[ \overline{\mathfrak{a}}^+ = \{ X = \mathrm{diag}(r+s, r,s,0,-s,-r,-r,-s) \in \g_2' \mid r\geq s \geq 0\}.\]
Identify $\mathfrak{p} = T_o\X$, so that $\mathfrak{a} \subset T_o\X$. 
The \textbf{Cartan projection} $\mu: T\X \rightarrow \overline{\mathfrak{a}}^+$ is the map sending $X$ to $gX$, where $gX $ is the unique $G$-translate of $X$ in the Weyl chamber $\overline{\mathfrak{a}}^+$. Let $\alpha$, $\beta$ be the long and short roots of $\Delta(\g_2', \mathfrak{a})$, taken again to be $\alpha = r^*-s^*$ and $\beta = s^*$. 

Finally, following the standard terminology, we say that $0 \neq X \in T\X$ \emph{points towards} $\Ein^{2,3}$ when $\alpha(\mu(X))=0$ and that $X$ \emph{points towards} $\Pho^\times$ when $\beta(\mu(X)) = 0$. \medskip 

We now describe a more geometric way to detect when $X \in T\X$ points towards $\Ein^{2,3}$ or $\Pho^\times$. To this end, a tangent vector $X \in T_P \X$ naturally corresponds to a homomorphism $P \to P^\perp$, denoted abusively by $X$. We may switch between these perspectives. 

First, we observe that $X_{\alpha} = \mathrm{diag}(2,2,0,0,0,-2,-2) \in \overline{\mathfrak{a}}^+$ is the unique element in $\overline{\mathfrak{a}}^+$, up to positive scalars, pointing towards $\Pho^\times$. When viewed as a map $X_{\alpha}:o \rightarrow o^\bot$, then $X_{\alpha}$ has rank two. In fact, $X \in T\X$ points towards $\Pho^\times$ if and only if $\rank(X) = 2$. We have shown one direction, as the $G$-orbit of a rank two map is rank two. To see the converse, note that if $X \neq 0$ and $\alpha(X) = 0$, then $\rank(X) = 3$. 

Next, we note that $X_{\beta} =\mathrm{diag}(2,1,1,0,-1,-1,-2)$ is the unique element in $\overline{\mathfrak{a}}^+$, up to positive scalars, pointing towards $\Ein^{2,3}$. Now, more generally, suppose that $X \in T_{P}\X$ and that $X = \mathrm{diag}(c_i)_{i=3}^{-3}$ under some identification whereby $X \in \mathfrak{a}'\subset \mathfrak{p}'$, where again $\g_2' = \mathfrak{k}'\oplus \mathfrak{p}'$ and $K' = \Stab_{\Gtwosplit}(P)$. Let $(\lambda_1, \lambda_2, 
\lambda _3,0 ,-\lambda_3, -\lambda_2,-\lambda_1)$ be the eigenvalues of $X$ in monotonically decreasing order. Then $X$ points towards $\Ein^{2,3}$ if and only if $\lambda_1 > \lambda_2$ and $\lambda_2 = \lambda_3$.  

We will apply these criteria in Section \ref{Sec:G2Sambarino}. \medskip

\subsection{Fibrations of \texorpdfstring{$\Flag(\Gtwosplit)$}{Flag(G2'} by Maximally Transverse 3-spheres}

In this section, we realize the full flag manifold $\mathcal{F}_{1,2}^\times$ of $\Gtwosplit$ as a principal $\Sp(1)$-bundle over the full flag manifold $\Flag(\R^3)$. The fibers of this bundle are maximally transverse 3-spheres in $\mathcal{F}_{1,2}^\times$. Since the relevant $\Sp(1)\cong \Spin(3)$-subgroup of $\Gtwosplit$ fixes $\R^{3,0}$ and acts on $\R^{0,4}$ as an irreducible  representation, this construction is a special case of the construction in Section \ref{Sec:Spinors}.

In the following theorem, we may denote for any vector subspace $W < \R^{3,4}$, the subset $\Pho^\times(W)\coloneqq \{ \omega \in \Pho^\times \mid \omega \subset W \}$ of $\Pho^\times$.

\begin{theorem}\label{Thm:G2Fibration}
Fix any point $P \in \mathbb{X}_{\Gtwosplit}$. For each $\Gtwosplit$-flag manifold $\mathcal{F}$, there is a principal $\mathsf{Sp}(1)$-bundle fibration of $\mathcal{F}$ naturally associated to $P$ with maximally transverse 3-sphere fibers in $\mathcal{F}$. 
Let $\pi_P \colon \R^{3,4} \rightarrow P$ be the orthogonal projection. The fibrations and fibers are as follows: 
\begin{enumerate}[label=(\alph*)]
    \item\label{item:S3 in Ein23} 
     $\pi_P \colon \Ein^{2,3} \rightarrow \Gr_1(P)$ by $\ell \mapsto \pi_P(\ell)$, with fibers $\Ein^{2,3}|_{U^1} = \Ein(U^1 \oplus P^\bot) \cong \Ein^{0,3}$.
    \item\label{item:S3 in Phox} 
     $\pi_P \colon \Pho^\times \rightarrow \Gr_2(P)$ by $\omega \mapsto \pi_P(\omega)$, with fibers $\Pho^\times|_{U^2} = \Pho^\times(U^2 \oplus P^\bot) \cong \Pho^\times (\R^{2,4})$.
    \item\label{item:S3 in F12x} 
    $\pi_P \colon \mathcal{F}_{1,2}^\times \rightarrow \Flag(P)$ by $(\ell, \omega) \mapsto (\pi_P(\ell), \pi_P(\omega) )$. For $(U^1,U^2) \in \Flag(P)$, the fiber $\mathcal{F}_{1,2}^\times|_{(U^1,U^2)}$ consists of all pointed annihilator photons that are graphs over $(U^1,U^2)$.
\end{enumerate} 
The $\mathsf{Sp}(1)$-action is given by  $H_P\coloneqq\Stab_{\Gtwosplit}^{pt}(P) = \{ \varphi \in \Gtwosplit\; |\; \varphi(x) = x,\; \forall x \in P\}$. 
\end{theorem}

\begin{proof}

\ref{item:S3 in Ein23} 
A copy of $\Ein^{0,3}$ is of the form $\mathcal{S} = \Ein(U^1 \oplus R)$ for subspaces $U^1 \in \Gr_{(1,0)}(\R^{3,4})$ and $R \in \Gr_{(0,4)}(\R^{3,4} \cap \ell^\bot)$. The set $\mathcal{S}$ is seen to be transverse as follows. Fix $x\in Q_+(U^1)$. Any two distinct points in $\Ein(U^1 \oplus R)$ are of the form $[x + y],\; [x + z]$ for $y \neq z \in Q_-(R)$. Then\[ \langle x + y,  x + z\rangle_{q_{3,4}} = +1+ \langle y , z \rangle_{q_{3,4}} > 0, \]
which verifies transversality. $\Ein^{0,3}$ is non-nullhomotopic in $\Ein^{2,3}$ and thus maximally transverse by Fact \ref{Fact:MaximallyTransverse}. Finally, observe that $\mathcal{S}$ is a fiber of the map $\pi_P \colon \Ein^{2,3} \to \Gr_1(P)$ and that $H_P$ preserves the fibers of $\pi_P$. The action of $H_P$ is simply transitive by Proposition \ref{Prop:GtwosplitFirstStiefel}. \medskip

\ref{item:S3 in Phox} Choose any $U=U^2 \in \Gr_2(P)$ and then form $\Pho^\times(U \oplus P^\bot)$. Every photon in $\Pho(\R^{2,4})$ is uniquely the graph over any spacelike two-plane in $\R^{2,4}$. Hence, $\omega \in \Pho^\times(U\oplus P^\bot)$ uniquely obtains the form of a graph of a map $\varphi_{\omega}: U \rightarrow P^\bot$. 
Fix any orthonormal basis $(u,v)$ for $U$ and define $z \coloneqq \varphi_{\omega}(u)$. By Proposition \ref{Prop:Annihilators}, 
\begin{align}\label{PhotonEqn1}
    \varphi_{\omega}(v) =-z\times (u\times v) = (u\times v ) \times z. 
\end{align}
Thus, $\omega = \spann_{\R} \{ u+z, \,v+ (u \times v)\times z\}.$ This means $(U, (u,z))$ determines $\omega$. 
Hence, we obtain a map $\Pho^\times(U\oplus P^\bot) \rightarrow \Ein([u] \oplus P^\bot)$ by $\omega \mapsto [u+ \varphi_{\omega}(u)]$, which is an $H_P$-equivariant bijection. Using that $\Ein([u]\oplus P^\bot)$ is transverse, we will show $\Pho^\times(U\oplus P^\bot)$ is transverse. 

Recall $\mathcal{C}_x|_{x^\bot}^{2} = -\mathrm{id}$ for $x \in Q_+(\R^{3,4})$ by \eqref{DCP}. 
Define $w \coloneqq u \times v$ and note that $\omega \in \Pho^\times( w^\bot)$ is a complex line in $(w^\bot ,\mathcal{C}_w) \cong \C^6$, since $v = w \times u$. Now, the relation $q( w\times x, w \times y) = q(x,y)$ for $x,y \in w^\bot$ implies $(w^{\bot}, \mathcal{C}_{w}, q|_{w^\bot}) \cong \C^{1,2}$. Thus, for any annihilator photons $\omega, \omega' \in \Pho^\times(w^\bot)$, if the pairing between the real lines $\R\{ u+ \varphi_{\omega}(u)\}, \R\{ u+ \varphi_{\omega'}(u)\}$ is non-degenerate, then the pairing between $\omega, \omega'$ is non-degenerate. This implies $\Pho^\times(U \oplus P^\bot)$ is transverse in $\Pho^\times$. 
Finally, we address maximal transversality of $\Pho^\times(U \oplus P^\bot)$. Choose any photon $\omega \in \Pho^\times \backslash \Pho^\times(U \oplus  P^\bot)$. Then $\pi_P(\omega) = U' \neq U \in \Gr_{2}(P)$. We find $\omega' \in \Pho^\times(U \oplus P^\bot)$ such that $\omega'$ is not transverse to $\omega$. To start, the two-planes $U, U'$ must intersect. Take $u \in Q_+(U \cap U')$ and form the line $\ell \coloneqq [u+ \varphi_{\omega}(u)] \subset \omega$. 
By Proposition \ref{Prop:Annihilators}, we may consider the following annihilator photon: $\omega' \coloneqq \graph(\varphi_{\ell, P}|_U) \in \Pho^\times(U \oplus P^\bot)$. That is, we restrict the graph map $\varphi_{\ell,P}$ of the annihilator of $\ell$ with respect to $P$ to $U$ to obtain an annihilator photon. Since $\ell \subset \omega' \cap \omega$, then $q|_{\omega +\omega'}$ is not non-degenerate. This means $\omega \not \pitchfork \omega'$. 

\ref{item:S3 in F12x} The proof of the remaining claims around $\pi_P: \mathcal{F}_{1,2}^\times \rightarrow \Flag(P)$ is essentially contained above. Indeed, the previous arguments show that the projection maps $\pi_{\Ein^{2,3}}: \mathcal{F}_{1,2}^\times \rightarrow \Ein^{2,3}$ and $\pi_{\Pho^\times}: \mathcal{F}_{1,2}^\times \rightarrow \Pho^\times$ are transverse when restricted to the fiber $\fullflag|_{(U^1,U^2)}$. This shows that $\fullflag|_{(U^1,U^2)}$ is transverse and then maximally transverse as well by the maximal transversality of $\Pho^\times(U^2\oplus P^\bot)$ or of $\Ein(U^1\oplus P^\bot)$. The group $H_P$ acts simply transitively on $\fullflag|_{(U^1,U^2)}$ by similar considerations as in \ref{item:S3 in Ein23}, \ref{item:S3 in Phox}. 
\end{proof}

In fact, the 3-spheres in $\mathcal{F}_{1,2}^\times$ are not only maximally transverse in each $\Gtwosplit$-flag manifold, but also in all the $B_3$-flag manifolds; see Corollary \ref{Cor:G2SphereMaxTransversality}.

\subsection{Non-surface Anosov subgroups of \texorpdfstring{$\Gtwosplit$}{G2'}}\label{Sec:G2Sambarino}

In this section, we demonstrate the existence of $\{\alpha\}$-Anosov and $\{\beta\}$-Anosov subgroups of $\Gtwosplit$ that are not virtually free or surface groups, obtaining a negative answer to the generalized form of Sabmarino's question for $G=\Gtwosplit$ and $\Theta =\{\beta\}, \Theta = \{\alpha\}$. 
The examples are constructed by finding appropriate real rank one subgroups of $\Gtwosplit$.

\subsubsection{For \texorpdfstring{$\{\alpha\}$}{alpha}-Anosov Representations}

The first counterexample comes from studying the stabilizer of a spacelike vector $x \in Q_+(\R^{3,4})$. This $\SU(1,2)$-subgroup is classically known \cite{Yok77}. Here, note that the 3-sphere in the following proposition is exactly one of the 3-spheres from Theorem \ref{Thm:G2Fibration} part (b). 

\begin{proposition}
    Fix any $P \in \mathbb{X}_{\Gtwosplit}$ and $x \in Q_+(P) $.
    Define $H_x\coloneqq \Stab_{\Gtwosplit}(x) \cong \mathsf{SU}(1,2)$. Then the sub-symmetric space $\mathbb{X}_{H_x} \cong_{\mathbf{Isom}} \mathbb{CH}^2$ has its visual boundary $\visualboundary \mathbb{X}_{H_x}$ identify in $\visualboundary\mathbb{X}_{\Gtwosplit}$ as the transverse $\mathbb{S}^3$ given by $\Pho^\times( x^\bot) \cong \Pho^\times(\R^{2,4})$.
\end{proposition}

\begin{proof}
We first recall that a tangent vector $X \in T_{P}\Gtwosymmetricspace$ points towards $\Pho^\times$ exactly when $X$ is a rank two map $X \colon P \rightarrow P^\bot$, as described in Section \ref{Sec:PointingTowardsG2'Flags}. We claim that, up to positive scalars, such a map must obtain the following form: 
\begin{align}\label{PointingTowardsPhox}
\begin{cases} 
    X(u) &= z \\ 
    X(v) &= (u\times v) \times z \\
    X(u \times v) &= 0, \end{cases} 
\end{align}
where $u,v \in P$ are orthonormal vectors in $P$ and $z \in Q_-(P^\bot)$. Define $w\coloneqq u\times v$, so $(u,v,w)$ is an oriented orthonormal basis for $P$. A tangent vector $X$ of the form \eqref{PointingTowardsPhox} has rank two and thus points towards $\Pho^\times$. In fact, $X$ points towards the annihilator photon $\omega \coloneqq \graph(X) = \spann \{ u+z, v+wz \}$. One then notes by \eqref{PointingTowardsPhox} that $X$ is equivalently the graph of the unique $\mathcal{C}_w$-holomorphic map $X: \spann \{ u,v\} \rightarrow P^\bot$ such that $X(u) = z$. 

We now show if $X: P \rightarrow P^\bot$ is any rank two map of operator norm one that is a derivation of $\times_{3,4}$, then $X$ obtains the form \eqref{PointingTowardsPhox}. Indeed, choose any $u \notin \ker(X)$ and write $X(u) = z$. Also, choose $w \in Q_+(\ker(X))$. The derivation condition on $X$ forces  
\[ X(w \times u) = X(w) \times u + w \times X(u) = w\times X(u).\]
The normalization condition (iv) on $q$ implies $q(w \times X(u) ) = q(w)q(X(u)) = q(X(u))$. In particular, the ratio $\frac{ |q(X(y))|}{q(y)}$ is independent of $y \in \spann_{\R} \{ u,\,w\times u \}$. Hence, $q(X(y))$ is independent of $y \in Q_+(P)$. By hypothesis, $X$ has operator norm $||X||_{\rm{op}} = 1$, so we must have $|q(z)| = 1$. We conclude that $X(u) =z, \,X(v) = w\times z$, meaning $X$ obtains the form \eqref{PointingTowardsPhox}.

Now, take any $x \in Q_+(P)$. Then the cross-product endomorphism $\mathcal{C}_x: \R^{3,4} \rightarrow \R^{3,4}$ by $\mathcal{C}_x(y) = x \times y$ satisfies $\mathcal{C}_x^{^\circ 2} |_{x^\bot} = -\id_{x^\bot} $. Hence, identify $(x^\bot, \mathcal{C}_{x}) \cong \mathbb{C}^3$ as a complex vector space. Now, the quadratic form $q$ on $\R^{3,4}$ reinterprets as a signature (1,2)-hermitian form on $x^\bot$. Thus, identify $(x^\bot, \mathcal{C}_x, q) \cong \C^{1,2}$. One can then swiftly show the map $\Stab_{\Gtwosplit}(x) \rightarrow \SU(x^\bot, \mathcal{C}_x, q)$ by $g \mapsto g|_{x^\bot}$ is a Lie group isomorphism by dimension count and the fact that $\SU(1,2)$ is a connected simple Lie group.

We can now assemble the pieces. Set $H_x: = \Stab_{\Gtwosplit}(x)$. Here, we may include $\mathbb{X}_{H_x} \hookrightarrow \mathbb{X}_{\Gtwosplit}$ as the $H_x$-orbit of the model point $P \in \Gtwosymmetricspace$. Now, take any tangent vector $X \in T_P\mathbb{X}_{H_x}$. By construction, $x \in P$. Hence, we may write $P = \spann \{ x, u, v\}$ for some orthonormal basis $(u,v,x)$ of $P$. As $H_x$ fixes $x$, then viewed as a map $X: P \rightarrow P^\bot$, we have $\rank(X) \leq 2$. On the other hand, if $\rank(X) \leq 1$, then $X =0$, which follows immediately from the derivation condition. Hence, any nonzero tangent vector $X \in T_P\mathbb{X}_{H_x}$ has rank two and thus points towards $\Pho^\times$. Finally, as $X(x) =0$, we see $X$ must be $\mathcal{C}_x$-holomorphic and hence $X$ points towards $\omega \in \Pho^\times(U \oplus P^\bot)$ for $U = x^\bot \cap P$. 
\end{proof}

We now conclude a negative answer to the generalized form of Sambarino's question for $\{\alpha\}$-Anosov representations. 
\begin{corollary}
    Let $\Gamma < H_x =\Stab_{\Gtwosplit}(x) \cong \mathsf{SU}(1,2)$ be a uniform lattice. 
    The inclusion $\iota \colon \Gamma \hookrightarrow \Gtwosplit$ is $\{\alpha\}$-Anosov, and the boundary map $\xi^2 \colon \partial \Gamma \cong \mathbb{S}^3 \rightarrow \Pho^\times $ has image $\Pho^\times( x^\bot)$.
\end{corollary}

\subsubsection{For \texorpdfstring{$\{\beta\}$}{beta}-Anosov Representations}

The negative example for $\{\beta\}$-Anosov representations comes from the $\SL(2,\R)$-subgroup associated to the short root. 
We recall the relevant Lie theory from \cite{DavEva25}. 
Define $H_{\beta}$ to be the analytic subgroup of $\Gtwosplit$ with Lie algebra $\mathfrak{h}_{\beta} = \spann \{ t_{\beta}, e_{\beta}, e_{-\beta} \}$. 
To describe $H_{\beta}$, fix a background $\R$-cross-product basis $(x_i)_{i=3}^{-3}$ and recall the $\alpha$-height grading of $\R^{3,4}$ as follows: 
\[ \R^{3,4} = \spann \{ x_3, x_2 \} \oplus \spann \{ x_1, x_0, x_{-1} \} \oplus \spann \{ x_{-2}, x_{-3} \} =: V_{\alpha,1} \oplus V_{\alpha, 0} \oplus V_{\alpha, -1}. \]
The Levi subgroup $L_{\alpha} < P_{\alpha}$ is precisely the subgroup of $\Gtwosplit$ preserving the $\alpha$-height splitting. Now, $L_{\alpha} \cong \GL(2,\R)$, with the restrictions to $V_{\alpha, \pm 1}$ each being Lie group isomorphisms. Finally, $H_{\beta } <L_{\alpha}$ is the subgroup preserving the volume form on the subspaces $V_{\alpha, \pm1}$. 
\begin{proposition}
The sub-symmetric space $\mathbb{H}^2_{\beta} \hookrightarrow \mathbb{X}_{\Gtwosplit}$ associated to $H_{\beta}$ has its visual boundary $\visualboundary\Ha^2_{\beta}$ identify in $\visualboundary\mathbb{X}_{\Gtwosplit}$ as a copy of $\Ein^{0,1}$. 
\end{proposition}

\begin{proof}
Fix a background $\R$-cross-product basis $(x_i)_{i=3}^{-3}$ for $\R^{3,4}$. We use the model basepoint 
\[p_0 = \spann \{ x_3 +x_{-3}, x_2+x_{-2}, x_1+x_{-1}\} \in \mathbb{X}_{\Gtwosplit}.\]
Consider the transformation $g_s = \exp(s\,t_{\beta})$, for $t_{\beta} = \mathsf{diag}(1,-1,2,0,-2,1,-1)$ the co-root of $\beta$ in the basis $(x_i)_{i=3}^{-3}$. Then 
\[g_s \cdot p_0 = \spann \{ e^s x_3 +e^{-s}x_{-3}, e^{-s} x_2+e^s x_{-2}, e^{2s} x_1+e^{-2s}x_{-1} \} .\]
Note that $X \coloneqq \frac{d}{ds}\big|_{s=0}(g_s\cdot p_0 )\in T_{p_0}\mathbb{X}_{\Gtwosplit}$, viewed as a map $X: p_0 \rightarrow p_0^\bot$, is a non-vanishing endomorphism. Since $\alpha(\mu(X)) = 0$, the tangent vector $X$ points towards $\Ein^{2,3}$ by the criteria in Subsection \ref{Sec:PointingTowardsG2'Flags}.  
Moreover, $X$ points towards the leading eigenline of $g_s$, which is $[x_1]$.
Thus, $[x_1] \in \visualboundary\mathbb{H}^2_\beta$. The embedding $\visualboundary\Ha^2_\beta \hookrightarrow \visualboundary\mathbb{X}_{\Gtwosplit}$ is an orbit map, which implies $\Ein^{0,1} \cong \Ein(V_{\alpha,0}) = H_{\beta} \cdot [x_1] \subset \visualboundary\mathbb{H}^2_\beta$. As $\visualboundary\mathbb{H}^2_\beta$ is a topological circle, this completes the proof. 
\end{proof}

We now conclude that there are $\{\beta\}$-Anosov representations of hyperbolic groups $\Gamma$ that are not virtually free or surface groups. 
\begin{corollary}\label{Cor:13AnosovCounterExamples}
    Let $\Gamma$ be a surface group. 
    Then there exists 
    \begin{enumerate}[label=(\alph*)]
        \item \label{item13Anosov_a} $\{\beta\}$-Anosov representations $\Gamma *\mathbb{Z} \rightarrow \Gtwosplit$.
        \item \label{item13Anosov_b} $\{1,3\}$-Anosov representations $\Gamma *\mathbb{Z} \rightarrow \SL(7,\R)$.
    \end{enumerate}
\end{corollary}

The corollary follows from \cite[Corollary 6.3]{DK23}, using the same reasoning from \cite[Proposition 5.2]{KT24}. Here is a quick sketch of the argument. To start, select $\{x^+, x^-\} \subset \Ein^{2,3}$ such that $x^+ \pitchfork x^-$ and $x^{\pm} \pitchfork \Ein(V_{\alpha,0})$, using that $\Ein(V_{\alpha, 0}) $ is contained in a transverse 3-sphere $\Ein^{0,3}$.  
Then there is a loxodromic element $\delta \in \Gtwosplit$ such that $\delta$ acts on $\Ein^{2,3}$ with attracting fixed point $x_+$ and repelling fixed point $x_-$. More precisely, one can ask the subgroup $\Z \cong\langle \delta \rangle $ to be $\{\beta\}$-Anosov with limit set $\{x_+, x_-\}$. Now, define $\Gamma'$ to be a Fuchsian, and hence $\{\beta\}$-Anosov, subgroup of $H_{\beta}$. One then finds a finite index subgroup $\Gamma < \Gamma'$ such that the free product $\Gamma * \Z$ is also $\{\beta\}$-Anosov by \cite[Corollary 6.3]{DK23}. As the subgroup $\Gamma$ is also a surface group, this proves \ref{item13Anosov_a}. Then \ref{item13Anosov_b} follows from studying the inclusion $\Gtwosplit \hookrightarrow \SL(7,\R)$. 

\begin{remark}
The inclusion $\Gamma \hookrightarrow H_{\beta} < \Gtwosplit$ of a Fuchsian subgroup $\Gamma < H_{\beta}$ is $(3,3,6)$ hyperconvex in the sense of \cite[Definition 6.1]{PSW21}. Indeed, for $x,y ,z \in \Ein(V_{\alpha, 0}) \cong \partial \Gamma$ pairwise distinct, one finds 
\[ (x^3+y^3)+z^1 = (x^1+y^1+V_{\alpha, 1} + V_{\alpha, -1})+z^1 = (x^1+y^1+z^1)+V_{\alpha,1} + V_{\alpha,-1} = \R^{3,4}.\]
Thus, while the $\{\beta\}$-Anosov representation $\Gamma * \Z \rightarrow \Gtwosplit$ from \ref{item13Anosov_a} is $\{1,3\}$-Anosov in $\SL_7(\R)$, this representation cannot be (3,3,6)-hyperconvex by \cite[Theorem 1.3]{PT24}. 
\end{remark}

\section{Transverse Spheres via Normed Division Algebras}\label{Sec:TransverseSpheresDivisionAlgebras}

In this section, we construct maximally transverse spheres in type $B_1, B_3, B_7, D_2, D_4, D_8$ full flag manifolds via \emph{normed division algebras} and their more general relatives \emph{composition algebras}. In particular, this construction is a generalization of the 3-sphere construction in $\Flag(\Gtwosplit)$ from Theorem \ref{Thm:G2Fibration} and a special case of the construction from Section \ref{Sec:Spinors}. 
We then show the resulting 1,3,7-spheres are not only maximally transverse in the relevant $B$, $D$ full flag manifolds, but also in every possible partial flag manifold.

\subsection{Composition Algebras \& Normed Division Algebras}\label{Subsec:CompositionAlgebras}
In this section, we recall some details on (real) normed division algebras, composition algebras, and the Cayley-Dickson process in preparation for the construction of maximally transverse spheres. This material is very well-known and can be referenced in \cite{HL82, Sch95, Bae02, Fon18, Kar20}. 

Let us now recall \emph{Hurwitz's theorem} on the classification of normed division algebras over $\R$.
\begin{definition} 
Let $\mathbb{A} = (V, \odot, q)$ be an $\R$-algebra with underlying vector space $V$, algebra product $\odot$, and quadratic form $q$. Then $\mathbb{A}$ is a \textbf{normed division algebra}, when $\mathbb{A}$ is a (real) unital algebra such that $q(x \odot y) = q(x)q(y)$, and $q$ is Euclidean. More generally, a real unital algebra $\mathbb{A} = (V, \odot ,q)$ such that $q(x \odot y) = q(x)q(y)$, with $q$ non-degenerate, is called a \textbf{composition algebra}. 
\end{definition}

Hurwitz's theorem asserts that the list of all normed division algebras is as follows: $\{ \R, \C ,\Ha, \Oct\}$ \cite{Hur98}.\footnote{There are very closely related results. A \emph{division algebra} $A$ is an algebra where $ax=b$ is uniquely solvable for any $a,b \neq 0\in A$. Over $\R$, the only alternative division algebras are $\R, \C, \Ha, \Oct$, and the only associative division algebras are $\R,\C, \Ha$. The latter result is \emph{Frobenius' theorem} \cite{Fro78} and the former result was proven by Zorn \cite{Zor31}.}  
In the larger family of composition algebras, there are three more members, namely the \emph{split} counterparts of the aforementioned algebras: $\C', \Ha', \Oct'$. This more general result is also referred to as Hurwitz's theorem.
We refer the reader to \cite[Remark 3.7]{Fon18}, which also explains Proposition \ref{Prop:CrossProductCompAlg}, namely, how composition algebras relate to cross-products. 
Recall that if $(V,q)$ a vector space with quadratic form $q$, then a \emph{cross-product} $\times:V \times V \rightarrow V$ on $(V,q)$ is a bilinear, alternating, $q$-orthogonal map satisfying $q(u\times v) = q(u)q(v)-q(u,v)^2$.  

\begin{proposition}[Cross-Products \& Composition Algebras] \label{Prop:CrossProductCompAlg}
Suppose $(V, q, \times_V)$ is a normed vector space with cross-product. Then $\mathbb{A} = (V \oplus \R, \odot, q_{\mathbb{A}})$ is a composition algebra, where $q_{\mathbb{A}}:= q_V \oplus q_{\R}$, $1 \odot u = u $, and $u\odot v \coloneqq u\times v-q(u,v)1$. 
\end{proposition}
In the converse direction, given a composition algebra $\mathbb{A}= (V, \odot, q)$, then define $\mathrm{Im}(\mathbb{A})$ as the orthogonal complement of $1 \in \mathbb{A}$ and form the linear map $\wedge: \mathrm{Im}(\mathbb{A}) \times \mathrm{Im}(\mathbb{A}) \rightarrow \mathrm{Im}(\mathbb{A})$ by $u \wedge v \coloneqq \mathrm{Im}(u\odot v)$. Then $\wedge$ is a cross-product on $(\mathrm{Im}(\mathbb{A}), q|_{\mathrm{Im}(\mathbb{A})})$. In this way, Hurwitz's theorem leads to a classification of the possible signatures in which a cross-product is possible, with $\Oct'$ corresponding to $\times_{3,4}$ on $\R^{3,4}$. See \cite[Section 3.3]{Kar20} for further details, focused on the case of normed division algebras. 

Next, we briefly recall the Cayley-Dickson construction, which builds all composition algebras from $\R$. To start, let $(\mathcal{A}, *_{\mathcal{A}}, q_{\mathcal{A}} )$ be a real unital algebra $\mathcal{A}$ with involution $*_{\mathcal{A}}$ satisfying $(xy)^* =y^*x^*$ and $\R \{1_{\mathcal{A}} \} = \mathrm{Fix}(*_{\mathcal{A}})$, and quadratic form $q_{\mathcal{A}} = xx^*$. The Cayley-Dickson and split Cayley-Dickson processes $\mathrm{CD}, \mathrm{CD}'$ produce a new real unital algebra $(\mathcal{B}, *_{\mathcal{B}},q_{\mathcal{B}})$ as follows. The vector space underlying $\mathcal{B}$ is $\mathcal{A} \oplus \mathcal{A}$, with multiplication defined by 
\begin{align}\label{CayleyDickson}
 (a,b)(c,d) \coloneqq (ac + \epsilon \, db^*, a^*d + cb) , 
\end{align}
where $\epsilon = -1$ is used for $\mathrm{CD}$ and $\epsilon = +1$ for $\mathrm{CD}'$. Then $q_{\mathcal{B}} \coloneqq xx^*$, and $*_{\mathcal{B}}=(\mathrm{*}_{\mathcal{A}}, -\id).$
This procedure amounts to an algebra extension of $\mathcal{A}$ by a new element $\ell$, corresponding to $ (0,1) \in \mathcal{A} \oplus \mathcal{A}$, such that  $\ell^2 = \epsilon$.\footnote{The product in \eqref{CayleyDickson} of $(a,b)$ and $(c,d)$ defines multiplication of $a+\ell b$ and $c+\ell d$. The convention $(a,b) \leftrightarrow a+b\ell$ yields a different multiplication formula. Our convention matches that of \cite{BH14, Sch95} and differs from \cite{HL82}.} In particular, $\mathcal{A} \cong \mathcal{A} \oplus 0$ is a subalgebra of $\mathcal{B}$ and $1_{\mathcal{B}}= (1,0)$. 
One finds the new quadratic form is given by $q_{\mathcal{B}} = q_{\mathcal{A}} \oplus (-\epsilon  q_{\mathcal{A}})$. 

We recall some well-known general properties of these processes.
Note that we do not assume $\mathcal{A}$ is a composition algebra below. The proof for $\CD$ is in \cite[Lemma A.11]{HL82}, while small modifications handle the case of $\CD'$. See \cite[page 24]{Sch95} for proof of (c). 
\begin{proposition}\label{prop:cayley dickson}
Let $(\mathcal{A}, q_{\mathcal{A}}, *_{\mathcal{A}})$ be a real unital algebra and set $\mathcal{B}\coloneqq\CD(\mathcal{A})$,  $\mathcal{B}'\coloneqq \CD'(\mathcal{A})$. Then $(\mathcal{B}, \mathcal{*}_{\mathcal{B}}, q_{\mathcal{B}})$ is a real unital algebra and 
\begin{enumerate}[label=(\alph*)]
    \item $\mathcal{B},\mathcal{B}'$ are commutative if and only if $\mathcal{A} = \R$. 
    \item $\mathcal{B},\mathcal{B}'$ are associative if and only if $\mathcal{A}$ is commutative and associative. 
    \item $\mathcal{B},\mathcal{B}'$ are alternative if and only if $\mathcal{A}$ is associative. 
 \end{enumerate}
\end{proposition}
Recall that an algebra $\mathcal{A}$ is \emph{alternative} when the subalgebra generated by two elements $x,y \in \mathcal{A}$ is associative.\footnote{Traditionally,  $\mathcal{A}$ is called alternative when the identities $x(xy)= x^2y$ and $(yx)x= yx^2$ hold for all $x,y \in \mathcal{A}$ \cite{Sch95}. It is a non-trivial theorem of Artin that this notion of alternative coincides with the definition given above. See \cite[Section III]{Sch95} for further details and exposition.}
Proposition \ref{prop:cayley dickson} clarifies how the Cayley-Dickson processes erode structure. 
In particular, applying either Cayley-Dickson process to one of the alternative algebras $\Oct$ or $\Oct'$, we do not obtain an alternative algebra. We emphasize that the alternative property leads to useful identities like $L_x^2 = L_{x^2}$ and is exactly why $\Oct, \Oct'$ are tractable. 

The Cayley-Dickson processes build all composition algebras from $\R$. More specifically, one can realize each composition algebra via the following sequences.
\begin{itemize}
    \item $\R \stackrel{\CD}{\longrightarrow}\C \stackrel{\CD}{\longrightarrow}\Ha \stackrel{\CD}{\longrightarrow}\Oct $.
    \item $\R \stackrel{\CD'}{\longrightarrow}\C'$.
    \item $\R \stackrel{\CD}{\longrightarrow}\C \stackrel{\CD'}{\longrightarrow}\Ha'$.
    \item $\R \stackrel{\CD}{\longrightarrow}\C \stackrel{\CD}{\longrightarrow}\Ha \stackrel{\CD'}{\longrightarrow}\Oct' $.
\end{itemize}
Note that there may be more than one realization of a composition algebra by such a sequence. For example, $\Ha'$ is realized also by $\R \stackrel{\CD'}{\longrightarrow}\C '\stackrel{\CD}{\longrightarrow}\Ha'$. As noted earlier, some Clifford algebras are actually composition algebras. In particular, we have the following $\R$-algebra isomorphisms: $\Cl(\R^{1,0}) \cong \C$, $\Cl(\R^{2,0})\cong \Ha$, $\Cl(\R^{0,1}) \cong \C' \cong \R \oplus \R$, $\Cl(\R^{1,1}) \cong \Cl(\R^{0,2} ) \cong \Ha' \cong \Mat_2(\R)$. 

\subsection{Constructing the Spheres}
We now use the structure of normed division algebras to construct transverse spheres in the $B_1, D_2$, $B_3, D_4$, $B_7, D_8$ full flag manifolds. 
After the proof, we remark on a cousin of the construction with composition algebras, which leads back to the 3-spheres in $\Flag(\Gtwosplit)$ from Theorem \ref{Thm:G2Fibration}. 

\begin{theorem}\label{Thm:DivisionAlgebras}
Let $\mathbb{A} \in \{\C, \Ha, \Oct\}$. Consider the map $\rho: Q_+(\mathbb{A}) \rightarrow \SO(\mathbb{A})$ given by 
$x \mapsto L_x$ and write $a = \dim_{\R}(\mathbb{A})$. Let $\epsilon \in \{-1,0\}$. Define $\iota: \SO(\mathbb{A}) \hookrightarrow \SO_0(a+\epsilon,a)$ by $g\mapsto \id_{\R^{a+\epsilon}} \oplus g$. Then for any flag $F \in \Flag(\R^{a+\epsilon,a})$, the orbit $\image(\iota \circ \rho)\cdot F$ is a maximally transverse $(a-1)$-sphere in the full flag manifold $\Flag(\R^{a+\epsilon,a})$. 
\end{theorem}

\begin{proof}
If $\mathbb{A} = \mathbb{C}$, then the argument is quick. In this case, $\Iso_1(\R^{2+\epsilon, 2}) = \Ein^{1+\epsilon,1}$ and $Q_+(\C)$ acts freely on $\Ein^{1+\epsilon,1}$ and each orbit is a copy of $\Ein^{0,1}$, which is a maximally transverse circle. 

Now, suppose $\mathbb{A} \in \{ \Ha, \Oct\}$. For notational convenience, we shall work with $\epsilon = -1$, namely the $B_{a-1}$ case. The same arguments work for the $D_a$ case. 

Fix a point $P \in \Gr_{(a-1,0)}(\R^{a-1,a})$. Choose any flag $F \in \Iso_{\dbrack{a-1}}(\R^{a-1,a})$. 
Fix any orthonormal basis $U= (u_i)_{i=1}^{a-1}$ of $P$ 
and choose an orthonormal tuple $(z_i)_{i=1}^{a-1}$ of $P^\bot$ such that $F^i  = \spann_{\R} \{ u_j + z_j \}_{j=1}^i$. For $x \in Q_+(\mathbb{A})$, write $Z_i(x) = u_i + xz_i$. Then the flag $x \cdot F$ has component subspaces $x \cdot F^i = \spann \{ u_j+xz_j \}_{j=1}^i$. We need to prove that for $x \neq y \in Q_+(\mathbb{A})$ and any indices $1 \leq i \leq a-1$ that the $i$-planes $x \cdot F^i$ and $y \cdot F^i$ are transverse. The remaining argument is a variation on the proof from Theorem \ref{Thm:TransverseSpinorSphere}. 

We now wield the algebraic gadgetry that $\mathbb{A}$ carries: the quadratic form $ q= q_{\mathbb{A}}$ that is multiplicative over the algebra product. 
In particular, by polarizing the identity $q(vw)=q(v)q(w)$, as in \cite[Appendix IV.A]{HL82}, one finds that for any $x, y \in Q_+(\mathbb{A})$ and $w \in Q_+(\mathbb{A})$, we have 
\begin{align}\label{DivisionAlgebraSymmetry}\langle xw, yw\rangle_{\mathbb{A}} = \langle x, y \rangle_{\mathbb{A}}.
\end{align}
Moreover, we recall that for $x \in Q_+(\mathbb{A})$, we have $L_x^* = L_{x^*}$ and $R_{x}^* = R_{x^*}$ \cite[(A.2)]{HL82}. Here, $x^*$ is the $\mathbb{A}$-conjugate of $x$ and $g^*$ denotes the $q_{\mathbb{A}}$-adjoint. Now, we suppose that $x,y \in Q_+(\mathbb{A})$ are arbitrary and $u, w \in Q_+(\mathbb{A})$ are orthogonal. Then we see that 
\begin{align}\label{DivisionAlgebraAntiSymmetry}
\langle x u, yw\rangle = \langle (xu)w^*, y \rangle =_{(\star)} -\langle (xw)u^*, y\rangle = -\langle xw, yu\rangle.
\end{align}
Here, $(\star)$ follows from the general fact that $(\alpha \beta)\gamma^* = -(\alpha \gamma)\beta^*$ for any $\alpha \in \mathbb{A}$ and $\beta, \gamma \in \mathbb{A}$ that are orthogonal \cite[(A.7)]{HL82}. 

We have all we need to prove the transversality. Let $1 \leq m \leq a-1$ be any index. Assemble the matrix $A_m = (Z_j(x) \cdot Z_k(y))_{j, k= 1}^m$. We will show that $\det(A_m) \neq 0$, which proves that $x \cdot F^m$ and $y \cdot F^m$ are transverse. 
Write $A_m = (a_{ij})$. Then we see that 
\[ a_{ii} = \langle u_i+xz_i, u_i+yz_i\rangle_{\R^{a,a}} = 1-\langle xz_i, yz_i \rangle_{\mathbb{A}} = 1-\langle x, y\rangle_{\mathbb{A}}. \] 
Thus, $a_{ii} = a_{11} > 0$ for all $1\leq i \leq m$. On the other hand, 
for indices $i \neq j$, we see that 
\[ a_{ij} = \langle u_i + xz_i, u_j + yz_j\rangle_{\R^{a,a}} = \langle xz_i ,\,yz_j\rangle_{\mathbb{A}} =_{\eqref{DivisionAlgebraAntiSymmetry}} - \langle yz_i, \,xz_j\rangle_{\mathbb{A}} = -a_{ji}.\]
Hence, $A_m = a_{11} \id + B$, where $B$ is skew-symmetric. By Lemma \ref{Lem:ElementaryLinAlg}, $\det(A_m) \neq 0$. We conclude that $x \cdot F^m$ is transverse to $y \cdot F^m$. Thus, $x \cdot F$ is transverse to $y \cdot F$. 

The maximal transversality follows from the fact that $\image(\iota \circ \rho) \cdot F^1$ is a copy of $\Ein^{0,a-1}$, which is maximally transverse in $\Ein^{a-2,a-1}$ by Fact \ref{Fact:MaximallyTransverse}.
\end{proof}

\begin{remark}\label{Remk:CompositionAlgebra}
If one considers the split counterparts $\C', \Ha', \Oct'$, which are, as normed vector spaces, $\C'\cong \R^{1,1}, \Ha' \cong \R^{2,2}, \Oct' \cong \R^{4,4}$, then there is a natural variation of the construction in Theorem \ref{Thm:DivisionAlgebras}. Write in order $(\C', \Ha', \Oct') = (\mathbb{A}'_n)_{n=2}^4$ and $(\R, \C, \Ha, \Oct) = (\mathbb{A}_n)_{n=1}^4$. Consider the canonical splitting $\mathbb{A}_n' =\mathbb{A}_{n-1} \oplus \mathbb{A}_{n-1}^\bot$ coming from split Cayley-Dickson. As in Theorem \ref{Thm:G2Fibration}, we can consider the subgroup $H_n \coloneqq \Stab_{\mathcal{G}_n'}^{pt}(\mathbb{A}_{n-1})$ for $\mathcal{G}_n'\coloneqq \Aut_{\R-\mathrm{alg}}(\mathbb{A}_n')$. Define again $a_n = \dim_{\R}(\mathbb{A}_n)$. In particular, $H_n \cong_{\mathrm{Lie}} \mathbb{S}^{a_n-1}$ for $n \in \{2,3,4\}$. One finds that $H_n$ acts freely on $\Flag(\R^{a_n,a_n})$ and with orbits that are maximally transverse $(a_n-1)$-spheres, for $n \geq 3$. In the case of $n =4$, applied to $\R^{3,4}$ rather than $\R^{4,4}$, this construction describes the fibration of $\Flag(\Gtwosplit)$ by maximally transverse 3-spheres in Theorem \ref{Thm:G2Fibration}. 
\end{remark}

In fact, the transverse spheres in Theorem \ref{Thm:DivisionAlgebras} are not only maximally transverse in the full flag manifold, but also in every partial flag manifold. 
\begin{corollary}\label{Cor:CompAlgMaxTransverlity}
Let $a \in \{2,4,8\}$ and $\epsilon \in \{0,-1\}$ and $\mathrm{pr}_{\Theta}: \Flag(\R^{a-\epsilon,a}) \rightarrow \Iso_{\Theta}(\R^{a+\epsilon,a})$ be the natural projection to any self-opposite $\SO_0(a+\epsilon, a)$-flag manifold. 
If $\Lambda \subset \Flag(\R^{a+\epsilon,a})$ is a transverse $(a-1)$-sphere from Theorem \ref{Thm:DivisionAlgebras}, then $\mathrm{pr}_{\Theta}(\Lambda)$ is maximally transverse.
\end{corollary}

\begin{proof}
We first note a general fact. Suppose $\Theta' \subset \Theta \subset \Delta(\g, \mathfrak{a})$. If $\mathcal{S}' \subset \mathcal{F}_{\Theta'}$ is maximally transverse and $\mathcal{S} \subset \mathcal{F}_{\Theta}$ is transverse and $\mathcal{S}$ projects to $\mathcal{S}'$ in $\mathcal{F}_{\Theta'}$, then $\mathcal{S}$ is maximally transverse in $\mathcal{F}_{\Theta}$. By this general fact, it suffices to prove the claim for all of the partial flag manifolds of the maximal parabolics, namely in $\Iso_{i}(\R^{a+\epsilon,a})$ and $\Iso_{a}^{\pm}(\R^{a,a})$. \medskip 

The claim for $i=1$ was addressed already by the maximal transversality of $\Ein^{0,a-1}= \mathrm{pr}_1(\Lambda)$. 
Thus, we may suppose $a > 2$. For the rest of the proof, we handle the $B_{a-1}$-case, but the same proof applies to handle $D_a$, including $\Iso_a^+(\R^{a,a})$. \medskip 

Now, suppose $2 \leq i \leq a-1$. 
Choose $T^i \in \Iso_i(\R^{a-1,a}) \backslash \mathrm{pr}_i(\Lambda)$. Recall the map $\iota \circ \rho: Q_+(\mathbb{A}) \rightarrow \SO_0(a-1,a)$. Write $\Lambda = \mathrm{image}(\iota\circ \rho) \cdot F$ for some flag $F \in \Iso_{\dbrack{a-1}}(\R^{a-1,a})|_f$, i.e., $F$ is a graph over $f \in \Flag(\R^{a-1})$. Write as well $f= (f^i)_{i=1}^{a-1}$. We consider two cases. Let $P \coloneqq \R^{a-1,0}$ be the positive definite subspace in which we split $\R^{a-1,0} \oplus \R^{0,a}$ for $\rho$. Define $U^i = \pi_{P}(F^i)$, for $\pi_P:\R^{a-1, a} \rightarrow P$ the orthogonal projection. 

\textbf{Case 1}: $U^i \cap f^i \neq \{0\} $. Then choose a vector $u \in Q_+(f^i \cap T^i)$. By definition, we then have a line $\ell = [u+z] \subset T^i$, for some $z \in Q_-(P^\bot)$. On the other hand, recall that $\mathrm{image}(\iota \circ \rho) \cong \mathbb{S}^{a-1}$ ``acts'' simply transitively on $Q_-(P^\bot)$. Now, for  any flag $F_{(0)} \in \Lambda$, for some $w \in Q_-(P^\bot)$, we have $[u+w] \subset F_{(0)}^i$. Hence, by the action of $\mathbb{S}^{a-1}$, there is some flag $F \in \Lambda$ such that $\ell \subset F^i$. We then conclude that $\ell \subset F^i \cap T^i$, so $F^i \not \pitchfork T^i$.  

\textbf{Case 2}: $U^i \cap f^i = \{0\}$. Write $T^i = \graph(\varphi_i)$ for a linear map $\varphi_i: U^i \rightarrow P^\bot$ and set $W^i \coloneqq \image(\varphi_i)$. Then take any vector $u \in Q_+(f^i \cap (U^i)^\bot)$ and $z \in Q_-(P^\bot \cap (W^i)^\bot)$. Then as before, we find some other flag $F_{(0)} \in \Lambda$ such that $[u+z] \in F_{(0)}^i$. However, we then see that $[u+z]$ is orthogonal to all of $T^i$. Hence, $F^i \cap (T^i)^\bot \neq \{0\}$, so $F^i \not \pitchfork T^i$. 
\end{proof}

Now, by combining Remark \ref{Remk:CompositionAlgebra} with Corollary \ref{Cor:CompAlgMaxTransverlity} and Theorem \ref{Thm:G2Fibration}, we conclude the following. 

\begin{corollary}\label{Cor:G2SphereMaxTransversality}
Let $\Lambda < \mathcal{F}_{1,2}^\times$ be a transverse 3-sphere in the $\Gtwosplit$-full flag manifold from Theorem \ref{Thm:G2Fibration}. Then $\mathrm{pr}_i(\Lambda) \subset \Iso_i(\R^{3,4})$ is maximally transverse for $1 \leq i \leq 3$. 
\end{corollary}

\subsection{Fibrations of \texorpdfstring{$B_3$}{B3}-Full Flag Manifold by Maximally Transverse 3-spheres}\label{Sec:B3Fibration}

As a result of Theorem \ref{Thm:DivisionAlgebras}, in the case of $G = \SO_0(3,4)$, there are some geometrically interesting fibrations of the full flag manifold $\Iso_{\dbrack{3}}(\R^{3,4})$, corresponding to $\mathbb{A} = \Ha$.

To start, fix $P \in \mathbb{X}_{\SO_0(3,4)} \cong \Gr_{(3,0)}(\R^{3,4})$. Consider the pointwise stabilizer $H_P \coloneqq \Stab_{G}^{pt}(P) \cong \SO(4)$. There is a fibration $\pi_P:\Iso_{\dbrack{3}}(\R^{3,4})\rightarrow \Flag(P)$ that is also a principal $H_P$-bundle, seen as follows. The map $\pi_P$ orthogonally projects the component subspaces of an isotropic full flag to $P$. We will denote $\Iso_{\dbrack{3}}(\R^{3,4})|_f$ as the fiber over $f \in \Flag(P)$. 
Now, fix a point $f \in \Flag(P)$ as well as an orthonormal basis $B=(u_1, u_2,u_3)$ for $P$ such that $f$ is the standard flag in the basis $B$. Then every flag $F \in \Iso_{\dbrack{3}}(\R^{3,4})|_{f}$ obtains the form $F = (F^1, F^2, F^3)$ with $F^i = \spann_{\R} \{ u_j +z_j \}_{j =1}^i $. Let us then write $F^i= F^i(\bm{z})$, for $\bm{z} = (z_1,z_2,z_3)$. One is then led to consider the Stiefel manifold 
\[ V_3(P^\bot) = \{ \,(z_1, z_2, z_3) \in Q_-(P^\bot)^3 \; |\; z_1 \cdot z_2 = z_1 \cdot z_3 = z_2 \cdot z_3 = 0\} . \] 
In summary, relative the choice of $B$, we obtain a diffeomorphism $V_{3}(P^\bot) \cong \Iso_{\dbrack{3}}(\R^{3,4})|_f$ given by $(z_1, z_2, z_3) \mapsto F(\bm{z})$, where $F(\bm{z}) =( F^1(\bm{z}), F^2(\bm{z}), F^3(\bm{z}))$ as above. It is evident that $H_P$ preserves the fiber $\Iso_{\dbrack{3}}(\R^{3,4})|_f$. Now, since $\SO(P^\bot) \cong H_P$ acts simply transitively on the Stiefel manifold $V_3(P^\bot)$, we conclude that $H_P$ acts simply transitively on the fiber $\Iso_{\dbrack{3}}(\R^{3,4})|_f$.

Following the above reasoning, one can prove more generally that there is an $\SO(n+1)$-principal bundle fibration $\Iso_{\dbrack{n}}(\R^{n,n+1}) \rightarrow \Flag(\R^n)$. That is, the $B_n$-full flag manifold fibers over the $A_{n-1}$-full flag manifold with $V_n(\R^{n+1})\cong_{\mathrm{Diff}} \SO(n+1)$ fibers. 

\begin{proposition}\label{Prop:BnAnFibration}
Fix $P \in \Gr_{(n,0)}(\R^{n,n+1})$. 
The fiber bundle $\Iso_{\dbrack{n}}(\R^{n,n+1}) \rightarrow \Flag(P)$ is a principal $H_P$-bundle for $H_P = \Stab_{\SO_0(n,n+1)}^{pt}(P) \cong_{\mathrm{Lie}} \SO(n+1)$.
\end{proposition}

We remark that the exact same argument proves an analogous result, namely that there is a principal bundle fibration $\SO(n) \rightarrow \Flag(\R^{n,n}) \rightarrow \Flag(\R^n )$. 
\begin{proposition}\label{Prop:DnAnFibration}
Fix $P \in \Gr_{(n,0)}(\R^{n,n})$.  
The fiber bundle $\Iso_{\dbrack{n}}(\R^{n,n}) \rightarrow \Flag(P)$ is a principal $H_P$-bundle for $H_P = \Stab_{\SO_0(n,n)}^{pt}(P) \cong_{\mathrm{Lie}} \SO(n)$.
\end{proposition}

While the fibration in Proposition \ref{Prop:BnAnFibration} is a general result, there is additional structure when $n=3$. 
In this case, the $\SO(4)$-fiber $\Iso_{\dbrack{3}}(\R^{3,4})|_f$ is itself a principal $\Sp(1)$-bundle over $\RP^3$ as follows. 
We now pick out a `model' flag from each $\Sp(1)$-orbit in  $\Iso_{\dbrack{3}}(\R^{3,4})|_f$. Select a basis $U= (u_i)_{i=1}^3$ of $P$ such that $f$ is the standard flag in the basis $U$.
We will show that relative to $U$, there is a natural identification $\Iso_{\dbrack{3}}(\R^{3,4})|_f \cong V_2(\mathrm{Im} \, \Ha) \times \Sp(1)$.
Here, $V_2(\mathrm{Im} \, \Ha)$ is the Stiefel manifold of ordered orthonormal pairs $\bm{z} = (z_2, z_3) \in Q_+(\mathrm{Im} \, \Ha)^2$. Now, define a flag $F(\bm{z}) \in \Iso_{\dbrack{3}}(\R^{3,4})$ for $\bm{z} \in V_2(\mathrm{Im} \, \Ha)$ as follows:  
$$ F(\bm{z}) = \left [ \, \langle u_1+1_{\Ha}\rangle  \subset (F^1 \oplus  \spann \{\,u_2 + z_2 \}\,) \subset \,\left( F^2 \oplus \spann \{\,u_3+z_3 \} \right) \, \right].$$
It is evident that the $\Sp(1)$-orbits of $F(\bm{z})$ and $F(\bm{z}')$ do not intersect for $\bm{z} \neq \bm{z} \in V_2(\mathrm{Im} \, \Ha)$. Moreover, every $\Sp(1)$-orbit in $\Iso_{\dbrack{3}}(\R^{3,4})$ contains a point of the form $F(\bm{z})$. Thus, we have identified $\Iso_{\dbrack{3}}(\R^{3,4})|_f \cong V_2(\mathrm{Im} \, \Ha) \times \Sp(1)$ as desired. One notes that $V_2(\mathrm{Im} \, \Ha) \cong \SO(3) \cong \RP^3$.  

\begin{figure}[ht]
\centering
\includegraphics[width=.65\textwidth]{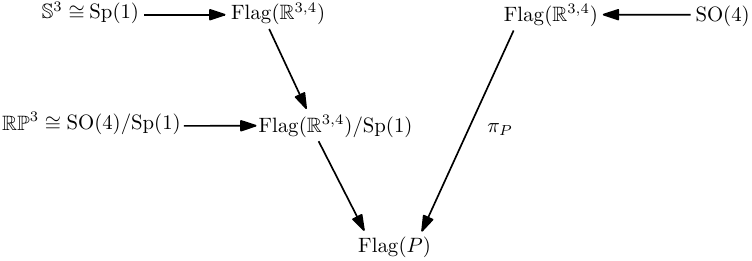}
\caption{\emph{\small{Fiber bundle structure of $\Iso_{\dbrack{3}}(\R^{3,4})$ relative to choice of $P \in \mathbb{X}_{\SO_0(3,4)}$, including two principal bundle fibrations. A similar diagram applies for $\Flag(\R^{4,4})$ after changing the base to $\Flag(\R^4)$.}}}
\label{Fig:B3FibrationSimple}
\end{figure}

Stitching together these two perspectives in the case $n =3$ on the total space $\Iso_{\dbrack{3}}(\R^{3,4})$ yields the diagram in Figure \ref{Fig:B3FibrationSimple}. In particular, each $\SO(4)$-orbit in $\Iso_{\dbrack{3}}(\R^{3,4})$ consists of an $\RP^3$-family of maximally transverse 3-spheres. 

\subsection{Fibrations of \texorpdfstring{$D_4$}{D4}-Full Flag Manifold by Maximally Transverse \texorpdfstring{$3$}{3}-spheres}\label{Subsec:D4Fibration}

In this section, we highlight the remarkable structure of $\Iso_{\dbrack{3}}(\R^{4,4})$ as a consequence of Theorem \ref{Thm:DivisionAlgebras}. In particular,  $\Iso_{\dbrack{3}}(\R^{4,4})$ is very nearly the product of four transverse 3-spheres. \medskip 

The maximal compact subgroup $K \cong \SO(4)\times \SO(4)$ of $G = \SO_0(4,4)$ acts transitively on the full flag manifold $\Iso_{\dbrack{3}}(\R^{4,4})$ with finite stabilizer. Crucial to our analysis of this special case is the isomorphism $\Spin(4) \cong \Sp(1)\times \Sp(1)$. This isomorphism is induced by the surjective 2-1 map $\Sp(1) \times \Sp(1) \rightarrow \SO(4) $ by $(x, y)\mapsto L_x \circ R_{y^{-1}}$, where we identify $\Sp(1) \cong Q_+(\Ha)$. 
As a consequence, $\Iso_{\dbrack{3}}(\R^{4,4}) \cong (\mathbb{S}^3)^4/Z$ for some finite group $Z$. We will show there are four distinct free actions of $\Sp(1) \cong \mathbb{S}^3$ on $\Iso_{\dbrack{3}}(\R^{4,4})$ whose orbits are copies of maximally transverse 3-spheres that pairwise intersect at most finitely many times.\medskip 

To select the desired copies of $\Sp(1)$, fix a linear identification $\R^{4,4} \cong \Ha^{1,0} \oplus \Ha^{0,1}$, where $\Ha^{k,l}$ denotes $\Ha^{k+l}$ equipped with a signature $(k,l)$ quaternion-hermitian form. The subgroup $K <\SO_0(4,4)$ stabilizing this splitting is a copy of the maximal compact of $G$. The relevant $\Sp(1)$-subgroups shall be denoted $H_L^{\pm}, H_{R}^{\pm}$. Recall the standard isomorphism $Q_+(\Ha) \cong \Sp(1)$ by $x \mapsto L_x$. With this perspective, we can build the map $\varphi: Q_+(\Ha)^4 \rightarrow K$ by 
\[ \varphi(u,v,w,x) = L_u\circ R_{v^{-1}} \oplus  L_w\circ R_{x^{-1}}.\]
The map $\varphi$ is a 4--1 cover of $\Sp(1)^4$ onto $K$. The subgroup $H_{L}^+$ corresponds to $Q_+(\Ha) \times \{\mathrm{id}\} \times \{\mathrm{id}\}\times \{\mathrm{id}\}$, $H_{L}^-$ corresponds to 
$\{\mathrm{id}\} \times \{\mathrm{id}\} \times Q_+(\Ha) \times \{\mathrm{id}\}$, and similarly for $H_R^{\pm}$. 

\begin{corollary}\label{Thm:D4Fibration}
Via the map $\varphi$, the group $\Sp(1)^4$ acts transitively on $\Iso_{\dbrack{3}}(\R^{4,4})$, with finite stabilizer. Moreover, each $\Sp(1)$-subgroup $H_{L/R}^{\pm}$ has maximally transverse orbits. 
\end{corollary}

\begin{proof}
The only remaining claim to be proven is that each $\Sp(1)$-subgroup has maximally transverse orbits. However, this follows immediately from the proof of Theorem \ref{Thm:DivisionAlgebras}. Indeed, in the case of $\Iso_{\dbrack{3}}(\R^{4,4})$, there is a symmetry between the spacelike and timelike 4-planes, and acting by $Q_+(\Ha)$ on the left/right on either $\R^{4,0}$ or $\R^{0,4}$ will produce a transverse 3-sphere in $\Iso_{\dbrack{3}}(\R^{4,4})$ by the same argument. 
\end{proof}

\subsection{Obstructing Transverse Spheres in Low Dimensions}

In this subsection, we reformulate arguments of Canary-Tsouvalas from \cite{CT20} to obstruct transverse $k$-spheres, for $k \in \{2,4,8\}$ in certain partial flag manifolds, including the full flag manifolds of $\Flag(\R^{2k})$. 

\begin{proposition}[Obstructing Transverse Spheres]\label{Prop:FibrationsSpheres}
    Let $k \geq 2$ be a positive integer and $\xi\colon \mathbb{S}^k \rightarrow \Gr_{k}(\R^{2k})$ a continuous transverse map. 
    Then $\xi$ does not lift continuously to $\Flag_{\{1, k\}}(\R^{2k})$. 
    In particular, if $\{1,k\} \subset \Theta$, then there is no transverse $k$-sphere in $\Flag_{\Theta}(\R^{2k})$. 
\end{proposition}

\begin{proof}
The reasoning is identical to part of the proof of \cite[Theorem 1.3]{CT20}, though their result is not stated as such. 
We recall the argument for completeness. 

Given such a transverse map $\xi$, form the fiber bundle $\mathbb{RP}^{k-1} \rightarrow E_{\xi} \rightarrow \mathbb{S}^k$ by $E_{\xi} \coloneqq \xi^*\mathbb{P}(E^k)$, where $\tau\colon E^k\rightarrow \Gr_{k}(\R^{2k})$ is the tautological vector bundle. 
Transversality of $\xi$ implies that if $\ell_i \in (E_{\xi})_{x_i}$ for $x_1 \neq x_2 \in \mathbb{S}^k$, then $\ell_1 \neq \ell_2$. 
Hence, we may conflate the total space of $E_{\xi}$ with its image in $\mathbb{RP}^{2k-1}$ under the tautological map $(p,\ell) \mapsto \ell$, which is injective by the previous sentence. 
Observe that $E_{\xi}$ is evidently closed in $\mathbb{RP}^{2k-1}$. 
On the other hand, $E_{\xi}$ is open by dimension count. 
Thus, the total space of $E_{\xi}$ is homeomorphic to $\mathbb{RP}^{2k-1}$.
It follows that $E_{\xi}$ exhibits a fibration $\mathbb{RP}^{k-1} \rightarrow \mathbb{RP}^{2k-1} \rightarrow \mathbb{S}^k$. 

We claim that $\pi: E_{\xi} \rightarrow \mathbb{S}^k$ admits no section, which implies $\xi$ does not lift to $\Flag_{\{1,k\}}(\R^{2k})$. This follows immediately by the fact that if $s: \mathbb{S}^k \rightarrow \mathbb{RP}^{2k-1}$ were a section, then $(\id)_{*} = \pi_*\circ s_*$. This equality is impossible since $H_k(\mathbb{RP}^{2k-1}) \in \{0, \Z_2\}$, so $\pi_*$ is trivial, yet on the other hand  $\id_{*}$ is the identity on $H_{k}(\mathbb{S}^k) \cong \Z$.  
\end{proof}

\begin{corollary}
Let $p \in \{2,4,8\}$. If $\{1,p\} \subset \Theta$, there is no transverse $p$-sphere in $\Iso_{\Theta}(\R^{p,p})$.  
\end{corollary}

\begin{remark}
The result in Proposition \ref{Prop:FibrationsSpheres} is stated in unnecessary generality. 
As remarked in \cite{CT20} and is well-known, the only possible fibrations $\mathbb{S}^p \rightarrow \mathbb{S}^q \rightarrow \mathbb{S}^r$ of spheres by spheres occur for triples $(1,3,2), (3,7,4), (7,15,8)$. 
The existence of such fibrations is demonstrated by the Hopf fibrations of the normed division algebras $(\mathbb{A}_n)_{n=1}^4 =( \R,  \C, \Ha , \Oct)$. 
Set $a\coloneqq \dim_{\R}(\mathbb{A})$. 
Using $\mathbb{P}^1(\mathbb{A}) \cong \mathbb{S}^{a}$, one finds fibrations $\mathbb{S}^{a-1} \rightarrow \mathbb{S}^{2a-1} \rightarrow \mathbb{S}^a$ by 
\[ Q_+(\mathbb{A}) \rightarrow Q_+(\mathbb{A}^2) \rightarrow \mathbb{P}^1(\mathbb{A}), \;\]
where the projection $Q_+(\mathbb{A}^2) \rightarrow \mathbb{P}^1(\mathbb{A}) \cong \mathbb{A} \cup \{\infty\} $ is given by $(u,v)\mapsto u/v$. 
Projectivizing, one obtains an associated fibration 
\[ \mathbb{RP}^{a-1} \rightarrow \mathbb{RP}^{2a-1} \rightarrow \mathbb{S}^a,\]
with projection $\mathbb{RP}^{2a-1} \rightarrow \mathbb{S}^a$ now by $[u:v]\mapsto u/v$. 
Thus, the essence of Proposition \ref{Prop:FibrationsSpheres} is contained in the cases $k \in \{1,2,4,8\}$.

Additionally, as observed in \cite{CT20}, the projective lines $\mathbb{P}^1(\mathbb{A})$ yield noteworthy transverse spheres. In particular, for $\mathbb{A} \in \{\R,\C, \Ha, \Oct\}$, we have a transverse map $\mathbb{P}^1(\mathbb{A}) \stackrel{\pitchfork}{\hookrightarrow} \Gr_{a}(\R^{2a})$, for $a = \dim_{\R}(\mathbb{A})$. 
Proposition \ref{Prop:FibrationsSpheres} implies that this transverse sphere cannot be lifted to a transverse sphere in $\Flag(\R^{2a})$. 
\end{remark}

\section{Maximal Transversality via the Atiyah-Bott-Shapiro Isomorphism}\label{Sec:MaximalTransversality}

The first goal of this section is to prove the following theorem.

\begin{theorem}[Maximal Transversality of Spinor Spheres]\label{Thm:ABSMaximallyTransverse}
    Let $n \geq 2$ be an integer. 
    The transverse spheres $\Lambda \cong \mathbb{S}^{n-1}$ in $\Flag(\R^{d-1, d})$ and in $\Flag(\R^{d,d})$ constructed in Theorem \ref{Thm:TransverseSpinorSphere} with $p \in \{d-1, d\}$ are maximally transverse if and only if $n \in \{0,1,2,4\} \bmod 8$. 
\end{theorem}

As a consequence we will be able to show that many transverse spheres in full flag manifolds of other split Lie groups are maximally transverse.
The proof rests on the powerful machinery of Atiyah-Bott-Shapiro \cite{ABS64}, which relates the reduced $K$-groups $\widetilde{\mathrm{KO}}(\mathbb{S}^{n+1})$, measuring a certain kind of non-triviality of real vector bundles over $\mathbb{S}^{n+1}$, with certain Grothendieck groups $\mathfrak{M}_n$ depending only on representation theory of the Clifford algebra $\Cl(n)$. 
In particular, leveraging the Atiyah-Bott-Shapiro isomorphism allows us to reduce a hard problem in algebraic topology to an easy calculation in representation theory. 

\subsection{Preliminaries on Topological \texorpdfstring{$K$}{K}-theory}

In this subsection, we introduce the essential objects from topological $K$-theory necessary for the proof of Theorem \ref{Thm:ABSMaximallyTransverse}. 
Additional background on Grothendieck group completion can be found in Appendix \ref{Appendix:GGC}.
\medskip

We first briefly recall the groups $\mathrm{KO}(X)$ and $\widetilde{\mathrm{KO}}(X)$ for $X$ a compact Hausdorff space. 

For further exposition on topological $\mathrm{K}$-theory, we refer the reader to \cite[Chapter 1, Section 9]{LM89}. 
Here, we follow the treatment of \cite{Hat17}.

The elements of $\KO(X)$ and $\widetilde{\KO}(X)$ represent certain isomorphism classes of real vector bundles on $X$. 
All vector bundles going forwards will be real unless otherwise specified. 
Now, let $\epsilon^i$ denote the trivial rank $i$ vector bundle over $X$. 
Consider the equivalence relation $E_1 \sim_s E_2$ when $E_1 \oplus \epsilon^i \cong E_2 \oplus \epsilon^i$. 
We say $E_1 $ and $E_2$ are \textbf{stably isomorphic} when $E_1 \sim_s E_2$. 
Then the group $\KO(X)$ consists of stable isomorphism classes of vector bundles, with group operation $[E_1] + [E_2] = [E_1\oplus E_2]$. 
It is evident $\KO(X)$ is a monoid as defined. 
The fact that any vector bundle $E \rightarrow X$ admits a complementary bundle $E^\bot \rightarrow X$ such that $E \oplus E^\bot = \epsilon^i$ for some $i$ (\cite[Ch 1, Corollary 9.89]{LM89}) implies $\KO(X)$ is actually a group. 

The group $\widetilde{\KO}(X)$ is a reduction of $\KO(X)$, and the two are related by $\KO(X) = \widetilde{\KO}(X) \oplus \Z$. 
We can realize $\widetilde{\KO}(X)$ as follows. 
Let $E_i$ be vector bundles on $X$ for $i \in \{1,2\}$. Define $E_1 \approx_s E_2$ when there exist non-negative integers $n, m \geq 0$ such that $E_1 \oplus \epsilon^n \cong E_2 \oplus \epsilon^m$. 
The group $\widetilde{\KO}(X)$ then consists of $\approx_s$ equivalence classes of vector bundles on $X$, also with $\oplus$ as the group operation. 
By definition, there is a natural surjective group homomorphism $\KO(X) \rightarrow \widetilde{\KO}(X)$.  

We denote by $\Vect(X)$ the set of isomorphism classes of real vector bundles on $X$, which is naturally a monoid under direct sum. 
Note that by definition there is a natural projection map $\Vect(X) \rightarrow \KO(X)$.  
The Grothendieck group construction of $\KO(X) \coloneqq \mathrm{K}(\Vect(X))$ is precisely to upgrade the monoid $\Vect(X)$ to an honest group. 
We include further basic details on Grothendieck groups, specifically related to $\KO(X),\widetilde{\KO}(X)$ as well as $\mathfrak{M}_n, \hat{\mathfrak{M}}_n$ (to be defined shortly), in Appendix \ref{Appendix:GGC}. \medskip 

\subsection{Grothendieck Group of Clifford Modules }

In this subsection, we recall from \cite{ABS64} the next important ingredient in the Atiyah-Bott-Shapiro isomorphism -- the Grothendieck groups of Clifford modules. 
These groups come in two flavors: ungraded and $\Z_2$-graded. 
We also prove a technical lemma needed for the proof of Theorem \ref{Thm:ABSMaximallyTransverse}. 

Let $\hat{M}_k$ denote the set of isomorphism classes of $\Z_2$-graded real $\Cl(k)$-modules, which is a monoid under direct sum. We shall study the group $\hat{\mathfrak{M}}_k \coloneq \mathrm{K}(\hat{M}_k)$, the Grothendieck group of $\hat{M}_k$. Similarly, let $M_k$ denote the monoid of isomorphism classes of ungraded real $\Cl(k)$-modules and $\mathfrak{M}_k \coloneqq \mathrm{K}(M_k)$ the Grothendieck group of $M_k$. The groups $\mathfrak{M}_k, \hat{\mathfrak{M}}_k$ are each free abelian groups  generated by a maximal set of mutually non-isomorphic irreducible $\Cl(k)$-modules that are ungraded and graded, respectively \cite{LM89}. Thus, each of $\mathfrak{M}_k, \hat{\mathfrak{M}}_k$ is isomorphic to either $\Z$ or $\Z \oplus \Z$, depending only on $k \bmod8$. See Appendix \ref{Appendix:GGC} for further details. 

Before the proof of maximal transversality, we shall need the following technical lemma, which is a slightly stronger version of Proposition \ref{prop:SpinIrreducible?}. 
That is, when $n \in \{0,1,2,4\} \bmod 8$, the unique irreducible real $\Cl(n)$-representation is actually $\Z_2$-graded. 
This result seems to only be partially proven in \cite{LM89}, so we provide a complete argument here.\footnote{To be more precise, \cite[Ch 1, Proposition 5.10]{LM89} and \cite[Ch 1, Proposition 5.12]{LM89} describe how an irreducible $\Cl(n)$-module decomposes into $\Cl^0(n)$-submodules, but does not produce a $\Z_2$-grading of the whole $\Cl(n)$-module. Now, when $n \in \{0,3\} \bmod 4$, \cite[Ch 1, Proposition 3.6]{LM89} produces the desired $\Z_2$-grading on any $\Cl(n)$-module, but there is no analogue of this result presented when $n \in \{1,2\} \bmod 4$.} 

\begin{lemma}[$\Z_2$-grading of $\Cl(n)$-irreps]\label{Lem:Z2GradedCliffordModules}
    Let $n \in \{0,1,2,4\} \bmod 8$. 
    Then any irreducible $\Cl(n)$ representation $\eta\colon \Cl(n) \rightarrow \End_{\R}(W)$ admits a $\Z_2$-grading $W^{\bullet}$ such that $\Cl^{i} \cdot W^k \subseteq W^{i+k}$. 
\end{lemma}

\begin{proof}
Suppose $n \in \{0,1,2,4\} \bmod 8$. 
Then there is a unique irreducible (real) representation of $\Cl(n)$ up to isomorphism \cite[Ch 1, Theorem 5.8]{LM89}. 
Thus, to prove every irreducible representation of $\Cl(n)$ is $\Z_2$-graded, it suffices to produce a single $\Z_2$-graded irreducible representation. 
Here, we produce such a grading first in the base cases $n \in \{1,2,4,8\}$, then use periodicity of $\Cl(n)$ to handle the general case.

Now, let us make an observation. 
If the representation $\Cl(V) \rightarrow \End_{\R}(W)$ is induced by a defining map $c\colon V \rightarrow \End_{\R}(W)$, then the $\Cl(V)$-module $W$ is $\Z_2$-graded if and only if there is a splitting $W = W^0\oplus W^1$ such that 
$c(x) = \begin{pmatrix} 0 & * \\ * & 0\end{pmatrix}$ for all $x \in V$ in this splitting. 
Armed with this observation, let us now handle the base cases with \emph{normed division algebras}. 
Some of the following representations are actually complex or quaternionic (i.e.\ the \emph{intertwiner algebra} is isomorphic to $\C$ or $\Ha$), but such structure is not needed presently. \medskip 

\textbf{Cl(1) model.} 
Set $V_1 \coloneqq \mathrm{Im} \C$ and $S_1 = \C$. 
The map $c_1\colon V_1 \rightarrow \End_{\C}(\C)$ by $c_1(x) = L_{x}$ satisfies $c_1(x)^2 = -q(x)\id$ and induces an isomorphism $\Cl(1) \cong \End_{\C}(\C) \cong \C$ of $\R$-algebras, as well as an irreducible (real) representation $c_1\colon \Cl(1) \rightarrow \End_{\R}(\C)$. 
Define the $\Z_2$-grading by $S_1^0 = \R$ and $S_1^1 = \mathrm{Im}\C$.
One quickly verifies $S_1^{\bullet}$ is a $\Z_2$-graded $\Cl(1)$-module.\medskip 

\textbf{Cl(2) model.}  
Set $V_2 = \spann_{\R} \{ i, j\} \subset \Ha$ and $S_2 = \Ha$. 
Next, consider the defining map $c_2\colon V_2 \rightarrow \End_{\R}(\Ha)$ given by $c_2(x) = L_x$. 
Since $\Cl(2) \cong \Ha$ as $\R$-algebras, the non-trivial map $c_2\colon \Cl(V_2) \rightarrow \End_{\Ha}(\Ha) \cong \Ha$ must be an isomorphism. 
Moreover, $c_2\colon \Cl(2) \rightarrow \End_{\R}(\Ha)$ is also irreducible. 
Now, set $S_2^{0} = \spann_{\R} \{1, k\}$ and $S_2^{1} = \spann_{\R} \{i ,j\}$. 
Since the above splitting satisfies $S_2^{i}S_2^j = S_2^{i+j}$, we do achieve a $\Z_2$-graded $\Cl(2)$-module $S_2^{\bullet}$. \medskip

The previous defining reps are well-known. 
One can use an idea originally by \cite{Bry20} for $n =4$ and $n= 8$, extended by \cite[(13.1)]{Esc18} and \cite[page 13]{San25} from different perspectives. \medskip 

\textbf{Cl(4) model}. 
Set $S_4 = \Ha $. 
Consider the defining map $c_4 \colon \Ha \rightarrow \End_{\R}(\Ha^2)$ by
\begin{align}\label{c4Def}
     c_4(x) = \begin{pmatrix} 0 & -L_{x^*} \\ L_x & 0 \end{pmatrix}.
\end{align}
Since $q_{\Ha}(x) =xx^*$, then associativity of $\Ha$ implies $c_4$ satisfies the Clifford identity. 
Since $\Cl(4) \cong \Mat_2(\Ha)$ as $\R$-algebras, $c_4$ induces an isomorphism $\Cl(4) \cong \End_{\Ha}(\Ha^2)$ by dimension count and simplicity of $\Cl(4)$. 
The splitting $\Ha \oplus \Ha$ is the desired $\Z_2$-grading. \medskip 

\textbf{Cl(8) model.} 
We proceed as with $\Cl(4)$, replacing $\Ha$ by $\Oct$. Consider the defining map $c_8: \Oct \rightarrow \End_{\R}(\Oct^2)$ via \eqref{c4Def}. Again, $c_8$ induces an isomorphism $\Cl(8) \cong \End_{\R}(\Oct^2)$. Here, $S_8 = \Oct \oplus \Oct$ is the $\Z_2$-grading preserved by $\Cl(8)$. \medskip 

\textbf{The general case}. 
Let $k \geq 1$ be an integer. 
Recall the isomorphism $\Cl(8k+r) \cong \Cl(8k) \otimes_{\R} \Cl(r)$ from Theorem \ref{thm:CliffordPeriodicity}. 
In particular, $\Cl(8k) \cong \Mat_{16^k}(\R)$. 
Choose any irreducible representation $c_{8k} \colon \Cl(8k) \rightarrow \End_{\R}(S_{8k})$, which is an $\R$-algebra isomorphism, for $S_{8k}$ a real $16^k$-dimensional vector space. 
Let $r \in \{1,2,4,8\}$, and $S_r$ be the explicit $\Z_2$-graded real $\Cl(r)$-module built above. 
We can define an irreducible representation 
$c_{8k+r}\colon \Cl(8k+r) \rightarrow \End_{\R}(S_{8k} \otimes S_{r})$ by identifying $\Cl(8k+r) \cong \Cl(8k)\otimes_{\R}\Cl(r)$ and then taking the linear map induced by 
$c_{8k+r}(x \otimes_{\R} y) = c_{8k}(x) \otimes_{\R} c_r(y)$. 
As noted in \cite[page 39]{LM89}, this representation is also irreducible. 
Define $S_{8k+r} \coloneqq S_{8k} \otimes_{\R} S_{r}$. 
Now, since $\Cl(8)$ admits a $\Z_2$-graded module $S_8$, we can suppose that the $\Cl(8k)$-module $S_{8k}$ is also $\Z_2$-graded. 
Next, one notes that under the isomorphism $\Cl(8k+r) \cong \Cl(8k)\otimes_{\R} \Cl(r)$, we have the induced $\Z_2$-grading:
\[ \begin{cases}
    \Cl^0(8k+r) = & \Cl^0(8k)\otimes_{\R} \Cl^0(r)+\Cl^1(8k) \otimes_{\R}\Cl^1(r) \\
    \Cl^1(8k+r) = & \Cl^1(8k)\otimes_{\R} \Cl^0(r)+\Cl^0(8k) \otimes_{\R}\Cl^1(r). 
\end{cases} \] 
Finally the $\Z_2$-grading $S_{8k+r}^{\bullet}$ is given by:
\[ 
\begin{cases}
    S_{8k+r}^{0} \coloneqq S_{8k}^{0} \otimes_{\R} S_{r}^0+S_{8k}^{1} \otimes_{\R} S_{r}^1\\
    S_{8k+r}^{1} \coloneqq S_{8k}^{1} \otimes_{\R}S_{r}^0 +S_{8k}^{0} \otimes_{\R} S_{r}^1.
\end{cases}\]
By checking each block $S_{8k}^{i} \otimes_{\R}S_{r}^j$, one finds the $\Z_2$-gradings $\Cl^{\bullet}(8k+r)$ and $S_{8k+r}^{\bullet}$ are compatible and hence $S_{8k+r}^{\bullet}$ is a $\Z_2$-graded $\Cl(8k+r)$-module. 
\end{proof}

\subsection{Proof of Maximal Transversality}

We now prove the maximal transversality of the spheres in the $B_{d-1}$ and $D_{d}$ cases. We can bootstrap off these cases after the theorem. 

\begin{proof}[Proof of Theorem \ref{Thm:ABSMaximallyTransverse}]
    The proof is nearly the same in the $B_{d-1}$ and $D_d$ cases. We argue in the $D_d$ case. 
     \medskip 

    Note that if $n \in \{3,5,6,7\} \bmod 8$, then the $(n-1)$-sphere $\Lambda$ is contained in a transverse $n$-sphere in $\Flag(\R^{d,d})$ by Corollary \ref{Cor:IsotropicTransverseSphere} and definition of $d =d(n)$ in Theorem \ref{SpinModuleDimensions}. 
    Thus, $\Lambda$ is not maximally transverse in these cases.\medskip 

    Now, suppose $n \in \{0,1,2,4\} \bmod 8$. We will show that the inclusion $\xi_n \colon \mathbb{S}^{n-1} \to \Flag(\R^{d,d})$ is not nullhomotopic and then invoke Fact \ref{Fact:MaximallyTransverse}. 

    Recall that the maximal compact subgroup $\SO(d) \times \SO(d) $ of $\SO_0(d,d)$ is a finite covering of $\Flag(\R^{d,d})$. 
    For a chosen basepoint $F \in \Flag(\R^{d,d})$, the covering map $orb \colon \SO(d)\times \SO(d) \rightarrow \Flag(\R^{d,d})$ takes the form $g \mapsto g \cdot F$. 
    Under this perspective, we have $\xi_n = orb \circ \hat{\xi}_n$, for $\hat{\xi}_n \colon \mathbb{S}^{n-1} \rightarrow \SO(d) \times \SO(d)$.
    In particular, since $n\ge 2$, the map $\hat{\xi}_n$ is not nullhomotopic if and only if the same is true for $\xi_n$. 
    By the construction of the transverse sphere, the lift $\hat{\xi}_n$ takes the form $\hat{\xi}_n=(\mathrm{id} ,f)$, where $f\colon \mathbb{S}^{n-1} \rightarrow \SO(d)$ is the map $f(x) = \eta(x_0x)$ and $\eta\colon \Spin(n) \rightarrow \SO(d)$ is a spinor representation induced by restricting an irreducible $\Cl(n)$-module $S$ to an irreducible $\Cl^0(n)$-submodule $S_n^+$. 
    Let $S = S_n^+ \oplus S_n^-$ be the decomposition of $S$ into irreducible $\Cl^0(n)$-submodules.
    Now, $\hat{\xi}_n$ is nullhomotopic if and only if $f$ is nullhomotopic. Furthermore, $f$ is nullhomotopic if and only if $\overline{f}\colon \mathbb{S}^{n-1} \rightarrow \mathrm{O}(S_n^+, S_n^-) \cong \mathrm{O}(d)$ by $\overline{f}(x) = \eta(x)|_{S_n^+}$ is nullhomotopic. 
    We will proceed to show $[\overline{f}] \in \pi_{n-1}(\mathrm{O}(d))$ is nontrivial. 

    Let $\Vect^k(X)$ denote the monoid of isomorphism classes of rank $k$ vector bundle on $X$. 
    Recall that by clutching functions, there is a natural bijection $\Vect^k(\mathbb{S}^n) \leftrightarrow  \pi_{n-1}(\mathrm{O}(k))$ \cite[Proposition 1.17]{Hat17}. 
    Thus, associated to $\overline{f}$ is an isomorphism class $[E_n] \in \Vect(\mathbb{S}^n)$. 
    Then $[E_n]$ is naturally associated to an element $\alpha_n \in \KO(X)$, which induces an element $\widetilde{\alpha_n}$ of $\widetilde{\KO}(\mathbb{S}^n)$.
    We will show that $\widetilde{\alpha_n} \neq 0$, which entails $[\overline{f}] \neq \id \in \pi_{n-1}(\mathrm{O}(d))$. 

    The remainder of the proof rests on the Atiyah-Bott-Shapiro isomorphism \cite[Theorem 11.5]{ABS64}. 
    This is an explicit ``geometric'' group homomorphism $\hat{\alpha}\colon \hat{\mathfrak{M}}_k \rightarrow \widetilde{\KO}(\mathbb{S}^k)$ that descends to an isomorphism $\hat{\alpha} \colon \hat{\mathfrak{M}}_k/\iota^*\hat{\mathfrak{M}}_{k+1} \rightarrow \widetilde{\KO}(\mathbb{S}^k)$. 
    The essential property is this: for $\eta\colon \Cl(V_k) \rightarrow \End(W^{\bullet})$ a $\Z_2$-graded $\Cl(V_k)$-representation describing a point in $\hat{\mathfrak{M}}_k/\iota^*\hat{\mathfrak{M}}_{k+1}$, the vector bundle representing $\hat{\alpha}(W^{\bullet})$ is precisely represented by the clutching function $[\beta] \in \pi_{k-1}(\mathrm{O}(d))$, where $\beta\colon Q_+(V_k) \rightarrow O(W)$ is given by $\beta(x) = \eta(x)$.
    
    Let us now define the map $\hat{\alpha}$ a bit more precisely, by unraveling \cite[Proposition 9.25]{LM89} and its relation to the definition of $\hat{\alpha}$ in \cite[Section 11]{ABS64}. 
    Let $S^{\bullet}$ be a $\Z_2$-graded $\Cl(k)$-module. 
    Then we form a vector bundle on $E(S)\rightarrow \mathbb{S}^k$ as follows: 
    let $\mathbb{D}_k$ be the open unit disk in $\R^k$ and $\overline{\R}^{k} = \R^{k} \cup \{\infty\} \cong \mathbb{S}^k$ the one-point compactification of $\R^{k}$. 
    Associated to $S^{\bullet}$, define the vector bundle $E(S^{\bullet}) \rightarrow \mathbb{S}^k$ by 
    \[ E(S^{\bullet}) \coloneqq \overline{\mathbb{D}}_k \times S^+ \cup_{\sigma} (\overline{\R}^k\backslash \mathbb{D}_k)\times S^-,\]
    where the gluing map $\sigma\colon \mathbb{S}^{k-1} \times S^+ \rightarrow \mathbb{S}^{k-1} \times S^-$ is given by 
    $\sigma(x, v) = x\cdot v$, where $x \cdot v$ denotes the action of $x \in \R^k \subset \Pin(k) <\Cl(k)$ on $v \in S$. 
    The above process defines a monoid homomorphism $\hat{\hat{\alpha}} \colon \hat{M}_k \rightarrow \Vect(\mathbb{S}^k)$ by $[S^{\bullet}] \rightarrow [E(S^{\bullet})]$. 
    By the universal property of $\hat{\mathfrak{M}}_k$, there is a unique induced group homomorphism $\hat{\alpha}:\hat{\mathfrak{M}}_k \rightarrow \Vect(\mathbb{S}^k)$ such that $\hat{\hat{\alpha}} = \hat{\alpha}\circ i$ for $i \colon \hat{M}_k \rightarrow \mathfrak{M}_k$. 
    A simple argument shows  $\iota^*\hat{M}_{k+1} \subset \ker(\hat{\hat{\alpha}})$, i.e., pullback of larger $\Z_2$-graded modules give trivial vector bundles (because their clutching data is homotopically trivial). Hence, the map $\hat{\alpha}$ descends to a homomorphism  $\hat{\mathfrak{M}}_k/\iota^*\hat{\mathfrak{M}}_{k+1} \rightarrow \Vect(\mathbb{S}^k)$. Then by \cite[Theorem 11.5]{ABS64}, post-composing by the projection to $\widetilde{\KO}(\mathbb{S}^k)$, we obtain finally an isomorphism $\hat{\alpha}: \hat{\mathfrak{M}}_k/\iota^*\hat{\mathfrak{M}}_{k+1} \rightarrow \widetilde{\KO}(\mathbb{S}^k)$.\footnote{In fact, Atiyah-Bott-Shapiro present the isomorphism slightly differently from us. 
    In the notation of \cite[page 65]{LM89}, associated to $S^{\bullet}$ is a  tuple $[V_0, V_1,\sigma]$, where $V_i =\underline{S^i} \rightarrow \overline{\mathbb{D}}_k$ are trivial bundles and  $\sigma: \underline{S}^0|_{\mathbb{S}^{k-1}} \rightarrow \underline{S}^1|_{\mathbb{S}^{k-1}}$ is the isomorphism written above, and to this tuple they associate an element of $\KO(\mathbb{S}^k)$ that is equivalent to what we wrote here. }
    
    Now, since $n \in \{0,1,2,4\} \bmod 8$, any irreducible $\Cl(n)$-module admits a $\Z_2$-grading by Lemma \ref{Lem:Z2GradedCliffordModules}.   
    The connection to the present situation is now clear: \cite[Theorem 11.5]{ABS64} shows that $\overline{f} \colon \mathbb{S}^{n-1} \to \mathrm{O}(d)$ is non-nullhomotopic as soon as we can prove that an irreducible $\Z_2$-graded $\Cl(n)$-module $W^{\bullet}$ (with $W^0 \cong \R^d$) represents a non-trivial element of $\hat{\mathfrak{M}}_n/\iota^*\hat{\mathfrak{M}}_{n+1}$.
    
    We can now pass back from $\Z_2$-graded modules to ungraded modules.
    Indeed, there is an equivalence between $\Z_2$-graded $\Cl(k+1)$-modules and ungraded $\Cl(k)$-modules, using that $\Cl^0(k+1) \cong \Cl(k)$ \cite[Proposition 5.20]{LM89}. 
    Recall that $M_k$ denotes the monoid of isomorphism classes of $\Cl(k)$-modules under direct sum and $\mathfrak{M}_k$ be the Grothendieck group of $M_k$. 
    We obtain an isomorphism $\hat{\mathfrak{M}}_{k}/\iota^*\hat{\mathfrak{M}}_{k+1} \rightarrow \mathfrak{M}_{k-1}/\iota^*\mathfrak{M}_{k}$, which is induced 
    by the semigroup isomorphism $\hat{M}_{k} \rightarrow M_{k-1}$ that restricts $\Z_2$-graded $\Cl(k)$-modules to ungraded $\Cl^0(k) \cong \Cl(k-1)$-modules. 

    To summarize the discussion above: via the Atiyah-Bott-Shapiro isomorphism, showing that $\overline{f} \colon \mathbb{S}^{n-1} \to \mathrm{O}(d)$ is non-nullhomotopic reduces to proving that if $W_{n-1}$ is an irreducible $\Cl(n-1)$-module for $n \in \{0,1,2,4\} \bmod 8$, that $[W_{n-1}]$ is a non-zero element of $\mathfrak{M}_{n-1}/\iota^*\mathfrak{M}_{n} \cong \iota^*\hat{\mathfrak{M}}_{n}/\iota^*\hat{\mathfrak{M}}_{n+1}\cong \widetilde{\KO}(\mathbb{S}^n)$. 
    Thus, we need only pursue some simple representation theory to finish the proof. \medskip 
    
    \textbf{Case 1: $\bm{n \in \{0,4\}}$ mod 8}. In this case, $\mathfrak{M}_{n-1} \cong \Z \oplus \Z = \langle W_+\rangle \oplus \langle W_- \rangle$, seen by Table \ref{Table:Clifford}, with $W_{\pm}$ chosen irreducible representations that generate $\mathfrak{M}_{n-1}$. Let $W_n$ be an irreducible representation of $\Cl(n)$ generating $\mathfrak{M}_n$. Since $n \equiv 0 \bmod 4$, the module $W_n$ is $\Z_2$-graded as a $\Cl(n)$-module and $\Cl^0|_{W_n^+}$ and $\Cl^0|_{W_n^-}$ are distinct irreducible representations of $\Cl^0(n)$ by \cite[Proposition 5.10]{LM89}. Hence, $\iota^*([W_n]) = [W_-] + [W_+]$. In particular, this means that $\mathfrak{M}_{n-1}/\iota^*\mathfrak{M}_n \cong \Z(W_+) \oplus \Z(W_-)/(1,1) \cong \Z$, and, more importantly, $[W_+], [W_-] \neq 0$ in the quotient $\mathfrak{M}_{n-1}/\mathfrak{M}_n$. This handles the cases $n \equiv 0 \bmod 4$. \medskip 

    \textbf{Case 2: $\bm{n \in \{1,2\}}$ mod 8}. In this case, $\mathfrak{M}_{n-1} \cong \Z \cong \mathfrak{M}_n$. Let $W_{k} \in \mathfrak{M}_{k}$ be generators associated to irreducible $\Cl(k)$-representations for $k \in \{n-1, n\}$. Then $\iota^*W_n = 2 W_{n-1}$, so $0 \neq [W_{n-1}]\in \mathfrak{M}_{n-1}/\iota^*\mathfrak{M}_n\cong \Z_2$.
    \end{proof}

\begin{table}[ht] 
\renewcommand{\arraystretch}{1.3}
\centering
\resizebox{0.9\textwidth}{!}{
\begin{tabular}{ |c|c|c|c|c| } 
\hline
Dimension $n \bmod 8$ & Grothendieck group $\mathfrak{M}_n$ & $\mathfrak{M}_{n-1}/\iota^*\mathfrak{M}_{n}$ & $\widetilde{\KO}(\mathbb{S}^n)$ & Stable Homotopy Group $\pi_{n-1}(\mathbf{O})$\\
\hline
0 & $\Z$ & $\Z$ & $\Z$ & $\Z$ \\
\hline 
1 & $\Z$ & $\Z_2$ & $\Z_2$ & $\Z_2$ \\ 
\hline 
2 & $\Z$ & $\Z_2$ & $\Z_2$ & $\Z_2$\\ 
\hline 
3 &  $\Z\oplus \Z$ & 0 & 0 & 0  \\ 
\hline
4 & $\Z$ & $\Z$ & $\Z$ & $\Z$  \\ 
\hline
5 & $\Z$ & 0 & 0 & 0 \\ 
\hline
6 & $\Z$ & 0 & 0 & 0 \\ 
\hline
7 & $\Z \oplus \Z$ & 0 & 0 & 0\\ 
\hline
\end{tabular} }
\caption{\small{\emph{The Grothendieck groups $\mathfrak{M}_n$ and Bott-periodicity in a few incarnations. The isomorphism $\mathfrak{M}_{n-1}/\iota^*\mathfrak{M}_n \cong \widetilde{\KO}(\mathbb{S}^n)$ comes from \cite{ABS64}.}}}
\label{Table:KTheoryBottPeriodicity}
\end{table}

Along the way in Theorem \ref{Thm:ABSMaximallyTransverse}, we proved something even stronger than homotopic non-triviality. 
Since the map we constructed represents a stably non-trivial vector bundle, the inclusion of $f\colon \mathbb{S}^{n-1} \hookrightarrow \mathrm{SO}(d) \hookrightarrow \mathrm{SO}(d+k)$ will remain homotopically non-trivial for all non-negative integers $k \geq 0$, even when $n+1 \ge d$, which can happen when $n\le 8$.
Thus, we obtain the following: 
\begin{corollary}\label{Cor:MaxlTransversalityABS}
Let $n \geq 2$ be an integer and $d= d(n)$ as in \eqref{SpinModuleDimensions}. The transverse $(n-1)$-spheres $\Lambda$ in $\Flag(\R^{d-1,d})$ and $\Flag(\R^{d,d})$ built in Theorem \ref{Thm:TransverseSpinorSphere} include maximally transversely in the $A_{2d-2}$, $A_{2d-1}$, $A_{2d}$, $B_{d}$, $D_{d+1}$ full flag manifolds. 
\end{corollary}

\begin{proof}
Note that for $d> n+1$, if a map $f:\mathbb{S}^{n-1} \rightarrow \mathrm{SO}(d)$ is homotopically non-trivial, then so too is $\iota \circ f$, for $\iota: \mathrm{SO}(d) \hookrightarrow \mathrm{SO}(d+k)$ for any $k \geq 0$. Indeed, this follows from the fact that the inclusion $\iota$ induces an isomorphism 
$\iota_{*}:\pi_{n-1}(\mathrm{SO}(d)) \hookrightarrow \pi_{n-1}(\SO(d+k))$ for any $k \geq 0$ because $d>n+1$ means we are in the stable range. 
The long exact sequence of homotopy groups from the principal bundle fibration $\mathrm{SO}(m) \rightarrow \mathrm{SO}(m+1) \rightarrow \mathbb{S}^{m}$ verifies this claim. 

Now, we relate the inclusions to transversality-preserving embeddings of flag manifolds. We need only consider (equivariant) transversality-preserving embeddings of full flag manifolds in the following cases from Example \ref{Ex:GeneralChains}: $(B_{n}, A_{2n})$, $(D_{2n}, A_{4n-1})$, $(D_{2n},B_{2n})$. In particular, transversality of the inclusion of $\Lambda$ follows immediately. The remainder of the proof handles maximal transversality. 

Suppose $\tau\colon \Flag(G) \hookrightarrow \Flag(G')$ is a transversality-preserving embedding of full flag manifolds induced by an inclusion $\iota\colon G \hookrightarrow G'$, with corresponding maximal compact subgroups $K$ and $K'$, respectively. 
We obtain a commutative diagram 
\[\begin{tikzcd}
	K \arrow{r}{{\iota|_{K}}} \arrow[swap]{d}{\pi} & K' \arrow{d}{\pi} \\
	\Flag(G) \arrow{r}{\tau} & \Flag(G').
	\end{tikzcd} \]
Applying this reasoning in the three aforementioned cases, we see the  map $\iota|_K$ is the (standard) block-diagonal inclusions $\SO(n)\times \SO(n+1) \hookrightarrow \SO(2n+1)$, $\SO(2n) \times \SO(2n) \hookrightarrow \SO(4n)$, $\SO(2n) \times \SO(2n) \hookrightarrow \SO(2n)\times \SO(2n+1)$, respectively. 

Now, write $\Lambda= \mathrm{image}(\xi_n)$, where
$\xi_n(x)= f(x) \cdot F$, for $f \colon \mathbb{S}^{n-1} \rightarrow \SO(d+\epsilon) \times \SO(d)$ and $\epsilon \in \{-1,0\}$ in the $B_{d-1}, D_d$ cases, respectively.
In these cases, $f\colon \mathbb{S}^{n-1} \rightarrow \SO(d+\epsilon)\times \SO(d)$ obtains the form $f = (\mathrm{id}, \overline{f})$, for $\overline{f}\colon\mathbb{S}^{n-1} \rightarrow \SO(d)$, and so we need only consider the map $\overline{f}$ to the second factor of the maximal compacts of $\SO_0(d-1,d)$ and $\SO_0(d,d)$, respectively. 
In each of the three respective cases, the claim follows by the fact that the standard reducible inclusions $\SO(d+1) \hookrightarrow \SO(2d+1)$, $\SO(2d) \hookrightarrow \SO(4d)$ and $\SO(2d) \hookrightarrow \SO(2d+1)$ induce isomorphisms on $\pi_{n-1}$. 
In particular, $\Lambda$ includes homotopically non-trivially in the desired full flag manifolds. 
\end{proof}

We observe an important corollary that will allow us to build more maximally transverse spheres with a direct sum construction later on.  
\begin{corollary}\label{Cor:ABSMaximalDirectSum}
Suppose $n\geq 4$ satisfies $n \in \{0,1,2,4\} \bmod 8$ and $\eta: \Cl(n) \rightarrow \End_{\R}(S)$ is an irreducible real representation. Define the map $f: \mathbb{S}^{n-1} \rightarrow \mathrm{O}(d)$ by $f(x) = \eta(x)$. For $j \in \Z_+$, the map $f_{\oplus 2j+1}:\mathbb{S}^{n-1} \rightarrow O((2j+1)d)$ by $f_{\oplus 2j+1} = \mathrm{diag}(f,f,\dots, f)$ is homotopically nontrivial. 
Moreover, if $n \equiv 0\bmod 4$, then $[f_{\oplus j}] $ is homotopically nontrivial for any positive integer $j$. 
\end{corollary}

\begin{proof}
Set $k\coloneqq 2j+1$. 
Observe that $f_{\oplus k} $ factors via $f_{\oplus k} = i_k \circ \Delta \circ f$, where $\Delta$ and $i_k$ are 
\[ \begin{cases}
 \Delta: \mathrm{O}(d) \rightarrow \prod_{i=1}^k \mathrm{O}(d), \; &A \mapsto (A,A,\dots, A)\\
 i_k: \prod_{i=1}^k \mathrm{O}(d)\rightarrow \mathrm{O}(kd), \;&(A_1, \dots,A_k) \mapsto \mathrm{diag}(A_1,A_2,\dots, A_k). 
\end{cases}\]
Identifying $\pi_{n-1}\left(\prod_{i=1}^k \mathrm{O}(d) \right) \cong \prod_{i=1}^k\pi_{n-1}(\mathrm{O}(d))$, the induced map 
\[ (i_k)_*: \prod_{i=1}^k\pi_{n-1}(\mathrm{O}(d)) \rightarrow \pi_{n-1}(O(kd))\] 
is given by $([f_1],[f_2], \dots, [f_k]) \mapsto [\mathrm{diag}(f_1, \dots, f_k)]$. Next, let $g_k=\mathrm{diag}(\underbrace{I, I, \dots, I}_{k-1}, f,I,\dots I)$. 
Then observe that $(i_k)_*([g_\ell] )= (i_k)_*([g_m]) \in \pi_{n-1}(\mathrm{O}(d))$ for any indices $1 \leq \ell,m\leq k$. We use additive notation in the abelian group $\pi_m(X)$ when $m \geq 2$. 
Set $x \coloneqq(i_k)_*[g_1]$. Thus, $\Delta_*([f]) = \sum_{i=1}^k [g_i] = kx$. There is a canonical surjection $\pi_{n-1}(\mathrm{O}(d)) \twoheadrightarrow \pi_{n-1}(\mathbf{O})\cong \widetilde{\KO}(\mathbb{S}^n)$, under which $x \neq 0 \in \widetilde{\KO}(\mathbb{S}^n)$ by \cite{ABS64}, as in Theorem \ref{Thm:ABSMaximallyTransverse}. Thus, $x \neq 0 \in \pi_{n-1}(\mathrm{O}(d))$ as well. Now, by Table \ref{Table:KTheoryBottPeriodicity}, we have $\pi_{n-1}(\mathrm{O}(d))\in \{\Z, \Z_2\}$ since $n \in \{0,1,2,4\} \bmod 8$. In particular, $ [f_{\oplus k}] = (i_k)_*(\Delta_*([f])) = k x \neq 0 $ since $k$ is odd. 

The final claim follows from the fact that $\pi_{n-1}(\mathrm{O}(d))=\Z$ if $n \equiv 0 \bmod 4$ and $d> n+1$. 
\end{proof}

We now show that stable homotopy classes in $\pi_k(\mathbf{O})$ correspond to maximally transverse spheres. 

\begin{corollary}
Let $k \geq 2$ be an integer and $0 \neq \omega \in \pi_{k-1}(\mathbf{O})$. Then $\omega =[f] $ for a representative $f: \mathbb{S}^{k-1} \rightarrow \SO(d)$ such that the orbit map $\xi: \mathbb{S}^{k-1} \rightarrow \Flag(\R^d)$ by $\xi(x) = f(x) \cdot F$, for some flag $F \in \Flag(\R^d)$, is a maximally transverse sphere in $\Flag(\R^d)$. 
\end{corollary}

\begin{proof}
By \cite{ABS64}, there is a generator $\eta: \mathbb{S}^{k-1} \rightarrow \mathrm{O}(d)$ of $\pi_{k-1}(\mathbf{O})$ such that $\eta$ is the restriction of an irreducible Clifford algebra representation $\hat{\eta}: \Cl(k) \rightarrow \End_{\R}(S)$ to $Q_+(\R^{k,0})\subset \Cl(k)$. Set $d = \dim_{\R}(S^+)$, where $S^+$ is an irreducible $\Spin(k)$-submodule. Define $f: \mathbb{S}^k \rightarrow \SO(d)$ by $\eta(xx_0)$ for $x_0 \in \R^{k,0}$ arbitrary. Note that $[f] = [\eta] \in \pi_k(\mathbf{O})$. The map $\phi^+=(\id, f): \mathbb{S}^k \rightarrow \SO(d)\times\SO(d)$ yields a maximally transverse and homotopically non-trivial $k$-sphere $\Lambda$ in $\Flag(\R^{d,d})$ by Theorem \ref{Thm:TransverseSpinorSphere} and Theorem \ref{Thm:ABSMaximallyTransverse}. 

Now, by Table \ref{Table:KTheoryBottPeriodicity}, $k \in \{0,1,2,4\} \bmod 8$. There are two subcases. Recall by Fact \ref{Fact:MaximallyTransverse} we need only show a transverse sphere is homotopically non-trivial to prove it is maximally transverse. 

\textbf{Case 1}: $\bm{k \in \{1,2\}}$ \textbf{mod 8}. Then $\pi_{k-1}(\mathbf{O}) \cong \Z_2$. Thus, $[\eta]  =[\omega] \in \pi_{k-1}(\mathbf{O})$ and then the inclusion of $(\id, f)$ to $\SO(2d)$ suffices by Example \ref{Ex:D2nA4n-1} and Corollary \ref{Cor:ABSMaximalDirectSum}. 

\textbf{Case 2}: $\bm{k \in \{0,4\}}$ \textbf{mod 8}. Then $\pi_{k-1}(\mathbf{O})\cong \Z$. In this case, write $\omega = j [\eta]$. If $j \in \Z_+$, then the map $\phi^+_{\oplus j}: \mathbb{S}^{k-1} \rightarrow (\SO(d)\times \SO(d))^j$ includes to $\SO(2dj)$ as the desired map by Corollary \ref{Cor:IsotropicDirectSumMaps} and Corollary \ref{Cor:ABSMaximalDirectSum}. 

If instead $j < 0$, there is a simple fix. 
Let $\eta^{\pm}: \Spin(k) \rightarrow \End(S^{\pm})$ be the two distinct irreducible representations of $\Cl^0(k)$. 
Define $f^{\pm}: \mathbb{S}^{k-1} \rightarrow \SO(d)$ by $f^{\pm}(x) = \eta^{\pm}(x_0x)$ as before. 
Here is the key: by \cite{ABS64}, as in the end of the proof of Theorem \ref{Thm:ABSMaximallyTransverse}, the maps 
$f^{\pm}: \mathbb{S}^{k-1} \rightarrow \SO(d)$ satisfy $[f^+] = -[f^-] \in \pi_{k-1}(\mathbf{O})$. 
Thus, trading out the original representation $\eta$ for a non-isomorphic irreducible representation of $\Cl^0(k)$ reduces to the previous case $j > 0$, completing the proof. 
\end{proof}

Since deformations preserve our criterion for maximal transversality, the deformed spheres from Section \ref{Sec:Deformations} remain maximally transverse. 

\begin{corollary}[Maximal Transversality of Deformed Spinor Spheres]
    Let $n \in \{4,8\}$ and $\Lambda \subset \Flag(\R^{d,d})$ be a transverse $(n-1)$-sphere from Theorem \ref{thm:Full37Spheres}. Then $\Lambda$ is maximally transverse in $\Flag(\R^{d,d})$ and also in $\Flag(\R^{2d})$.
\end{corollary}

\begin{proof}
Write $d= nk$ for $n \in \{4,8\}$. Let $\Lambda' \subset \Flag(\R^{d,d})$ be an $(n-1)$-sphere from Theorem \ref{thm:Full37Spheres}. The sphere $\Lambda'$ lifts to $\hat{\Lambda}' \subset K$, the maximal compact of $\SO(d,d)$, given by the map $\hat{\xi}$ in \eqref{DeformedSphereLift}. We use the notation from Theorem \ref{thm:Full37Spheres}. In particular, the maps $F_j: \mathbb{S}^{n-1} \rightarrow \SO(n)$ are each contractible, and each map $\psi_j$ is homotopic to the identity map $\id_{\mathbb{S}^{n-1}}$. Hence, $\hat{\xi}$ is homotopic to the map $(\id, \id, \dots, \id, f, f, \dots, f)$. 
By Corollary \ref{Cor:ABSMaximalDirectSum} and Fact \ref{Fact:MaximallyTransverse}, one finds $\Lambda'$ is maximally transverse in $\Flag(\R^{d,d})$ and $\Flag(\R^{2d})$. 
\end{proof}

In some low-dimensional cases, stronger results follow from work of Tsouvalas-Zhu.
While Corollary \ref{Cor:SpecialMaximalTransversality} \ref{itemTZ_c} below was already proven in Theorem \ref{Thm:G2Fibration}, the proof here rests only on \cite{TZ24}, which does not require higher homotopy groups. 

\begin{corollary}[Special Cases of Maximal Transversality]
\label{Cor:SpecialMaximalTransversality}
Let $k \in \{3,7\}$. 
\begin{enumerate}[label = (\alph*)]
    \item \label{itemTZ_a} If $\{k, k+1\} \subseteq \Theta$, then any transverse $k$-sphere in $\Flag_{\Theta}(\R^{2k+1})$ is maximally transverse. 
    \item \label{itemTZ_b} If $k \in \Theta$, then any transverse $k$-sphere in $\Iso_{\Theta}(\R^{k,k+1})$ is maximally transverse. 
    \item \label{itemTZ_c} If $\beta \in \Theta$, then any transverse 3-sphere in $\Flag_{\Theta}(\Gtwosplit)$ is maximally transverse. 
\end{enumerate}
\end{corollary}

\begin{proof}
\ref{itemTZ_a} \cite[Lemma 5.3]{TZ24} proves the result for $\Flag_{\{k,k+1\}}(\R^{2k+1})$, which completes the proof by Fact \ref{fact:TransverseProjection}. 

\ref{itemTZ_b} One uses \ref{itemTZ_a}, the transversality-preserving map $\Iso_{k}(\R^{k,k+1}) \rightarrow \Flag_{\{k,k+1\}}(\R^{2k+1})$, and Fact \ref{fact:TransverseProjection}.

\ref{itemTZ_c} Apply \ref{itemTZ_b} with the transversality-preserving map $(\Gtwosplit, \Ein^{2,3}) \rightarrow (\SO_0(3,4),\Iso_3(\R^{3,4}))$.
\end{proof}

\section{Transverse Spheres in Isotropic Flag Manifolds}  \label{Sec:IsotropicFlags}

In this section, we combine the transverse spheres found in Section \ref{Sec:SpheresWithSpinors} with a direct sum construction for isotropic flags to build transverse spheres in the full flag manifolds of types $B_{4n-1}$, $B_{4n}$, $D_{4n}$, $D_{4n+1}$, and $D_{4n+3}$, when $n \ge 1$.
More generally, we determine when a self-opposite partial flag manifold of type $B_n$ or $D_n$ contains a transverse $2$-sphere.

\subsection{Isotropic Direct Sum}

We will need to consider certain decompositions of $\R^{p,q}$ into non-degenerate subspaces. 

\begin{definition}
We call a vector space with quadratic form $(V, q)$ \textbf{strictly pseudo-Euclidean} when $(V,q)\cong \R^{n,k}$, with $n, k \geq 1$. 

Let $V \cong \R^{n,k}$ be a strictly psuedo-Euclidean vector space. We call a splitting $V = V_1 \oplus V_2$ \textbf{admissible} when $V_1$ and $V_2$ are each strictly psuedo-Euclidean.
\end{definition}

We now present the direct sum construction for isotropic flags. 
\begin{definition}\label{Defn:FlagSum}
Let $V = V_1 \oplus V_2$ be an admissible direct sum decomposition of a strictly pseudo-Euclidean vector space $V$. Suppose the signature of $V_i$ is $(p_i, q_i)$. Let  $\Iso_{\bm{n}_i}(V_i)$ be self-opposite $\SO_0(p_i, q_i)$-flag manifolds, with indices $\bm{n}_i\coloneqq(n_j^i)_{j=1}^{k_i}$.  
Take flags $F_i \in \Iso_{\bm{n}_i}(V_i)$ and define $F \coloneqq F_1 +_f F_2 \in \Iso_{\bm{m}}(V)$, with index $\bm{m}\coloneqq (n_1^1,\dots, n_{k_1}^{1}, n_{k_1}^{1}+ n^2_1, \dots, n^1_{k_1}+n^2_{k_2})$, by 
\begin{itemize}
    \item $F^{n_j^1} \coloneqq (F_1)^{n^1_{j}}$ for $1 \leq j \leq k_1$. 
    \item $F^{n_{k_1}^1 + n^2_{l}} = (F_1)^{n^1_{k_1}} \oplus (F_2)^{n^2_{l}}$ for $1 \leq l \leq k_2 $.
\end{itemize}
We shall write $\bm{m}\coloneqq \bm{n}_1 +_f \bm{n_2}$ for the new index set, so that  $F \in \Iso_{\bm{m}}(V)$. 
\end{definition}

One can regard Definition \ref{Defn:FlagSum} as producing a map 
\begin{align}\label{IsotropicDirectSumMapping}
+_f: \Iso_{\bm{n}_1}(V_1) \times \Iso_{\bm{n}_2}(V_2) \rightarrow \Iso_{\bm{n}_1+_f \bm{n}_2}(V_1\oplus V_2). 
\end{align} 
We emphasize a special case that is permitted under the definition. 
\begin{remark}
In Definition \ref{Defn:FlagSum}, it is crucial that the index $p^+$ is permitted to appear in $\bm{n}_i$ in the case that $(p_i,q_i) = (2p,2p)$, to permit the self-opposite flag manifold $\Iso_{p}^+(\R^{2p,2p})$. In this case, Definition \ref{Defn:FlagSum} can be used to define a direct sum of flags $\Iso_p^+(\R^{p,p}) \times \Iso_p^+(\R^{p,p}) \rightarrow \Iso_{2p}^+(\R^{2p,2p})$. 
\end{remark}

The following lemma motivates Definition \ref{Defn:FlagSum} and is the crux of the inductive construction of the transverse spheres of interest. Note that if $(G_i, \mathcal{F}_{\Theta_i})$ are flag manifolds, for $i \in \{1,2\}$, then transversality in $(G_1\times G_2, \mathcal{F}_{\Theta_1} \times \mathcal{F}_{\Theta_2})$ is captured by transversality in both factors. 

\begin{lemma}\label{Lem:IsotropicDirectSum}
Let $V \cong \R^{p,q}$ be a strictly pseudo-Euclidean vector space, with admissible splitting $V = V_1 \oplus V_2$ such that $\sig(V_i)=(p_i, q_i)$. Let $\Iso_{\bm{n}_i}(V_i)$ be $\SO_0(p_i,q_i)$-self-opposite partial flag manifolds. 
The map $+_f$ in \eqref{IsotropicDirectSumMapping} is an equivariant transversality-preserving embedding. 
\end{lemma}

The following proof also shows the inverse of $+_f$ is transversality-preserving when restricted to its image, though we shall not need this property. 

\begin{proof}
We use notation for indices as in Definition \ref{Defn:FlagSum}. First, observe that the $\bm{n}_1$-transversality conditions on $F_1, F_1'$ are precisely that, inside of $V_1$, we have 
\begin{align}\label{Transversality1}
 (F_1)^{n^1_j} + \left[ ((F_1')^{n^1_j})^\bot \subset V_1 \right]= V_1 . 
\end{align}
Observe that \eqref{Transversality1} is equivalent to the following in $V = V_1 \oplus V_2$: 
\[ (F_1)^{n^1_j} + \left [ ((F_1')^{n^1_j})^\bot \subset V \right]= V_1 \oplus V_2 . \] 

Next, we examine the $(n^1_{k_1}+\bm{n}_2)$-transversality conditions on $F, F'$. To this end, we make an observation. Take any pair $W,W'$ of flags in a self-opposite flag manifold $\mathcal{F}$ of $\SO_0(p,q)$.
We suppose $\mathcal{F}$ is a minimal self-opposite partial flag manifold, so $\mathcal{F}$ may be of the form $\Iso_k(\R^{p,q})$, $\Iso_p^\pm(\R^{p,p})$ with $p$ even.
In any case, $W \pitchfork W'$ if and only if $(q_V)|_{W \times W'} \colon W \times W' \rightarrow \R$ is a non-degenerate pairing. 
Moreover, if two such isotropic subspaces $W, W' < V$ are realized as direct sums by $W = W_1 \oplus W_2$ and $W' = W_1' \oplus W_2'$, with $W_i, W'_i\ <V_i$, then $W, W'$ have non-degenerate pairing if and only if $W_1, W_1'$ have non-degenerate pairing and $W_2, W_2'$ have non-degenerate pairing. From this argument, we conclude that 
\[ \big[\,(F_1)^{n^1_{k_1}}+(F_2)^{n^2_j} \,\big] \pitchfork \big[\,(F_1')^{n^1_{k_1}}+(F_2')^{n^2_j}\,\big] \iff \begin{cases}
    (F_1)^{n^1_{k_1}}\pitchfork (F_1')^{n^1_{k_1}}\\
    (F_2)^{n_j^2} \pitchfork (F_2')^{n_j^2}.
\end{cases}\, \]
It follows that $F \pitchfork F'$ if and only if $F_1 \pitchfork F_1'$ and $F_2 \pitchfork F_2'$. 

The map $+_f$ is equivariant under the map $\SO_0(p_1,q_1) \times \SO_0(p_2,q_2) \rightarrow \SO_0(p_1+p_2, q_1+q_2)$ by $(g_1, g_2) \mapsto g_1 \oplus g_2$ and evidently injective.  
\end{proof}

Recall that we write $\phi \colon X \stackrel{\pitchfork}{\longrightarrow} \mathcal{F}$ when a map $\phi$ to a self-opposite flag manifold $\mathcal{F}$ is \emph{transverse}, meaning it sends distinct points to transverse flags. 

\begin{corollary}[Direct Sum of Transverse Maps]
\label{Cor:IsotropicDirectSumMaps}
Suppose $\phi_i \colon X \stackrel{\pitchfork}{\longrightarrow} \Iso_{\bm{n}_i}(\R^{p_i,q_i})$ are  continuous transverse maps for $i \in \{1,2\}$. 
Then $\phi_1 +_f \phi_2 \colon X \rightarrow \Iso_{\bm{m}}(\R^{p_1+p_2, q_1+q_2})$ given by $(\phi_1+_f\phi_2)(x) \coloneqq \phi_1(x) +_f \phi_2(x)$ is a continuous transverse map, for $\bm{m} = \bm{n}_1 +_f \bm{n}_2$. 
\end{corollary}

Using the isotropic direct sum construction, we can prove the following theorem. 

\begin{theorem}\label{Thm:BnCase13}
Let $k \in \Z_{>0}$ and $\epsilon \in \{-1,0,1\}$. The full flag manifold $\Flag(\R^{4k,4k+\epsilon})$ contains a maximally transverse $(\rho(4k)-1)$-sphere. 
\end{theorem}

\begin{proof}
Write the Radon-Hurwitz decomposition $4k=(2 \ell+1)2^j16^m$ and let $n=2^j16^m \ge 4$. 
If $\ell = 0$, then the result follows immediately from Corollary \ref{Cor:IsotropicTransverseSphere} and Corollary \ref{Cor:MaxlTransversalityABS}.
Otherwise, suppose $\ell \in \Z_+$. 
Consider the decomposition 
    \[ \R^{4k,4k+\epsilon} = (\R^{n,n})^{2l} \oplus \R^{n,n+\epsilon} .\] 
By Corollary \ref{Cor:IsotropicTransverseSphere}, each of  $\R^{n,n}$ and $\R^{n,n+\epsilon}$ admit a transverse $(\rho(4k)-1)$-sphere, so inclusions provide transverse maps $\phi_i \colon \mathbb{S}^{\rho(4k)-1} \stackrel{\pitchfork}{\longrightarrow} \Flag_+(\R^{n,n})$, for $1 \le i \le 2 \ell$, and $\phi_{2\ell+1} \colon \mathbb{S}^{\rho(4k)-1} \stackrel{\pitchfork}{\longrightarrow} \Flag(\R^{n,n+\epsilon})$.
We let $\phi \colon \mathbb{S}^{\rho(4k)-1} \rightarrow \Flag(\R^{4k,4k+\epsilon})$ denote the direct sum of $(\phi_i)_{i=1}^{2\ell+1}$ and see that transversality of $\phi$ follows by applying Corollary \ref{Cor:IsotropicDirectSumMaps} inductively.

Now we show the spheres are maximally transverse. 
Suppose $f \colon \mathbb{S}^{\rho(4k)-1} \rightarrow \SO(n)$ is the map $f(x) = \eta(x_0x)$ corresponding to the Clifford action described in Theorem \ref{Thm:TransverseSpinorSphere}. 
Define $\hat{f} \colon \mathbb{S}^{\rho(4k)-1} \rightarrow \SO(n) \times \SO(n)$ by $\hat{f} = (\id ,f)$, so that the transverse spheres $\phi_i$ obtain the form $\phi_i(x) =\hat{f}(x) \cdot F$ for a fixed flag $F$. 
Let $\phi \colon \mathbb{S}^{\rho(4k)-1} \rightarrow \Flag(\R^{4k, 4k+\epsilon})$ again be the direct sum of $(\phi_i)_{i=1}^{2\ell+1}$. 
Since $\rho(4k)-1 \ge 3$ and $\SO(4k)\times \SO(4k+\epsilon) \rightarrow \Flag(\R^{4k,4k+\epsilon})$ is a covering map, the direct sum $\phi$ is homotopically nontrivial in $\pi_{\rho(4k)-1}(\Flag(\R^{4k,4k+\epsilon}))$ if and only if the lift  $\psi \colon \mathbb{S}^{\rho(4k)-1} \rightarrow \SO(4k) \times \SO(4k+\epsilon)$ of $\phi$ is homotopically nontrivial. 
The lift $\psi$ takes a slightly different form in each case, as follows. We have 
\[ \psi \coloneqq (\hat{f}_{\oplus 2\ell}, \hat{f_{\epsilon}} ) =  
    ( \underbrace{(\id, f), (\id, f), \dots, (\id,f)}_{2\ell}, \hat{f_{\epsilon}}), \]  
where the final map $\hat{f_{\epsilon}} \colon \mathbb{S}^{\rho(4k)-1} \rightarrow \SO(n,n+\epsilon)$ takes the form 
\[ \hat{f_{\epsilon}} = \begin{cases}
    (f,\id) & \epsilon = -1 \\
    (\id,f) & \epsilon = 0\\
    (\id,\iota f) & \epsilon = 1 \\
\end{cases} ,\]
and $\iota: \SO(4k) \rightarrow \SO(4k+1)$ is the inclusion $A \mapsto \mathrm{diag}(A, 1)$. \medskip 

In each case, using Corollary \ref{Cor:ABSMaximalDirectSum}, we find $[\psi] \neq 0 \in \pi_{\rho(4k)-1}(\SO(4k) \times \SO(4k+\epsilon))$. 
\end{proof}

As a consequence of the above proof, as well as Corollary \ref{Cor:ABSMaximalDirectSum}, there is a great deal of flexibility in constructing maximally transverse spheres, both in changing the ambient flag manifold $\mathcal{F}$, and even with the flag manifold $\mathcal{F}$ fixed. The following example illustrates the latter idea. 

\begin{example}\label{Ex:R4848}
There exist maximally transverse spheres of dimensions $1,3,7, 8$ in $\mathcal{F} = \Flag(\R^{48,48})$. 
We first explain how to construct homotopically distinct maximally transverse 8-spheres. 
Let $f\colon \mathbb{S}^{8} \rightarrow \SO(16)$ be the map $f(x)= \eta(x_0x)$ associated to an irreducible representation $\eta\colon \Cl^0(9) \rightarrow \End_{\R}(\R^{16})$. Let us define $\phi^- = (\id, f)$ and $\phi^+ = (f, \id)$. 
The maps 
$\phi_{\oplus 3}^{\pm}\colon \mathbb{S}^{8} \rightarrow \SO(48)\times \SO(48)$ yield maximally transverse 8-spheres in $\mathcal{F}$ that are not homotopic. 
Indeed, the maps $\phi_{\oplus 3}^{\pm}$ are distinguished by which factor of $\pi_{8}(\SO(48)\times \SO(48)) \cong \Z_2 \oplus \Z_2$ they are nontrivial in. 

Next, we discuss 3-spheres in $\mathcal{F}$. 
Let $f: \mathbb{S}^{3} \rightarrow \SO(4)$ be the map $f(x_0x)$ associated to an irreducible representation $\eta:\Cl^0(4) \rightarrow \End_{\R}(\Ha)$. 
There is even more flexibility in this case. 
In each of the 12 factors of the direct sum $\R^{48,48} = (\R^{4,4})^{12}$, we have a choice of whether to take the map $(\id, f)$ or $(f, \id)$ to build a transverse 3-sphere in $\Flag(\R^{4,4})$. We find that for every pair $(a,b) \in \Z_{\geq 0}^2$ such that $a+b =12$, we have a different version $\phi_{a,b}$ of the construction. The resulting maximally transverse sphere in $\mathcal{F}$ will correspond the element $ (a \iota_{*}[f], b\iota_{*}[f]) \in \pi_3(\SO(48)) \times \pi_3(\SO(48)) \cong \Z \oplus \Z$.
Recall that there are two non-isomorphic spinor representations $\eta^{\pm}: \Spin(3) \rightarrow \SO(4)$. 
Since $[\eta^+] = -[\eta^-] \in \pi_3(\mathbf{O}) \cong \pi_3(\SO(48))$, the maximally transverse 3-spheres $\phi_{a,b}^{\pm}: \mathbb{S}^3 \rightarrow \mathcal{F}$ associated to $(\eta^{\pm}, a, b)$ are all homotopically distinct.

Similar considerations apply to produce maximally transverse 7-spheres, so we omit further details about this case. 
Finally, the limit set of any Hitchin representation $\pi_1S \rightarrow \SO_0(48,48)$ is positive \cite{FG06}, so it is a maximally transverse circle in $\mathcal{F}$ by \cite[Corollary 2.11]{DGR24}. \medskip 
\end{example}

Using the above ideas, we have the following extension of Theorem \ref{Thm:BnCase13}, to build maximally transverse spheres of different dimensions in the same full flag manifold. 

\begin{corollary}\label{Cor:ManyMaxlSpheres}
    Let $k \in \Z_{>0}$ and $\epsilon \in \{-1,0,1\}$. 
    The full flag manifold $\Flag(\R^{4k,4k+\epsilon})$ contains a maximally transverse $(\rho(2^j)-1)$-sphere for any $j\geq 2$ such that $2^j$ divides $4k$. 
\end{corollary}

\begin{proof}
We follow the same construction as in Theorem \ref{Thm:BnCase13}. 
Factor $4k = n m$, for integers $n$ and $m=2^j$ with $j \ge 2$. 
If $n=1$, we are done.
Otherwise, write 
\[ \R^{4k,4k+\epsilon} = (\R^{m,m})^{n-1} \oplus (\R^{m,m+\epsilon}).\]
Set $f \colon \mathbb{S}^{\rho(m)-1} \rightarrow \SO(m)$ by $f(x) =\eta(x_0x)$ as usual, for $\eta: \Cl^0(\rho(m)) \rightarrow \End_{\R}(\R^m)$ an irreducible representation. Define $f^+ = (f, \id)$ and $f^- = (\id , f)$. 
Then we build the map 
$\psi \coloneqq (f^-_{\oplus (n-1)}, f^+)\colon \mathbb{S}^{\rho(m)-1} \rightarrow \SO(4k) \times \SO(4k+\epsilon)$, following the notation from Example \ref{Ex:R4848}. 
The sphere $\Lambda $ obtained by letting $\psi$ act on a full flag $F \in \Flag(\R^{4k, 4k+\epsilon})$, realized as a direct sum of full flags from the component blocks, is then transverse by Corollary \ref{Cor:IsotropicTransverseSphere} and Corollary \ref{Cor:IsotropicDirectSumMaps}. 
It follows that the sphere is maximally transverse by Corollary \ref{Cor:ABSMaximalDirectSum}. 
\end{proof}

\subsection{Transverse Spheres in Type \texorpdfstring{$B$}{B} Flag Manifolds}\label{Subsec:BnTransversality}

In this section, we discuss the construction of transverse spheres in type $B$ flag manifolds and resolve Question \ref{q:circles maximal?} for partial flag manifolds of $\SO_0(p,p+1)$. 
We recall the notation $\Flag(\R^{p,p+1}) = \Iso_{\dbrack{p}}(\R^{p,p+1})$ and $\mathcal{F}_{\Theta} \cong \Iso_{\Theta}(\R^{p, p+1})$ for $\Theta \subseteq \dbrack{p}$ from Subsection \ref{Subsec:TypeBFlags}.

\begin{theorem}[Transverse spheres in $B_p$-flag manifolds]
  \label{thm:BnTransverseSpheres}
  Let $p \geq 1$ and $\Theta \subseteq \dbrack{p}$.
  Suppose
  \begin{enumerate}[label = (\alph*)]
      \item \label{item:Bp p 0 mod 4}$p \equiv 0$ mod 4. There exists a maximally transverse $(\rho(2^j)-1)$-sphere in $\Flag(\R^{p,p+1})$ for any positive integer $2^j$, with $j \geq 2$, that divides $p$. 
      \item \label{item:Bp p 1 mod 4} $p \equiv 1$ mod 4. If $p \in \Theta$, then every transverse circle in $\Iso_{\Theta}(\R^{p,p+1})$ is locally maximally transverse. Otherwise, $\Theta \subseteq \dbrack{p-1}$, and there exists a transverse $(\rho(p-1)-1)$-sphere in $\Iso_{\Theta}(\R^{p,p+1})$. 
      \item \label{item:Bp p 2 mod 4}$p \equiv 2$ mod 4. If $p \in \Theta$, then every transverse circle in $\Iso_{\Theta}(\R^{p,p+1})$ is locally maximally transverse. Otherwise, $\Theta \subseteq \dbrack{p-1}$, and there exists a transverse 2-sphere in $\Iso_{\Theta}(\R^{p,p+1})$.
      \item \label{item:Bp p 3 mod 4}$p \equiv 3$ mod 4. There exists a maximally transverse $(\rho(2^j)-1)$-sphere in $\Iso_{\Theta}(\R^{p,p+1})$ for any positive integer $2^j$, with $j \geq 2$, that divides $p+1$. 
  \end{enumerate} 
  \end{theorem}

We emphasize that in items \ref{item:Bp p 0 mod 4} and \ref{item:Bp p 3 mod 4}, we construct transverse spheres of dimension at least $3$ in the full flag manifold $\Iso_{\dbrack{p}}(\R^{p,p+1})$. 

\begin{proof}
In the case of $\SO_0(1,2) \cong \PSL(2,\R)$, the full flag manifold is $\Ein^{0,1}$, in which any two distinct points are transverse. \medskip

Now suppose $p \geq 2$ and write $p =4k+\epsilon$, with $\epsilon \in \{0,1,2,3\}$.
Cases \ref{item:Bp p 0 mod 4} and \ref{item:Bp p 3 mod 4} follow immediately from Corollary \ref{Cor:ManyMaxlSpheres}, corresponding to $\epsilon =+1$ and $\epsilon = -1$, respectively. 
For the `otherwise' claim of \ref{item:Bp p 1 mod 4}, \ref{item:Bp p 2 mod 4}, we need only produce transverse spheres in $\Iso_{\dbrack{p-1}}(\R^{p,p+1})$ by Fact  \ref{fact:TransverseProjection}. \medskip

Case \ref{item:Bp p 1 mod 4}: By \cite[Theorem A]{KT24} and \cite[Lemma 2.7]{DGR24}, every transverse circle in $\Iso_p(\R^{p,p+1})$ is locally maximally transverse. 
The remaining claim follows from Theorem \ref{Thm:BnCase13}. 
Indeed, we use the transversality-preserving embedding from Example \ref{Ex:D2B2n}, along with modifications to the examples in Section \ref{Subsec:ConcreteModelsABD} to produce the following chain of transversality-preserving maps:
\[ \Flag_+(\R^{4k,4k}) \hookrightarrow \Iso_{\llbracket 4k \rrbracket}(\R^{4k,4k+1}) \hookrightarrow\Iso_{\dbrack{4k}}(\R^{4k+1,4k+1}) \hookrightarrow \Iso_{\dbrack{4k}}(\R^{4k+1,4k+2}).  \]  Post-composing the transverse $(\rho(4k)-1)$-sphere in $\Flag_+(\R^{4k,4k})$ from Theorem \ref{Thm:BnCase13} with the above chain yields the result. 

Case \ref{item:Bp p 2 mod 4}: By \cite[Theorem A]{KT24} and \cite[Lemma 2.7]{DGR24}, every transverse circle in $\Iso_p(\R^{p,p+1})$ is locally maximally transverse.
Now we construct transverse $2$-spheres by induction.
In the base case $p = 2$ we have $X_1\coloneqq \Ein^{0,2} \stackrel{\pitchfork}{\hookrightarrow} \Ein(\R^{2,3})$.
In the general case, we consider a transverse $3$-sphere in $\Flag(\R^{4k,4k})$ from case \ref{item:Bp p 0 mod 4} and choose a $2$-sphere $X_2$ inside it.
We split $\R^{4k+2, 4k+3}= \R^{4k, 4k} \oplus \R^{2,3}$ and apply Lemma \ref{Lem:IsotropicDirectSum} to form the transverse map $X_1 +_fX_2 \colon \mathbb{S}^2 \stackrel{\pitchfork}{\hookrightarrow} \Iso_{\dbrack{p-1}}(\R^{p,p+1})$.
\end{proof}

\begin{remark}\label{remark: better than 2spheres}
It perhaps seems striking in Theorem \ref{thm:BnTransverseSpheres} \ref{item:Bp p 2 mod 4} that we can achieve only a transverse 2-sphere in general, in comparison to the neighbor cases. 
While we can build a $(\rho(4k)-1)$-sphere in $\Iso_{\dbrack{p-2}}(\R^{p,p+1})$ when $p\equiv 2 \bmod 4$, we can do no better than a transverse 2-sphere in $\Iso_{\dbrack{p-1}}(\R^{p,p+1})$ with our methods. 
Indeed, to apply our direct sum approach, we must remove some combination of $\R^{4n-1,4n}, \R^{4n,4n-1}, \R^{4n,4n+1}, \R^{4n+1,4n}$, or $\R^{2n+1,2n+1}$ blocks to kill off the extra parity in $p \bmod 4$, along with some $\R^{4n,4n}$ blocks, but the former blocks prevent the transverse sphere from reaching all the way to $\Iso_{\dbrack{p-1}}(\R^{p,p+1})$. 
\end{remark}

\subsection{Transverse Spheres in Type \texorpdfstring{$D$}{D} Flag Manifolds}\label{Subsec:DnTransversality}

In this section, we resolve Question \ref{q:circles maximal?} for partial flag manifolds of $\SO(p,p)$, $p \ge 4$. Since $\SO(3,3)$ is locally isomorphic to $\SL(4,\R)$, we handle that case in a later section. \medskip
 
Before the theorem, we recall some notation.
Here, we label the simple roots of $\mathfrak{so}(n,n)$ according to the Dynkin diagram as $\Delta = \{ 1, 2, \dots, p-2\} \cup \{p^+, p^-\}$ as in Figure \ref{fig:DynkinDp} so that $\mathcal{F}_{\{i\}} \cong \Iso_{i}(\R^{p,p})$ and $\mathcal{F}_{\{p^{\pm}\}} \cong \Iso_{p}^{\pm}(\R^{p,p})$ are the flag manifolds associated to maximal parabolic subgroups, as discussed in detail in Subsection \ref{Subsec:TypeDFlags}. 
We have the identification $\Flag(\R^{p,p}) \cong \Iso_{\dbrack{p-1}}(\R^{p,p})$ and $\mathcal{F}_{\Delta \backslash \{p^+,p^-\}} \cong \Iso_{\dbrack{p-2}}(\R^{p,p}) $.
Recall also that each $\Iso_{p}^{\pm}(\R^{p,p})$ is self-opposite if and only if $p$ is even.  

\begin{theorem}[Transverse spheres in $D_p$ flag manifolds]\label{Thm:DnTransverseSpheres}
Let $p \geq 4$. 
Suppose:
    \begin{enumerate}[label =(\alph*)]
    \item \label{item:Dp 0} $p \equiv 0 \bmod 4$. There is a maximally transverse $(\rho(2^j)-1)$-sphere in the full flag manifold $\Iso_{\dbrack{p-1}}(\R^{p,p})$ for any integer $2^j$, with $j\geq 2$, that divides $p$. 
    \item \label{item:Dp 1} $p \equiv 1 \bmod 4$. There is a maximally transverse $(\rho(2^j)-1)$-sphere in the full flag manifold $\Iso_{\dbrack{p-1}}(\R^{p,p})$ for any integer $2^j$, with $j\geq 2$, that divides $p-1$. 
    \item \label{item:Dp 2} $p \equiv 2 \bmod 4$. If $\Theta \subseteq \Delta$ satisfies $\{p^+,p^-\} \cap \Theta \neq \emptyset$, then every transverse circle in $\mathcal{F}_{\Theta}$ is locally maximally transverse. Otherwise, $\Theta \subseteq \dbrack{p-2}$, and there is a transverse $(\rho(p-2)-1)$-sphere in $\Iso_{\Theta}(\R^{p,p})$.
    \item \label{item:Dp 3}$p \equiv 3 \bmod 4$. There is a maximally transverse $(\rho(2^j)-1)$-sphere in the full flag manifold $\Iso_{\dbrack{p-1}}(\R^{p,p})$ for any integer $2^j$, with $j \geq 2$, that divides $p+1$, as long as $2^j \neq p+1$. 
    \end{enumerate}
  \end{theorem}

\begin{proof}
\ref{item:Dp 0} 
The statement is immediate from Corollary \ref{Cor:ManyMaxlSpheres}, with $p=4k$ and $\epsilon = 0$. \medskip

\ref{item:Dp 1} 
We apply Corollary \ref{Cor:ManyMaxlSpheres} with $p-1=4k$ and $\epsilon=1$ to obtain a maximally transverse $(\rho(2^j)-1)$-sphere $\Lambda$ in $\Iso_{\dbrack{4k}}(\R^{4k,4k+1})$, which includes transversely in $\Iso_{\dbrack{4k}}(\R^{4k+1,4k+1})$ by Example \ref{Ex:BnDn+1}. 
The sphere $\Lambda$ is homotopically non-trivial by Corollary \ref{Cor:ManyMaxlSpheres}. 
Since $\rho(2^j)-1 < 4k$, the inclusion $\iota \colon \SO(4k) \hookrightarrow \SO(4k+1)$ induces an isomorphism on $\pi_{\rho(2^j)-1}$ and the included sphere in $\Iso_{\dbrack{4k}}(\R^{4k+1,4k+1})$ is still homotopically non-trivial, hence maximally transverse. \medskip 

\ref{item:Dp 2}  When $p \equiv 2$ mod 4, every transverse circle in $\mathcal{F}_{\Theta}$ is locally maximally transverse if $\Theta \cap \{p^+, p^-\} \neq \emptyset$ by \cite[Theorem A]{KT24} and \cite[Lemma 2.7]{DGR24}. 
Case \ref{item:Dp 0} and the inclusion of $\Flag_+(\R^{4k,4k}) \stackrel{\pitchfork}{\hookrightarrow} \Iso_{\dbrack{4k}}(\R^{4k+2,4k+2})$ implies the remaining claim of \ref{item:Dp 2} by Fact \ref{fact:TransverseProjection}. \medskip 

\ref{item:Dp 3} Suppose $p \equiv 3\bmod 4$ and that $m = 2^j$ properly divides $p+1$. Write $p+1= nm$, for some positive integer $n \geq 2$. We break $\R^{p,p}$ into blocks as follows: 
\[\R^{p,p} = (\R^{m,m})^{n-2} \oplus \R^{m-1, m} \oplus \R^{m,m-1}.\]
We follow the construction of Corollary \ref{Cor:ManyMaxlSpheres}. 
Define $f: \mathbb{S}^{\rho(m)-1} \rightarrow \SO(m)$ by $f(x)=\eta(x_0x)$ for an irreducible representation $\eta: \Cl^0(\rho(m)) \rightarrow \End_{\R}(\R^m)$ and we use the notation $f^{\pm}$ as in Corollary \ref{Cor:ManyMaxlSpheres}. Define the map $\psi: \mathbb{S}^{\rho(m)-1} \rightarrow \SO(p) \times \SO(p)$ by $\psi = (f^-_{\oplus n-2}, f^-, f^+)$. The sphere $\Lambda$ obtained by taking the `orbit' of $\psi$ on a flag $F \in \Flag(\R^{p,p})$, realized as a direct sum of full flags in the component subspaces, is transverse by Corollary \ref{Cor:IsotropicDirectSumMaps} and also maximally transverse by similar reasoning as in Corollary \ref{Cor:ManyMaxlSpheres}.  
\end{proof}

\section{Transverse Spheres in Type \texorpdfstring{$A$}{A} Flag Manifolds}\label{Sec:AnTransversality}

In this section, we construct transverse spheres in self-opposite partial flag manifolds associated to $\SL(d,\R)$.
We recall that, for $\Theta \subseteq \dbrack{d}$, the partial flag manifold $\Flag_\Theta(\R^d)$ is self-opposite exactly when $\Theta$ is symmetric, i.e., when $k \in \Theta$ implies $d-k \in \Theta$, see Subsection \ref{Subsec:TypeAFlags}.

\subsection{General Direct Sum of Flags}

We first need to generalize the flag direct sum techniques of the previous section. 
The following lemma provides the building blocks for the construction. 

\begin{lemma}[(Basic) Flag Direct Sum]\label{Lem:BasicFlagDirectSum}
Suppose $V = V_1 \oplus V_2$ with $\dim(V_i) = d_i$. For any integers $1 \leq n_i\leq \lfloor \frac{d_i}{2} \rfloor$, with $i \in \{1,2\}$, there is an
equivariant transversality-preserving map 
\[
    +_f: \Flag_{\{n_1, \,d_1-n_1\}}(V_1)\times \Flag_{\{n_2, d_2-n_2\} }(V_2)\rightarrow \Flag_{\{n_1+n_2, \,d_1+d_2-n_1-n_2\}}(V_1 \oplus V_2), \]
where $\mathcal{F} \coloneqq F_1+_f F_2$ is given by 
\[ \begin{cases}
        \mathcal{F} ^{n_1+n_2} &= (F_1)^{n_1}+(F_2)^{n_2}\\
    \mathcal{F}^{(d_1+d_2)-(n_1+n_2)} &= (F_1)^{d_1-n_1} +(F_2)^{d_2-n_2}.
    \end{cases} \] 
Moreover, $(F_1+_f F_2) \pitchfork (F_1' +_f F_{2}')$ if and only if $F_1 \pitchfork F_1'$ and $F_2 \pitchfork F_2'$.
\end{lemma}

\begin{proof}
Note that if $n_i= d_i-n_i$, then $\Flag_{\{n_i, d_i-n_i\} }(V_i) = \Gr_{n_i}(V_i)$, which presents no problem for the following argument. 
Now, define $\mathcal{F} \coloneqq F_1+_f F_2$ and $\mathcal{F}
'\coloneqq F_1'+_f F_2'$. 
The lemma follows from the following re-organization of the flag subspaces: 
\begin{align*}
    \mathcal{F}^{n_1+n_2}+(\mathcal{F}')^{d_1+d_2-n_1-n_2}&= \left [(F_1)^{n_1}+(F_2)^{n_2} \right]+ \left [(F_1')^{d_1-n_1}+(F_2')^{d_2-n_2} \right]\\
     &= \left[(F_1)^{n_1}+(F_1')^{d_1-n_1} \right] + \left[ (F_2)^{n_2} +(F_2')^{d_2-n_2}\right] .
\end{align*}
The same considerations apply with the roles of $\mathcal{F}$ and $\mathcal{F}'$ reversed. 
It follows that $\mathcal{F} \pitchfork \mathcal{F}'$ if and only if $F_1 \pitchfork F_1'$ and $F_2 \pitchfork F_2'$. 
\end{proof}

As a consequence of the previous lemma, we may define a direct sum operation for the pair of flag manifolds $ (\Flag_{\bm{n}_1}(V_1), \Flag_{\bm{n}_2}(V_2))$, for symmetric indices $\bm{n}_j$. 
The method we use here is a bit more flexible than the one we used in the isotropic case. 
The flexibility in the general direct sum construction can be pinned down by a choice of certain lattice path as follows. 

\begin{definition}\label{Defn:FlagPath}
Let $n, m$ be positive integers. 
We call a $\Z^2$-lattice path $(\zeta_t)_{t=0}^{n+m}$ from $(0,0)$ to $(n, m)$ an \textbf{$(n,m)$-flag path} when 
\begin{itemize}
    \item $\zeta$ is symmetric about the diagonal from $(0,0)$ to $(n, m)$. 
    \item $\zeta(t+1)-\zeta(t)\in \{ (1,0), (0,1) \}$, for $0 \leq t \leq n+m-1$. 
\end{itemize}
\end{definition}

Note that that while the indices in $\zeta$ always increase by exactly one in either $x$ or $y$, this does not necessarily mean we increment the dimension of subspaces by one; see Example \ref{Ex:Dim11Flags}.\medskip

\begin{remark}
Suppose $(n,m) \in \Z_+^2$ is a pair of positive integers. One can easily verify that an $(n,m)$ flag path exists if and only if $n$ or $m$ is even. 
\end{remark}

\begin{remark}
In defining a direct sum of flag manifolds $\Flag_{\bm{n}_i}(V_i)$, one might wish to increment the indices of both flags simultaneously. Allowing diagonal steps by (1,1) in a flag path accommodates this wish. With this greater flexibility, flag paths now exist when $n=m = 2k+1$ are both odd and equal. However, for $n \neq m$, the use of diagonal (1,1) steps only produces flags which are projections of the flags defined by flag paths in Lemma \ref{Lem:DeluxeFlagSum} with the current definition.  
\end{remark}

The $(n,m)$-flag paths are defined precisely to carry the combinatorics of direct sum constructions. Indeed, the following result shows every $(k_1,k_2)$-flag path yields a direct sum construction for the flag manifolds $\Flag_{\bm{n_1}}(\R^{d_1}), \Flag_{\bm{n}_2}(\R^{d_2})$, where $k_i = |\bm{n}_i|$.

\begin{lemma}[(Deluxe) Flag Direct Sum]\label{Lem:DeluxeFlagSum}
Suppose $V = V_1 \oplus V_2$ is a direct sum of vector spaces, with $\dim(V_i) = d_i$ and $\dim(V) =d$. 
Let $\bm{n}_i = (n^i_j)_{j=0}^{k_i}$ be $d_i$-symmetric multi-indices. For any $(k_1, k_2)$-flag path $\zeta_t =(x_t, y_t)$, there is an equivariant transversality-preserving embedding 
\begin{align}\label{GeneralFlagDirectSumMapping}
    +_{f,\zeta}: \Flag_{\bm{n}_1}(V_1) \times \Flag_{\bm{n}_2}(V_2) \rightarrow \Flag_{\bm{n}_1+_{f, \zeta} \,\bm{n}_2}(V_1 \oplus V_2),
\end{align}
where the symmetric index $\bm{n} \subset \dbrack{d}$ and output flag $\mathcal{F} \coloneqq F_1+_{f,\zeta} F_2$ are as follows: 
\begin{align}
    n_t &\coloneqq n^1_{x_t} + n^2_{y_t} \\
    \mathcal{F}^{n_t} &\coloneqq(F_1)^{n^1_{x_t}} + (F_{2})^{n^2_{y_t}}.
\end{align}
Moreover, $(F_1+_{f, \zeta}F_2 )\pitchfork (F_1'+_{f, \zeta}F_2' )$ if and only if $F_1 \pitchfork F_2$ and $F_1'\pitchfork F_2'$. 
\end{lemma}

\begin{proof}
Note that the set $\bm{n} \subset \dbrack{d}$ is symmetric because $\zeta$ is symmetric about the diagonal. 
The subspaces are nested, $\mathcal{F}^{n_{t}} \subsetneq \mathcal{F}^{n_{t+1}}$, because the indices $x_t, y_t$ are monotonically increasing, with $z_t \coloneqq x_t+y_t$ strictly monotonically increasing in $t$. 
Thus, $\mathcal{F} \in \Flag_{\bm{n}}(V)$. 

For the transversality, set $\mathcal{F} \coloneqq F_1 +_{f, \zeta} F_2$ and $\mathcal{F}' \coloneqq F_1' +_{f, \zeta} F_2'$. Then $\mathcal{F} \pitchfork \mathcal{F}'$ if and only if for each $j \in \bm{n}$, $\mathcal{F}^j + (\mathcal{F}')^{d-j} = \R^d$. 
Now, for each $j \in \bm{n}$, writing uniquely $j = x_t+y_t$, then (the proof of) Lemma \ref{Lem:BasicFlagDirectSum} shows that $\mathcal{F}^j + (\mathcal{F}')^{d-j} = \R^d$ if and only if $(F_1)^{x_t} + (F_1')^{d_1-x_t} = V_1$ and $(F_2)^{y_t} + (F_2')^{d_2-y_t} = V_2$. 
The claim on equivalence of transversality follows. 
The equivariance is obvious, and injectivity follows since the projection $\pi_{V_i}: V \rightarrow V_i$ with respect to $V=V_1\oplus V_2$ of $F_1 +_{f,\zeta}F_2$ recovers the component subspaces of $F_i$. 
\end{proof}

\begin{corollary}\label{cor:sum of An flag maps}
Suppose $\phi_i: X \rightarrow \Flag_{\bm{n}_i}(V_i)$ are continuous transverse maps. Let $\zeta$ be a  $(k_1,k_2)$-flag path, and define $\bm{n} $ as in the previous lemma. Then $\phi: X \rightarrow \Flag_{\bm{n}}(V_1\oplus V_2)$ given by 
$\phi(x) \coloneqq \phi_1(x)+_{f,\zeta} \phi_2(x)$ is a continuous transverse map. 
\end{corollary}

See Figure \ref{Fig:LatticePath11} for an example of a $(7,2)$-flag path. 
This flag path is used in Example \ref{Ex:Dim11Flags} to apply Lemma \ref{Lem:DeluxeFlagSum}. 

\subsection{Constructing the Spheres}\label{Subsec:AnTransversality}

We now resolve Question \ref{q:circles maximal?} for self-opposite partial flag manifolds associated to $\SL(d,\R)$.
Recall that $\rho(d)$ is the Radon-Hurwitz number of $d \in \Z_+$ from \eqref{RadonHurwitz}. 

\begin{theorem}[Transverse spheres in $A_{d-1}$ flag manifolds]
\label{Thm:AnTransverseSpheres}
    Fix an integer $d \ge 2$. 
    Write $d= 8k+\epsilon$, for $k \geq 1$ and $\epsilon \in \{-1,0,1,2, \dots, 6\}$. Suppose 
     \begin{enumerate}[label=(\alph*)]
        \item \label{item:An 701} $\epsilon \in \{-1,0,1\}$. Write $4k =n2^j$ for some integers $j\geq 2$, $n \geq 1$, such that $j \in \{2,3\} \bmod 4$ or $n$ is odd. Then there is a maximally transverse $(\rho(2^j)-1)$-sphere in $\Flag(\R^d)$. 
        \item \label{item:An 2} $\epsilon=2$. If $\frac{d}{2} \in \Theta$, then every transverse circle in $\Flag_{\Theta}(\R^d)$ is maximally transverse. Otherwise, $\Theta$ is a symmetric subset of $\dbrack{d} \setminus \{\frac{d}{2}\}$ and there exists a transverse $(\rho(4k)-1)$-sphere in $\Flag_{\Theta}(\R^d)$. 
        \item \label{item:An 345} $\epsilon \in \{3,4,5\}$. Every transverse circle in $\Flag(\R^d)$ is locally maximally transverse. Otherwise, $\Theta$ is a proper symmetric subset of $\dbrack{d}$, and there exists a transverse $2$-sphere in $\Flag_\Theta(\R^d)$.
        \item \label{item:An 6} $\epsilon =6$. If $\frac{d}{2} \in \Theta$, then every transverse circle in $\Flag_{\Theta}(\R^d)$ is maximally transverse. Otherwise, $\Theta$ is a symmetric subset of $\dbrack{d} \setminus \{\frac{d}{2}\}$ and there exists a transverse $(\rho(4k+4)-1)$-sphere in $\Flag_{\Theta}(\R^d)$.
    \end{enumerate}
\end{theorem}

The proof essentially combines Theorem \ref{Thm:BnCase13}, Example \ref{Ex:BnA2n}, Example \ref{Ex:DnA2nMinusn}, and Corollary \ref{cor:sum of An flag maps} to construct transverse $m$-spheres with $m \ge 2$ in the various cases. 
The maximal transversality uses the ideas of Corollary \ref{Cor:ABSMaximalDirectSum}, but has some additional complications in casework that did not appear in the type $B$, $D$ settings. 
The remaining claim that the mentioned transverse circles are (locally) maximally transverse is established in \cite{tsouvalas2020borel} in cases \ref{item:An 2} and \ref{item:An 6} and in \cite{Dey25} in case \ref{item:An 345}. 
Before the proof, we illustrate the key ideas through two examples.

\begin{example}\label{Ex:Dim11Flags}
    Let us consider partial flag manifolds of $\R^{11}$. 
    We have a transverse $\mathbb{S}^3$ in $\Flag_{\dbrack{11} \backslash \{5,6\}}(\R^{11})$ by Theorem \ref{thm:BnTransverseSpheres}\ref{item:Bp p 1 mod 4} and Example \ref{Ex:DnA2nMinusn}.
    The inclusion $\mathrm{SL}(2,\C) \hookrightarrow \SL(4, \R)$ yields a transverse $\mathbb{CP}^1 \cong \mathbb{S}^2$ in $\Gr_2(\R^4)$ by considering the visual boundary of the sub-symmetric space of $\SL(2,\C)$. 
    With the transverse $\mathbb{S}^2$ in $\Gr_2(\R^4)$, we can construct a transverse $\mathbb{S}^2$ in any flag manifold $\Flag_{\llbracket 11 \rrbracket \setminus \{k,11-k\}}(\R^{11})$, $1 \le k \le 5$.
    Let $\eta \colon \mathbb{S}^2 \to \Flag_{\llbracket 7 \rrbracket}(\R^7)$ denote a transverse $\mathbb{S}^2$ we obtain by restricting our transverse $\mathbb{S}^3$, 
    and let $\zeta \colon \mathbb{S}^2 \to \Gr_{\{0,2,4\}}(\R^4)$ denote the transverse $\mathbb{S}^2$ mentioned above.
    Here, for notational convenience, we record the trivial $0$- and full-dimensional parts of the flags. Associated to $k$, we define $\xi_k \colon \mathbb{S}^2 \to \Flag_{\llbracket 11 \rrbracket \setminus \{k,11-k\}}(\R^{11})$ by writing $\R^{11} = \R^7 \oplus \R^{4}$ and then setting
    \[ \xi_k^i(x) \coloneqq 
        \begin{cases} 
          \eta^i(x) + \zeta^0, & 0\leq i < k \\
          \eta^{i-2}(x) +\zeta^2(x), &k < i < 11-k\\
          \eta^{i-4}(x)+\zeta^4, & 11-k < i \leq 11. 
        \end{cases}                         
    \]
   In fact, $\xi_k$ is equivalently described Lemma \ref{Lem:DeluxeFlagSum}, with $\bm{n}_1=(0,1,2,\dots, 7), \bm{n}_2 = (0,2,4)$ and $(\zeta_t)$  the $(7,2)$-flag path 
   $ (k-1)\cdot \rightarrow, \uparrow, (9-2k)\cdot \rightarrow,\uparrow, (k-1)\cdot \rightarrow $. For example, we illustrate the lattice path when $k =4$ in Figure \ref{Fig:LatticePath11}. In all cases, $\xi_k$ is transverse by Lemma \ref{Lem:DeluxeFlagSum}. 
    \begin{figure}[ht]
    \centering
    \includegraphics[width=.7\textwidth]{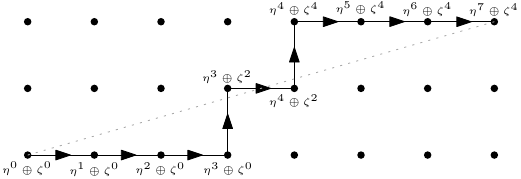}
    \caption{\emph{\small{The relevant lattice path in Example \ref{Ex:Dim11Flags} for the transverse map $\zeta_4: \mathbb{S}^2 \rightarrow \Flag_{\dbrack{11} \backslash \{4,7\}}$. The lattice path is symmetric about the diagonal. }}}
    \label{Fig:LatticePath11}
    \end{figure}

\end{example} 

\begin{example}
Consider the maximally transverse 8-sphere $\Lambda$ in $\Flag(\R^{16,16})$ from Theorem \ref{Thm:BnCase13}. Then a direct sum construction using Corollary \ref{Cor:IsotropicDirectSumMaps} builds a transverse 8-sphere in $\Flag(\R^{32, 32})$ that includes transversely in $\Flag(\R^{64})$ by Example \ref{Ex:D2nA4n-1}. If the direct sum construction is performed carefully, as in Corollary \ref{Cor:ManyMaxlSpheres}, then the resulting transverse 8-sphere in $\Flag(\R^{32,32})$ is maximally transverse. However, regardless of how the direct sum is performed, following the reasoning in Corollary \ref{Cor:ABSMaximalDirectSum}, one finds the resulting sphere in $\Flag(\R^{64})$ is homotopically trivial. It is currently unclear whether a transverse subset $\Lambda \subset \mathcal{F}_{\Theta}$ that is nullhomotopic in $\mathcal{F}_{\Theta}$ can be maximally transverse in $\mathcal{F}_{\Theta}$. So, with our current methods, we cannot guarantee the existence of a maximally transverse $8$-sphere in $\Flag(\R^{64})$. Similar phenomena occur for other transverse $(8k)$-spheres and $(8k+1)$-spheres in $\Flag(\R^{d})$ produced by even direct sums. This is the reason for the caveat in Theorem \ref{Thm:AnTransverseSpheres} \ref{item:An 701} that did not appear in Theorems \ref{thm:BnTransverseSpheres} or \ref{Thm:DnTransverseSpheres} in the $B$, $D$ cases. 
\end{example}

\begin{proof}[Proof of Theorem \ref{Thm:AnTransverseSpheres}]

We note here the considerations for $2 \leq d \leq 5$. 
When $d= 2$, the only flag manifold is the full flag manifold, which is already a circle. 
By \cite{Dey25}, for $3\leq d \leq 5$, every transverse circle in $\Flag(\R^d)$ is locally maximally transverse, which settles the case $d=3$. 
For the ``otherwise'' claim, proceed as follows. 
For $d= 4$, the maps $\Ein^{0,2} = \Ein(\R^{1,3}) \stackrel{\pitchfork}{\hookrightarrow} \Flag_{\{1,3\}}(\R^4)$ and 
$\CP^1 \stackrel{\pitchfork}{\hookrightarrow} \Gr_2(\R^4)$ provide transverse 2-spheres in the self-opposite partial flag manifolds. 
The case $d=5$ is handled similarly: 
the maps $\Ein^{0,2} \stackrel{\pitchfork}{\hookrightarrow} \Ein(\R^{2,3}) \stackrel{}{\hookrightarrow} \Flag_{\{1,4\}}(\R^5)$, and $ \CP^1 \stackrel{\pitchfork}{\hookrightarrow} \Gr_2(\R^4)\stackrel{}{\hookrightarrow} \Flag_{\{2,3\}}(\R^5)$ provide transverse 2-spheres, where the transversality-preserving map $\Gr_2(\R^4) \stackrel{}{\hookrightarrow} \Flag_{\{2,3\}}(\R^5)$ is given by $P \mapsto (P, P\oplus L)$, writing $\R^5 = \R^4 \oplus L$. \medskip 

Now, suppose $d \geq 6$. \medskip 

\ref{item:An 701} Write $m \coloneqq 2^j$. There is a  transverse $(\rho(m)-1)$-sphere in $\Flag(\R^{4k,4k+\epsilon})$ by Corollary \ref{Cor:ManyMaxlSpheres}. Post-composing with the transversality-preserving embedding $\Flag(\R^{4k,4k+\epsilon}) \stackrel{}{\hookrightarrow} \Flag(\R^{8k+\epsilon})$, for $\epsilon \in \{-1,0,1\}$, from Example \ref{Ex:GeneralChains} produces a transverse $(\rho(m)-1)$-sphere in $\Flag(\R^{8k+\epsilon})$. There are now two sub-cases. \medskip 

\textbf{Suppose $\bm{j \in \{2,3\}}$ mod 4.} Then $\pi_{\rho(m)-1}\SO(d) =\Z$ by Table \ref{Table:KTheoryBottPeriodicity} and \eqref{RadonHurwitzInversion}. Thus, any version of the direct sum construction, as in Corollary \ref{Cor:ManyMaxlSpheres} and Corollary \ref{Cor:ABSMaximalDirectSum}, will produce a maximally transverse sphere. 

\textbf{Otherwise, $\bm{j \in \{0,1\}}$ mod 4.} Again,  write $4k = 2^jn$ for some positive integer $n$. In this case, $\pi_{\rho(m)-1}(\SO(d))= \Z_2$ by \eqref{RadonHurwitzInversion} and Table \ref{Table:KTheoryBottPeriodicity}. We then see the resulting sphere, under any version of the direct sum construction, is homotopically nontrivial if and only if $n$ is odd. The result for \ref{item:An 701} follows. 

\ref{item:An 2} The proof of \cite[Theorem 1.1]{tsouvalas2020borel} establishes that a transverse circle in $\Gr_{\frac{d}{2}}(\R^d)$ is maximally transverse when $d \equiv 2 \mod 4$. 
By Fact \ref{fact:TransverseProjection}, the same holds for any $\Flag_\Theta(\R^d)$ where $\Theta$ contains $\frac{d}{2}$.
For the remaining claim, we use Theorem \ref{Thm:DnTransverseSpheres}\ref{item:Dp 1} and the transversality-preserving embedding 
$\Flag(\R^{4k+1, 4k+1})\stackrel{}{\hookrightarrow} \Flag_{\dbrack{8k+2} \backslash \{4k+1\}}(\R^{8k+2})$ of Example \ref{Ex:DnA2nMinusn}. This construction suffices by Fact  \ref{fact:TransverseProjection}. 

\ref{item:An 345} The claim about local maximal transversality of transverse circles is a special case of \cite[Corollary B]{Dey25}. 
The remaining claim is handled as follows. 
Let $k \geq 1$, $d = 8k+\epsilon$, and $\epsilon \in \{3,4,5\}$. 
To prove the claim, it suffices to consider $\Theta$ of the form $\Theta = \dbrack{d} \backslash \{j, d-j\}$, for some integer $1 \leq j \leq \lfloor \frac{d}{2} \rfloor$, by Fact \ref{fact:TransverseProjection}.

Suppose that $\epsilon \in \{3,5\}$. By Theorem \ref{Thm:BnCase13}, there is a continuous transverse map $\eta: \mathbb{S}^2 \rightarrow \Flag(\R^{8k+\epsilon-4})$.
Next, take $\zeta^2 \colon \mathbb{S}^2 \rightarrow \Gr_{2}(\R^{4})$ to be any continuous transverse map and write $\zeta = (\zeta^0, \zeta^2, \zeta^4)$, where $\zeta^0 =\{0\}, \zeta^4 =\R^4$. 

We define $\xi_{j} \colon \mathbb{S}^2 \rightarrow \Flag_{\Theta}(\R^{d})$ as follows:
\[ \xi_j^i(x) \coloneqq 
        \begin{cases} 
          \eta^i(x) + \zeta^0, &0 \leq i < j \\
          \eta^{i-2}(x)+\zeta^2(x), & j < i < d-j\\
          \eta^{i-4}(x)+\zeta^4,  & d-j < i \leq d.
        \end{cases}                         
    \]
One finds $\xi_j$ is defined by an appropriate application of Corollary \ref{cor:sum of An flag maps} and hence is transverse.   

In the case $\epsilon =4$, one can perform an analogous construction as above if $j < 4k+2$.
In the case $j=4k+2$, we proceed slightly differently. 
Decompose $\R^{8k+4} = \R^{8k} + \R^{4}$ and take $\eta \colon \mathbb{S}^2 \rightarrow \Flag(\R^{8k})$ and $\zeta \colon \mathbb{S}^2 \rightarrow \Flag_{\{1,3\}}(\R^4)$ to be continuous and transverse. 
The latter map $\eta$ can be built by taking $\mathbb{S}^2 \cong \Ein^{0,2} = \Flag(\R^{1,3}) \stackrel{\pitchfork}{\hookrightarrow} \Flag_{\{1,3\}}(\R^{4})$. 
Then we use the flag direct sum with $\bm{n}_1 = (0,1, 2, \dots, 8k), \, \bm{n}_2=(0,1,3,4)$ and the $(8k,3)$-flag path $(4k)\cdot\rightarrow , \uparrow, \uparrow, \uparrow, (4k)\cdot \rightarrow$. 
This construction yields the following continuous transverse map $\xi_{4k+2} \colon \mathbb{S}^2 \rightarrow \Flag_{\dbrack{8k+4}\backslash \{4k+2\}}(\R^{8k+4})$, given by
\[ \xi_{4k+2}^i(x) \coloneqq 
        \begin{cases} 
          \eta^i(x) + \zeta^0, & 0 \leq i < 4k+1 \\
            \eta^{4k}(x) +\zeta^{1}(x),  &  i=4k+1\\
         \eta^{4k}(x)+\zeta^{3}(x),  &  i=4k+3\\
         \eta^{i-4}(x) + \zeta^4, &  4k+3 < i \leq 8k+4.
        \end{cases}                         
    \]

\ref{item:An 6} 
When $\Theta$ contains $\{\frac{d}{2}\}$, a transverse circle in $\Flag_\Theta(\R^d)$ is maximally transverse by  \cite[Theorem 1.1]{tsouvalas2020borel}.
For the remaining claim on transverse spheres, by Fact \ref{fact:TransverseProjection}, it suffices to produce a transverse $(\rho(4k+4)-1)$-sphere in $\Flag_{\dbrack{8k+6}\backslash \{4k+3\}}(\R^{8k+6})$. 
By Theorem \ref{Thm:TransverseSpinorSphere}, there is a transverse $\rho(4k+4)-1$-sphere in $\Iso_{\dbrack{4k+2}}(\R^{4k+2,4k+4})$, which then includes transversely in $\Flag_{\dbrack{8k+6}\backslash \{4k+3\}}(\R^{8k+6})$. 
\end{proof}

\begin{remark}
For some partial flag manifolds in case \ref{item:An 345} of Theorem \ref{Thm:AnTransverseSpheres}, it is possible to obtain transverse $k$-sphere for $k > 2$, c.f.\ Remark \ref{remark: better than 2spheres}. 
However, we write only the most general result to avoid diverging into further cases.
\end{remark}

\section{Transverse Spheres in Flag Manifolds of Exceptional Groups}\label{sec:exceptional}

By abuse of notation, we let $F_4$ (resp.\ $E_6,E_7,E_8$) denote the split real (adjoint) form of Cartan-Killing type $F_4$ (resp.\ $E_6,E_7,E_8$).\medskip 

\begin{theorem}\label{thm:F4,E6,E8}
    There exists a transverse $3$-sphere in the full flag manifold associated to $F_4$ as well as the full flag manifold associated to $E_6$.
    There exists a transverse $7$-sphere in the full flag manifold associated to $E_8$.
\end{theorem}

\begin{proof}
    There is a $B_4$ root system in the $F_4$ root system, inducing a local embedding of $\SO_0(4,5)$ into $F_4$. 
    Since the root systems have the same rank and trivial opposition involution, Lemma \ref{lem:transverse map criterion} induces a transversality-preserving map of full flag manifolds and Theorem \ref{Thm:MainTheoremTypeB} produces the transverse $\mathbb{S}^3$.

    In the $E_6$ root system, the opposition involution is non-trivial, but its fixed point set is precisely $F_4$. 
    Moreover, this copy of $F_4$ contains regular elements by inspection. 
    Thus, Lemma \ref{lem:transverse map criterion} induces a transversality-preserving map of full flag manifolds and we obtain a transverse 3-sphere in $\Flag(E_6)$.

    There is a $D_8$ root system in the $E_8$ root system, inducing a local embedding of $\SO_0(8,8)$ into $E_8$.
    Since the root systems have the same rank and trivial opposition involution, Lemma \ref{lem:transverse map criterion} induces a transversality-preserving map of full flag manifolds and Theorem \ref{Thm:DivisionAlgebras} completes the proof.
\end{proof}

\begin{figure}[h]
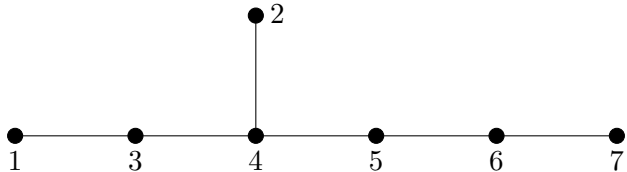

    \centering
    \begin{equation*}
        \begin{dynkinDiagram}[text style/.style={scale=1},
            edge length=.8cm,
            scale=2,
            labels={1,2,3,4,5,6,7},
            label macro/.code={{{#1}}}
            ]E{7}
        \end{dynkinDiagram}
    \end{equation*}
    \caption{\emph{\small{Dynkin diagram of type ${\rm E}_7$. Let $\Delta$ denote the set of simple roots of $E_7$ compatible with the labeling here.}}}
    \label{fig:dynkin diagram E7}
\end{figure}

For $E_7$, we determine when a transverse circle in a partial flag manifold is maximally transverse, and in the other cases we construct a transverse $3$-sphere.
We set some notation for the $E_7$ root system, see Figure \ref{fig:dynkin diagram E7}.

\begin{theorem}\label{thm:E7}
    Let $\Delta$ denote the simple roots of $E_7$, labeled as in Figure \ref{fig:dynkin diagram E7}.
    There exists a transverse $3$-sphere in the flag manifold associated to $\Delta \setminus\{7\}$.
    On the other hand, if $\Theta$ contains $7$, then every transverse circle in $\mathcal{F}_\Theta$ is maximally transverse.
\end{theorem}

\begin{proof}
    The inclusion of $E_6$ into $E_7$ induces a transversality-preserving map from the full flag manifold of $E_6$ to the partial flag manifold $\mathcal{F}_{\Delta \setminus\{7\}}$ of $E_7$ by Lemma \ref{lem:transverse map criterion}, so Theorem \ref{thm:F4,E6,E8} produces a transverse $\mathbb{S}^3$ in $\mathcal{F}_{\Delta \setminus\{7\}}$.

    On the other hand, it is well-known that the fundamental representation of $E_7$ associated to $\varpi_7$ is minuscule, $56$-dimensional, and preserves a symplectic form \cite{cartan1914,Tits1971}.\footnote{For more modern treatment, see Bourbaki. \cite[Ch VIII, \S 3]{Bou05} shows that $\varpi_7$ is minuscule.} 
    It is straightforward to verify that the representation is proximal when restricted to the coweight $\varpi^\vee_7$. 
    Indeed, the weights are $\frac{3}{2},\frac{1}{2},-\frac{1}{2}, -\frac{3}{2}$ of multiplicity $1,27,27,1$ respectively. 
    Lemma \ref{lem:transverse map criterion} then produces a transversality-preserving map $(E_7, \mathcal{F}_{\{7\}}) \rightarrow (\Sp(56, \R), \RP^{55})$. 
    Thus, a transverse circle in $\mathcal{F}_{\{7\}}$ yields a transverse circle in $(\Sp(56,\R),\RP^{55})$ that is maximally transverse by \cite{PT24}. 
\end{proof}

\begin{corollary}\label{Cor:E7Sambarino}
    If $\Gamma < E_7$ is $\{7\}$-Anosov, then $\Gamma$ is virtually isomorphic to a free or surface group.
\end{corollary}

\subsection{Maximal Transversality}

The following proposition uses well-known ideas regarding cohomology and homotopy groups of compact Lie groups (cf. \cite[Ch.  V, \S 12]{Bre93}). 

\begin{proposition}
    Let $K$ be a simply connected compact Lie group and let $f \colon \Sp(1) \to K$ be a nontrivial Lie group homomorphism. 
    Then $[f]$ is nonzero in $\pi_3(K)$.
\end{proposition}

\begin{proof}
    Let $B$ denote the Killing form of $K$, which is necessarily non-degenerate. 
    Then $\omega \in \Lambda^3(\mathfrak{k}^*)$ defined by $\omega(X,Y,Z) = B(X,[Y,Z])$ is $\mathrm{Ad}(G)$-invariant. 
    Considering left-invariant vector fields on $K$, then $\omega$ yields a closed 3-form $\omega \in \Omega^3(K)$ by the same formula  \cite[Proposition 12.6]{Bre93}.
    A closed $3$-form $\omega_1$ is defined on $\Sp(1)$ by the same process.
    By naturality of the definition, the pullback $f^\ast \omega$ agrees with $\omega_1$ up to rescaling by a positive scalar. 
    One checks that $\int_{\mathbb{S}^3} \omega_1$ is nonzero by a direct computation. 
    Hence, $[f^*\omega] \in H^3_{dR}(\mathbb{S}^3)$ is non-trivial, so $f$ is not nullhomotopic. 
\end{proof}

We deduce:
\begin{corollary}
There exists a maximally transverse 3-sphere in $\Flag(F_4)$ and $\Flag(E_6)$. 
\end{corollary}

\medskip

By understanding representations $\mathfrak{so}(8) \rightarrow \mathfrak{so}(16)$, we will be able to verify maximal transversality of 7-spheres in $\Flag(E_8)$. 
The following proposition is the key to this end. 

\begin{proposition}\label{Prop:D8E8}
    Let $\phi \colon D_8 \to E_8$ be a nontrivial Lie group homomorphism. 
    Consider the restriction of the differential to the maximal compact
        \[ d\phi \colon \mathfrak{so}(8) \times \mathfrak{so}(8) \to \mathfrak{so}(16) .\]
    Then $d\phi$ preserves a decomposition $\R^{16} = \R^8 \oplus \R^8$ so that each factor of $\mathfrak{so}(8) \times \mathfrak{so}(8)$ acts non-trivially on precisely one factor of $\R^8 \oplus \R^8$.
\end{proposition}

\begin{proof}
The differential $d\phi$ is injective since $D_8$ is simple and $\phi$ is non-trivial.
Denote by $\rho_i\colon \mathfrak{so}(8) \rightarrow \mathfrak{so}(16)$ the restriction to the $i^{th}$ factor of $\mathfrak{so}(8)$, for $i \in \{1,2\}$. 
Note that $\rho_1 = W_1 \oplus W_2$ splits as a direct sum of two invariant 8-dimensional subspaces. 
One can see this by verifying that $\mathfrak{so}(8)$ has exactly three nontrivial irreducible representations of dimension at most $16$, and these are its three non-isomorphic representations of dimension $8$. 
Since $\rho_1$ is non-trivial, at least one of $W_1, W_2$ is an irreducible sub-representation. 
Without loss of generality, suppose $W_1$ is irreducible. 
Then $W_2$ cannot also be irreducible, or else $\rho_2$ would be trivial, since $\rho_1, \rho_2$ must commute. 
It follows that $\rho = W_1 \oplus W_2$, with $\rho_i$ irreducible on $W_i$. 
\end{proof}

\begin{corollary}
There exists a maximally transverse 7-sphere in $\Flag(E_8)$. 
\end{corollary}

\begin{proof}
Let $\hat{\xi} = (\id, f)\colon \mathbb{S}^7 \rightarrow \SO(8) \times \SO(8)$ be the usual map, with $f\colon \mathbb{S}^7 \rightarrow \SO(\mathbb{O})$ given by $f(x) = L_x$, defining a transverse 7-sphere in $\Flag(D_8)$ by taking the orbit of a point. 
Let $K < E_8$ be the maximal compact subgroup, and $\phi\colon D_8 \rightarrow E_8$ a nontrivial homomorphism. 
The 7-sphere $\Lambda \subset \Flag(E_8)$ defined by $\phi \circ \hat{\xi}$ is transverse by the proof of Theorem \ref{thm:F4,E6,E8};  we need only verify maximal transversality.  

If $\rho\colon \mathfrak{so}(8) \rightarrow \mathfrak{so}(8)$ is any Lie algebra isomorphism, then the unique Lie group homomorphism $\hat{\rho} \colon \mathrm{Spin}(8) \rightarrow \mathrm{Spin}(8)$ with $d\hat{\rho}=\rho$ is an isomorphism on $\pi_7$. 
Applying Proposition \ref{Prop:D8E8}, the inclusion $\{\id\} \times \SO(8) \hookrightarrow K$ is, up to finite cover, and pre-composition by a diffeomorphism of $\SO(8)$, the same as the standard inclusion $\iota\colon \SO(8) \hookrightarrow \SO(16)$. 
Now, $\iota$ satisfies $\iota_*[f] \neq 0 \in \pi_7(\SO(16))$ since $f$ defines a stably non-trivial vector bundle by clutching. 
We conclude that $\phi \circ \hat{\xi}\colon \mathbb{S}^7 \rightarrow K$ is homotopically non-trivial. 
By Fact \ref{Fact:MaximallyTransverse}, the sphere $\Lambda$ maximally transverse in $\Flag(E_8)$. 
\end{proof}

\section{Obstructing Limit Sets}\label{Sec:Obstructions}

In this section, we provide obstruction arguments for subsets of various transverse spheres built in this paper. 
These arguments apply only to the undeformed spheres, i.e., not to the spheres from Section \ref{Sec:Deformations}.  

\subsection{Obstructing Spinor Spheres}\label{Subsec:ObstructSpheres}

In this subsection, we prove many of the higher dimensional spheres constructed in full flag manifolds in Sections \ref{Sec:G2Transversality}, \ref{Subsec:BnTransversality}, \ref{Subsec:DnTransversality},
\ref{Subsec:AnTransversality}
are not realized as Anosov flag limit sets. 

The following result shows that any $m$-sphere in $\Flag(\R^{n})$, for $m \geq 2$, which lifts a copy of $\Ein^{0,m}$ cannot be realized as the flag limit set of a Borel Anosov representation.

\begin{proposition}\label{Prop:ObstructSphere}
Suppose $\Lambda \subset \Flag(\R^{d})$ is a transverse $m$-sphere, $m \geq 2$, and $\mathrm{pr}_{1}: \Flag(\R^d) \rightarrow \mathbb{P}(\R^d)$ is the natural projection. If $\mathrm{pr}_1(\Lambda) =\Ein^{0,m}$, for some quadratic form $q$ on $\R^{d}$, then there is no Borel Anosov representation $\eta: \Gamma \rightarrow \SL(d,\R)$ whose flag limit set is $\Lambda$.
\end{proposition}

\begin{proof}
Suppose, for contradiction, that $\eta: \Gamma \rightarrow \SL(d,\R)$ was such a representation.  
The linear span of $\Ein^{0,m}$ is a vector subspace $V$ of $\R^{d}$ of dimension $m+2$ invariant by $\Gamma$. 
Recall by \cite[Proposition 1.8]{GGKW17a} that a representation $\Gamma \rightarrow G$ is $\Theta$-Anosov if and only if its semisimplification is. Thus, up to replacing $\eta$ by $\eta^{ss}$, we can assume $V$ has a $\eta(\Gamma)$-invariant complement $W$. 
Now, by \cite[Proposition 1.2]{KP22}, a representation $\Gamma \rightarrow \GL(V)$ is Borel Anosov if and only if it has uniformly large eigenvalue gaps for $1 \leq i \leq \dim(V)-1$. By this criterion, we see that $\eta|_V:\Gamma \rightarrow \GL(V)$ must be Borel Anosov since $\eta$ itself is Borel Anosov. 

Now, $\Ein^{0,m}$ is bounded in an affine chart of $\mathbb{P}(V)$, and the interior of the convex hull of $\Ein^{0,m}$ is the contractible open convex subset $\Omega \coloneqq\mathbb{P}Q_+(\R^{1,m+1}) \subset \mathbb{P}(V)$. Thus, $\Gamma$ must preserve $\Omega$ and then act cocompactly on $\Omega$ by \cite[Lemma 8.7]{DGK23}. 
In other words, $\eta|_V: \Gamma \rightarrow \GL(V)$ is a Benoist representation. Up to taking a finite index subgroup $\Gamma'<\Gamma$, we can assume $\eta|_V$ takes values in $\GL_+(V)$. Since a representation is $k$-Anosov if and only if its restriction to any finite index subgroup is \cite[Corollary 1.3]{GW12}, this step presents no issues. Then, as in \cite[Proposition 7.2]{CT20}, we can alter $\eta|_V$ to take values in $\SL(V)$ without altering the Anosov condition. Indeed, define $\eta': \Gamma \rightarrow \SL(V)$ by $\eta'|_V:= \frac{\eta|_V}{\det(\eta|_V)^{1/\dim(V)}}$. Since the ratios of eigenvalues for $\eta|_V$ and $\eta'$ are the same, then $\eta'$ is Borel Anosov by \cite{KP22}. However, now we conclude $\eta|_V: \Gamma \rightarrow \SL(V)$ is not $k$-Anosov for $2 \leq k \leq \frac{\dim(V)}{2}$ by \cite[Corollary 1.4]{CT20}, which is a contradiction since $\dim(V) \geq 4$. 
\end{proof}

The following corollary shows many of the non-deformed higher-dimensional spheres constructed in this paper in full flag manifolds are not realizable as flag limit sets of Borel Anosov representations. 

\begin{corollary}\label{Cor:ObstructSphere}
    Let $\Lambda$ be a transverse undeformed spinor $m$-sphere for $m\geq 2$ from Theorem \ref{Thm:G2Fibration} (c), Theorem \ref{thm:BnTransverseSpheres} \ref{item:Bp p 0 mod 4},\ref{item:Bp p 3 mod 4}, \ref{Thm:DnTransverseSpheres} \ref{item:Dp 0},\ref{item:Dp 1},\ref{item:Dp 3}, or \ref{Thm:AnTransverseSpheres} \ref{item:An 701}. 
    Then $\Lambda$ is not the flag limit set of a Borel Anosov representation $\Gamma \rightarrow G$, for appropriate $G$. 
\end{corollary}

\begin{proof}
In the case $G= \Gtwosplit$ and $G= \SO_0(p,p+1)$, the result is immediate from Proposition \ref{Prop:ObstructSphere} by the fact that the inclusions $\Gtwosplit \hookrightarrow \SL(7,\R)$ and $\SO_0(p,p+1) \hookrightarrow \SL(2p+1, \R)$ preserve the Borel Anosov condition.

The inclusion of a Borel Anosov subgroup of $\SO_0(p,p)$ into $\SL(2p,\R)$ may fail to be Borel Anosov (see Remark \ref{rem:Borel Dn in An}), so we need an alternative argument in this case.
Returning to the argument of Proposition \ref{Prop:ObstructSphere}, we find a subspace $V$ of signature $(1,m+1)$ preserved by $\Gamma$. 
Then the inclusion $\Gamma \to \SO_0(p,p)$ factors through a representation $\Gamma \to O(1,m+1) \times O(p-1,p-m-1)$, and therefore cannot be Borel Anosov.
\end{proof} 

We make a brief remark about one sub-case. 

\begin{remark}\label{rem:Borel Dn in An}
    The inclusion of $\SO_0(n,n)$ into $\SL(2n,\R)$ does not take Borel Anosov subgroups to Borel Anosov subgroups in general.
    This may happen even when $n$ is even and there is a transversality-preserving map of full flag manifolds. 
    It is easy to find cyclic examples. 
    Strikingly, Danciger-Zhang \cite[Theorem 1.3]{DZ19} prove a stronger result: no $\SO_0(n,n)$-Hitchin representation $\eta \colon \pi_1 S \rightarrow \SO_0(n,n)$ has the property that the inclusion $\iota \circ \eta \colon \pi_1S \rightarrow \SL(2n,\R)$ is Borel Anosov. 
\end{remark}

Next, we observe an even stronger restriction in the case of $\Gtwosplit$. 

\begin{proposition}
The stabilizer in $\Gtwosplit$ of $\Ein^{0,3}$ in $\Ein^{2,3}$  is compact.
In particular, every discrete group preserving $\Ein^{0,3}$ is finite.
\end{proposition}

\begin{proof}
It suffices to show that the stabilizer is contained in some copy of the maximal compact subgroup $K < \Gtwosplit$. 
Note that if a subgroup $\Gamma< \Gtwosplit$ preserves $\Ein^{0,3}$, then $\Gamma$ preserves the linear span of $\Ein^{0,3}$ in $\R^{3,4}$, which is a copy of $\R^{1,4}$. 
Thus, $\Gamma <O(3,4) $ preserves the complementary subspace $V\cong \R^{2,0}$. 
However, 
\[ W\coloneqq V\oplus (V \times_{3,4} V) = (V \oplus \R^{1,0}) \in \Gr_{(3,0}^\times(\R^{3,4})= \X_{\Gtwosplit}.\]
Since $\Gamma $ preserves $\times_{3,4}$, then $\Gamma \in \Stab_{\Gtwosplit}(W)$, which is a copy of the maximal compact in $\Gtwosplit$. 
\end{proof}

Since the transverse 3-sphere $ \mathbb{S}^3 \cong \Lambda <\Flag(\Gtwosplit)$ from Theorem \ref{Thm:G2Fibration} projects to a copy of $\Ein^{0,3}$, the only possible $\{\beta\}$-Anosov subgroup $\Gamma < \Gtwosplit$ preserving $\Lambda$ is a finite group. 

\subsection{Obstructing Subsets of Spinor Spheres}\label{Subsec:ObstructCircles}

In this section, we prove that there is no Borel Anosov representation of a (closed) surface group into $\SL(d,\R)$, for $d \in \{-1,0,1\} \bmod 8$, whose flag limit set is contained in one of the (non-deformed) higher-dimensional spheres. In fact, we prove a slightly stronger result via the following technical intermediary lemma. 

\begin{lemma}\label{Lem:TechnicalObstruction}
Let $\Lambda \subset \Flag(\R^d)$ be a transverse subset containing an embedded circle. 
Let $\mathrm{pr}_k:\Flag(\R^d) \rightarrow \Gr_k(\R^d)$ be the natural projection, $V_k \subset \R^d$ be the linear span of $\mathrm{pr}_k(\Lambda)$, and $\mathcal{F}$ the associated filtration:
\[ \mathcal{F}= \big [ V_1 \subseteq V_2\subseteq \cdots\subseteq V_{d} = \R^d \big]. \]  
If the associated graded decomposition $\R^d \cong V_1 \oplus V_2/V_1 \oplus \cdots \oplus V_d/V_{d-1}$ has at least two subspaces of dimension two, then there is no Borel Anosov representation $\Gamma \rightarrow \GL(d,\R)$ whose flag limit set is $\Lambda$.  
\end{lemma}

We follow the ideas of the proof Proposition \ref{Prop:ObstructSphere}, but with a few changes. 

\begin{proof}
Suppose, for contradiction, $\eta: \Gamma \rightarrow \GL(d,\R)$ were such a representation. Then $\eta$ must preserve the linear span of $\mathrm{pr}_k(\Lambda)$, so $\eta$ must preserve the filtration $\mathcal{F}$ of $\R^d$. Now, the semisimplification $\eta^{ss}:\Gamma \rightarrow \GL(d,\R)$ is a still Borel Anosov by \cite{GGKW17a}. By hypothesis, $\eta^{ss}$ has (at least) two distinct sub-representations $\eta_i$ for $i \in \{1,2\}$ on vector subspaces $V_i$ of dimension two. Then up to taking a finite-index subgroup of $\Gamma$, we can assume that $\eta_i$ obtains the form $\eta_i: \Gamma \rightarrow \GL_+(V_i)$, where $\GL_+(V) = \{ A \in \GL(V) \mid \det(A) > 0\}$.
Then let us update the original representation $\eta$ by replacing $\eta|_{V_1\oplus V_2}$ to force this sub-representation to take the form $\eta: \Gamma \rightarrow \SL(V_1\oplus V_2)$ as in the proof of Proposition \ref{Prop:ObstructSphere}. Under this change, $\eta|_{V_1 \oplus V_2}$ is still Borel Anosov. 
Hence, $\Gamma$ must be a surface group by \cite[Theorem 1.6]{CT20}. 
However, the sub-representation on $V_1 \oplus V_2$ is reducible, which contradicts \cite[Theorem 7.2]{CT20}.
\end{proof}

\begin{corollary}\label{Cor:ObstructCircles}
    Let $\Lambda \subset \Flag(G)$ be a subset containing an embedded circle. Suppose $\Lambda$ is contained in an (undeformed) transverse sphere $\Lambda'$ from Theorems \ref{thm:BnTransverseSpheres} \ref{item:Bp p 0 mod 4},\ref{item:Bp p 3 mod 4}, \ref{Thm:DnTransverseSpheres} \ref{item:Dp 0},\ref{item:Dp 1},\ref{item:Dp 3},  \ref{Thm:AnTransverseSpheres} \ref{item:An 701}. Then $\Lambda$ is not the flag limit set of a Borel Anosov subgroup $\Gamma < G$, for appropriate $G$. 
\end{corollary}

\begin{proof}
The point is that we can verify the hypotheses of Lemma \ref{Lem:TechnicalObstruction}. 

Set $d\coloneqq 8k+\epsilon$, for $\epsilon \in \{-1,0,1\}$ and we write $n = 4k$ and $m =4k+\epsilon$. Set also $l \coloneqq \rho(4k)-1$. Then by construction, the 1-Anosov limit set satisfies $\image(\xi_1) \subset \Ein^{0,m_1}$ for some integer $1 \leq m_1 \leq l$. In particular, the subspace $V_1$ spanned by $\image(\xi_1)$ is of dimension $\R^{2+m_1}$. Under the background splitting $\R^{d} = \R^{m} \oplus \R^{n}$ into space + time, under which $\Lambda'$ was constructed, we know $V_1 \cong \R^{1,1+m_1}$. Similarly, $\image(\xi_2) \subset \Pho(\R^{2, m_2})$ for some integer $m_2$ satisfying $m_1 \leq m_2 \leq 2l$. Thus, the subspace $V_2$ spanned by $\image(\xi_2)$ satisfies $V_2 \cong \R^{2,m_2}$. Continuing in this fashion, the ambient space $\R^{8k+\epsilon}$ decomposes into an $\eta(\Gamma)$-invariant filtration $V_1 \subsetneq V_2 \subsetneq \cdots \subsetneq V_{4k-1} \subseteq \cdots \subseteq \R^{d}$. In particular, $V_j \subsetneq V_{j+1}$ for $1 \leq j\leq 4k-2$, as the spacelike dimension always increases.

Note that $\dim(V_1) \geq 3$. We now finish verifying the hypotheses of Lemma \ref{Lem:TechnicalObstruction} with some cases, which completes the proof. \medskip 

\textbf{Case 1}: If $\epsilon \in \{0,1\}$, then we have $\geq 4k$ total sub-representations and hence at most $8k-2$ dimensions to distribute among at least $4k-1$ remaining blocks, after removing $V_1$. Every block must have size two or else we find two one-dimensional blocks, which is impossible. 

\textbf{Case 2}: If $\epsilon = -1$, then there are most $8k-4$ dimensions left after removing $V_1$, and at least $4k-2$ blocks remaining. The argument is finished with similar reasoning as in Case 1. 
\medskip 
\end{proof}

Now, we remark on an alternate proof that the higher dimensional spheres are not realizable as flag limit sets of Anosov subgroups.

\begin{remark}
The reasoning in Corollary \ref{Cor:ObstructCircles} gives another proof of Corollary \ref{Cor:ObstructSphere}. Indeed, by the same reasoning as in the corollary, one finds at least two invariant subspaces $V, W$ of $\R^d$ of size at most two, which already implies that $\Gamma$ is a free or surface group by \cite{CT20}. 
\end{remark}

We have a slightly more general result for other subsets of flag manifolds. 
Recall that the full flag manifolds $\Flag(\R^{n,n+1})$ and $\Flag(\R^{n,n})$ and fiber over $\Flag(\R^n)$ by Propositions \ref{Prop:BnAnFibration}, \ref{Prop:DnAnFibration}. 
More generally, for $n\geq 2$, the non-deformed higher dimensional transverse spheres $\mathbb{S}^{n} \stackrel{\pitchfork}{\hookrightarrow} \Flag(\R^{n,n+\epsilon})$ all sit inside a single fiber of the fibration $\Flag(\R^{n,n+\epsilon})\rightarrow \Flag(\R^n)$, for $\epsilon \in \{0,1\}$. We show any `sufficiently large' transverse subset contained in a single fiber of this fibration is not realized as the flag limit set of a Borel Anosov subgroup. 

\begin{corollary}
Let $n \geq 3$. Suppose that $\Lambda \hookrightarrow \Flag(\R^{n,n+\epsilon})$, for $\epsilon \in\{0,1\}$, is a transverse subset properly containing a circle and contained in a fiber of the fibration $\Flag(\R^{n,n+\epsilon}) \rightarrow \Flag(\R^n)$. Then there is no Borel Anosov representation $\Gamma \rightarrow \SO_0(n,n+\epsilon)$ with flag limit set $\Lambda$. 
\end{corollary}

\begin{proof}
One follows the reasoning in Corollary \ref{Cor:ObstructCircles}, one finds the hypotheses of Lemma \ref{Lem:TechnicalObstruction} are satisfied. In particular, we find the grading associated to $\Lambda$ has at least $n$ distinct blocks, so that at least two of these blocks have dimension at most two. Lemma \ref{Lem:TechnicalObstruction} finishes the proof. 
\end{proof} 

\section{Further Questions}

In this section, we include some open questions that are motivated by the present work. \medskip

We begin with an obvious pair of complementary questions.
Other than some low-dimensional results, nearly all cases remain open.
We pose the questions for the full flag manifold of $\R^d$; of course, analogous questions for other flag manifolds are interesting as well.

\begin{question}\label{q:n_for_d} 
    For each positive integer $n$, what is the smallest number $d=d(n)$ such that there is a transverse $n$-sphere in $\Flag(\R^d)$? 
\end{question}

The answer to Question \ref{q:n_for_d} for $n \in \{2,3\}$ is $d=7$ by \cite{Dey25, TZ24} and our results. 
By Proposition \ref{Prop:FibrationsSpheres}, $d(4) \neq 8$ and $d(8) \neq 16$.
Hence, $d(4) \in \{9,15\}$. 
For $5 \leq n \leq 7$, the answer is $d=15$.
For $n=8$, we prove $d=31$ is attainable. 
However, we cannot presently exclude the candidates $d \in \{17,23,24,25\}$. 

\begin{question}\label{q:d_for_n}
    For each positive integer $d$, what is the largest number $n=n(d)$ such that there is a transverse $n$-sphere in $\Flag(\R^d)$? 
\end{question}

By \cite{Dey25}, the answer to Question \ref{q:d_for_n} is $n=1$ when $d \in \{2,3,4,5,6\} \bmod8$, so the question is only interesting if $d \in \{-1,0,1\} \bmod 8$. 
We have $n(7) = n(8) =3$ and $n(15) = n(16) = 7$.
The existence of the spheres comes from Theorem \ref{Thm:IntroArbitrarilyLargeSpheres}, and the upper bounds come Tsouvalas-Zhu \cite[Theorem 5.3]{TZ24} and Proposition \ref{Prop:FibrationsSpheres}. 
All other cases remain open. 

\medskip

It is natural to ask for the generalization of Theorem \ref{Thm:IntroThmEquivalence} dropping the split condition.

\begin{question}\label{q:general property I equivalence}
    Let $\mathcal{F}_\Theta$ be a self-opposite flag manifold of a semisimple Lie group $G$ which does not satisfy Property (I).
    Does $\mathcal{F}_\Theta$ contain a transverse $2$-sphere?
\end{question}

In a related direction, we ask about a potential sharpening of Theorem \ref{Thm:IntroThmEquivalence}. 

\begin{question}\label{q:transverse circles in 2spheres}
    Suppose $C \subset \mathcal{F}_{\Theta}$ is a transverse circle that is not locally maximally transverse. 
    Is $C$ contained in a transverse 2-sphere? 
\end{question}

It makes sense to refine questions on the existence of transverse spheres by demanding the appearance of a fixed transverse triple $\{x,y,z\}$ in the sphere.

\begin{question}\label{q:maximalsphere_eachcomponent}
    Fix a semisimple Lie group $G$ and self-opposite flag manifold $\mathcal{F}_\Theta(G)$. 
    Suppose that $\{x,y,z\}$ is a transverse subset of $\mathcal{F}_\Theta(G)$. 
    What is the largest dimension of a transverse sphere $S \subset\mathcal{F}_{\Theta}$ such that $\{x,y,z\} \subset S$? 
\end{question}

Let us consider Question \ref{q:maximalsphere_eachcomponent} in the case of $\Gtwosplit$ and $\Theta = \Delta$. 
Returning to the notation in Section \ref{Sec:Anosov}, we have $y \in \Omega \coloneqq \mathscr{C}_x \cap \mathscr{C}_z$. 
Marsh-Rietsch compute that $\Omega$ has 11 connected components, and \cite[Figure 6]{MR02} shows that all but one is moved by the inversion map $\iota \colon \Omega \to \Omega$. 
If $y$ is in one of the 10 components moved by $\iota$, then the answer to Question \ref{q:maximalsphere_eachcomponent} is bounded above by $1$, and it is a nontrivial matter to determine if this upper bound is realized. 
In general, the answer is bounded above by $3$ (which can only be realized in the 11th component). \medskip 

Finally, we ask about filling transverse spheres to negatively curved subsets of symmetric spaces of non-compact type.
When $\phi \colon G' \to G$ is a nontrivial homomorphism of semisimple Lie groups and $G'$ has real rank $1$, the associated symmetric space $G'/K'$ is negatively curved and includes equivariantly as a totally geodesic subspace of $G/K$ \cite{Karpalevic53subgroup,Mostow55decomposition}.
Moreover, the visual boundary of $G'/K'$ includes as a transverse sphere into a partial flag manifold $\mathcal{F}_\Theta(G)$ for some $\Theta$ depending on $\phi$.
The transverse spheres we consider in this paper do not arise in this manner. 
On the other hand, we can ask if there is some abstract negatively curved space embedding into $G/K$ with flag limit set our transverse spheres. 
We offer the following precise formulation, in terms of ``$\Theta$-uniformly regular quasi-isometric embeddings," a notion introduced by Kapovich-Leeb-Porti \cite{KLP17}.
By the higher rank Morse lemma \cite{KLP18b}, the domain of such a map is necessarily Gromov hyperbolic.

\begin{question}\label{Question:Filling}
    Fix $G$ and $\Theta$.
    Let $\Lambda$ be a transverse sphere in $\mathcal{F}_\Theta(G)$.
    Is there a locally compact Gromov hyperbolic metric space $Z$ and a $\Theta$-uniformly regular quasi-isometric embedding $f \colon Z \to G/K$ whose $\Theta$-limit set is $\Lambda$?
\end{question}

By \cite{KLP18b}, a $\Theta$-uniformly regular quasi-isometric embedding $f \colon Z \to G/K$ of a locally compact space $Z$ extends continuously to a transverse map $\partial f  \colon \partial Z \to \Flag_\Theta(G)$, so if moreover every transverse circle in $\mathcal{F}_\Theta$ is locally maximally transverse, it follows that every embedded circle in the Gromov boundary of $Z$ is open, see \cite[Theorem F]{Dey25}.  
In particular, this precludes the possibility of $Z$ being a rank 1 symmetric space other than $\Ha^2$, see \cite[Corollary G]{Dey25}.

On the other hand, the situation when $\Flag_{\Theta}(G)$ contains higher-dimensional transverse spheres is completely different, and the possibility of filling remains, making Question \ref{Question:Filling} natural. 
Let us consider the case of $G= \Gtwosplit$.
There is a transverse spinor 3-sphere $\Lambda$ in $\Flag(\Gtwosplit)$, as well as interesting deformations that remain transverse. 
Moreover there is a $\CH^2$-sub-symmetric space in $\mathbb{X}_{\Gtwosplit}$, but this totally geodesic submanifold is not even $\beta$-regular, and so its 3-sphere boundary yields only a transverse subset of $\Pho^\times$. 
It is noteworthy here that there are no $\SO_0(4,1)$-subgroups of $\Gtwosplit$, or even $\SO_0(3,1)$ (cf. \cite{DMG25} and references therein). 
In other words, there is no totally geodesic $\Ha^4$-sub-symmetric space in $\mathbb{X}_{\Gtwosplit}$ to even attempt to modify to fill $\Lambda$. 
% In the $\SO_0(3,4)$ and $\SO_0(4,4)$ cases, filling of the transverse 3-dimensional spinor spheres, or their deformations, in the respective full flag manifolds $\Flag(\R^{3,4})$ and $\Flag(\R^{4,4})$ also remains an alluring challenge. 

\begin{appendices}

\section{Grothendieck Group Completion}\label{Appendix:GGC}

In this section, we recall some basic definitions and properties regarding the Grothendieck group completion $\mathrm{K}(S)$ of an abelian semigroup $S$, including the categorical perspective through its universal property, some standard constructions of $\mathrm{K}(S)$, and an explicit construction of $\KO(X) \coloneqq \mathrm{K}(\Vect_{\R}(X))$ for $X$ a compact Hausdorff topological space. We will assume $X$ is connected for simplicity. The material in this section is purely expository, and can mostly be found in either \cite{LM89} or \cite{Hat17}. \medskip 

Recall that a \emph{semigroup} $(S, +)$ is a set equipped with a binary operation $+:S\times S \rightarrow S$ that is associative. A \emph{monoid} is a semigroup with a two-sided identity element. 

\begin{definition}
Let $S$ be an abelian semigroup. Then the Grothendieck group completion of $S$ is a pair $(K(S),i)$, where $K(S)$ is an abelian group and $i: S \rightarrow K(S)$ is a semigroup homomorphism, satisfying the following universal property:
For any semigroup homomorphism $\phi: S\rightarrow A$ to an abelian group $A$, there is a unique group homomorphism $\psi: K(S) \rightarrow A$ such that $\phi = \psi\circ i$. 
\end{definition}
The word `the' is apt in the definition: one can easily verify using the universal property that if $(i, K)$ and $(i',K')$ are two Grothendieck group completions of $S$, then the unique group homomorphism $K \rightarrow K'$ associated to the map $i'\colon S \rightarrow K'$ is an isomorphism. 
Be advised that if $(i,K)$ is the Grothendieck group completion of $S$, the map $i\colon S \rightarrow K$ might not be an injection. 
See the examples below. 

There are various models for the abstract group $K(S)$. Here is one such construction \cite{LM89}. 

\begin{proposition}\label{Prop:GGCFirstModel}
    Let $S$ be an abelian semigroup. 
    Define $\mathcal{K}(S) \coloneqq (S\times S)/\sim$, where $(a,b) \sim (c,d)$ if there exists an element $(e,e) \in S^2$ such that $a+d+e=c+b+e$. 
    The pair $(\mathcal{K}(S),i)$, for $i\colon S \rightarrow \mathcal{K}(S)$ the map $s \mapsto [(s,0)]$, satisfies the universal property and hence is the Grothendieck group completion of $S$. 
\end{proposition}

Observe that if $S$ is actually a group, then $(a,b) \sim (c,d)$ if and only if $a-b=c-d$.
In this case, the map $S \rightarrow \mathcal{K}(S)$ by $s\mapsto [(s,0)]$ is a group isomorphism.

The motivating example for the Grothendieck group completion is $S = (\mathbb{N}, +)$. 
Following the construction, one finds $\mathcal{K}(S) \cong \Z$. 
The group homomorphism $\mathcal{K}(\mathbb{N}) \rightarrow \Z$ by $[(m,n)]\mapsto m-n$ is an isomorphism. \medskip 

There are two specific instances of the Grothendieck group completion we are interested in: $S = \Vect(X)$, the monoid of isomorphism classes of real vector bundles on $X$ (under direct sum), and $S=\hat{M}_k$ the monoid of isomorphism classes of $\Z_2$-graded real $\Cl(k)$-modules (under direct sum). 

In the former case, following \cite{Hat17}, we recall the following explicit geometric realization of $\mathrm{K}(\Vect(X))$. 
Here, $\epsilon^i$ denotes a trivial rank $i$ (real) vector bundle on $X$. 
\begin{proposition}
    Let $X$ be a compact Hausdorff space. 
    Then define $\KO(X)\coloneqq \Vect(X)/\sim$, where $[E_1] \sim [E_2]$ if there exists an integer $n\geq 0$ and an isomorphism $E_1 \oplus \epsilon^n \cong E_2 \oplus \epsilon^n$. 
    Then $\KO(X)$ is an abelian group. Denoting $\pi$ as the quotient map $\pi\colon \Vect(X) \rightarrow \KO(X)$, the pair $(\KO(X), \pi)$ satisfies the universal property for $\Vect(X)$. 
    Thus, $(\KO(X),\pi) \cong \mathrm{K}(\Vect(X))$. 
\end{proposition}

\begin{proof}
Define $\Vect_s(X) \coloneq \{ [E] \in \Vect(X) \mid  [E]\sim \epsilon^i\; \text{for some } i\}$. 
It is clear that $\Vect_s(X)$ is a sub-semigroup. 
Thus, $\KO(X) = \Vect(X)/\Vect_s(X)$ is naturally a semigroup. 
Note that $\Vect(X)$ is actually a monoid, with $[\epsilon^0]$ the identity.
Hence, $\pi(\Vect_s(X))$ is the identity in $\KO(X)$. 
The standard fact that every vector bundle $E \rightarrow X$ admits a `complementary' bundle $E'$ such that $E\oplus E' \cong \epsilon^i$ for some $i$ (c.f.\ \cite[Corollary 9.9]{LM89} implies that $[E'] =-[E]$ in $\KO(X)$, so $\KO(X)$ is a group. 

Next, we verify the universal property. 
Suppose $\phi\colon \Vect(X) \rightarrow A$ is a semi-group homomorphism and $\psi\colon \KO(X) \rightarrow A$ is a group homomorphism such that $\phi = \psi \circ i$. 
Observe that $\psi$ must obtain the form $\psi([E]) = \phi(E)$, since $\pi\colon \Vect(X) \rightarrow \KO(X)$ is a surjective monoid homomorphism. 
Thus, $\psi$ is unique if it exists. 
Conversely, one checks that $\phi\colon \Vect(X) \rightarrow A$ must satisfy $\Vect_s(X) \subset \ker(\phi)$, so that $\phi$ must descend to a semi-group homomorphism $\phi\colon \KO(X) \rightarrow A$, which is actually a group homomorphism. 
\end{proof}

\begin{corollary}\label{Cor:KOform}
Every element $[E] \in \KO(X)$ obtains the form $[E] = [\epsilon^i] -[E']$ for some $[E'] \in \KO(X)$. 
\end{corollary}

Now, we define the reduced $\KO$-group, still following Hatcher. 
\begin{definition}
Define $\widetilde{\KO}(X) \coloneqq \Vect(X)/\approx$, where $[E_1] \approx [E_2]$ when there exist integers $m,n \geq 0 $ such that $E_1 \oplus \epsilon^n \cong E_2 \oplus \epsilon ^m$. 
\end{definition}

It is clear that $\widetilde{\KO}(X)$ is a group in its own right. 
The following proposition shows the relationship between $\widetilde{\KO}(X)$ and $\KO(X)$. Recall that we assume $X$ is connected. 

\begin{proposition}
    The natural surjective map $\pi\colon \KO(X) \rightarrow \widetilde{\KO}(X)$ has $\ker(\pi)\cong \Z$. 
    Additionally, every element in $\ker(\pi)$ obtains the form $[\epsilon^n] -[\epsilon^m]$ for some $n, m \geq 0$. 
    Thus, $\KO(X) \cong \widetilde{\KO}(X) \oplus \Z$. 
\end{proposition}
\begin{proof}
By Corollary \ref{Cor:KOform}, taking inverses, every element $x \in \KO(X)$ obtains the form $x=[E]-[\epsilon^n]$ for some non-negative integer $n$. 
Thus, $\pi(x) = \pi([E])$, so $x \in \ker(\pi)$ if and only if $x = [\epsilon^m]-[\epsilon^n]$, for some $m \in \Z_{\geq 0}$. 
Then note that $[\epsilon^m]-[\epsilon^n] = [\epsilon^k] -[\epsilon^j]$ in $\widetilde{\KO}(X)$ if and only if $ m-n= j-k$, so that $\ker(\pi) \cong \Z$. The short exact sequence 
\[1 \rightarrow \Z \cong \ker(\pi) \rightarrow \KO(X) \stackrel{\pi}{\rightarrow}{\widetilde{\KO}(X)} \rightarrow 1\]
splits on the left. 
Indeed, consider the map $\deg: \KO(X) \rightarrow \Z$ by $[E] \mapsto \deg(E)$, well-defined by definition of $\KO(X)$. Then $\deg |_{\ker \pi}$ is an isomorphism. Hence, $\KO(X) \cong \widetilde{\KO}(X) \oplus \Z$. 
\end{proof}
It is important in $\mathrm{K}$-theory that $\KO(X)$ and $\widetilde{\KO}(X)$ can be turned into rings using tensor products of bundles.
However, we do not need the product structure, so we do not pursue the details here. \medskip 

We now turn to the study of (graded) Clifford modules and their Grothendieck group completions.
Recall that $\hat{M}_n$ (resp.\ $M_n$) denotes the monoid of isomorphism classes of $\mathbb{Z}_2$-graded (resp.\ ungraded) $\Cl(n)$-modules and $\hat{\mathfrak{M}}_n\coloneqq \mathrm{K}(\hat{M}_n)$ (resp.\ $\mathfrak{M}_n\coloneqq \mathrm{K}(M_n)$) denotes the Grothendieck group completion.

There is a natural isomorphism $\hat{M}_n \cong M_{n-1}$ \cite[Proposition 5.20]{LM89}. Indeed, the monoid isomorphism $\hat{M}_n \rightarrow M_{n-1}$ is given by the map $[S] \mapsto [S^0]$, where $S^{\bullet}$ is a $\Z_2$-graded $\Cl(n)$-module and $S^0$ is the associated $\Cl(n-1)\cong \Cl^0(n)$-sub-module.
Using this isomorphism, we can easily compute the groups $\mathfrak{M}_n\coloneqq \mathrm{K}(M_n)\cong \mathrm{K}(\hat{M}_{n+1}) \eqqcolon \hat{\mathfrak{M}}_{n+1}$ found in Table \ref{Table:KTheoryBottPeriodicity}. 

Here, we note a nice realization of $\mathfrak{M}_n$ from \cite{LM89}. 
In fact, $\mathfrak{M}_n$ is isomorphic to a free abelian group with one generator for each distinct irreducible representation of $\Cl(n)$. 
To see this, consider the model $\mathcal{K}(M_n)$ for $\mathrm{K}(M_n)$ constructed in Proposition \ref{Prop:GGCFirstModel}. 
Note that if $\omega_1, \dots, \omega_k$ are a maximal set of pairwise non-isomorphic irreducible representations, then $M_k \cong \Z_+^k$. 
Indeed, (pairwise non-isomorphic and irreducible) implies injectivity and maximality implies surjectivity of the map $\Z_+^k \rightarrow M_k$ by $(a^i) \mapsto a^i\omega_i$. 
Hence, $\mathrm{K}(M_k) \cong \Z^k$, with $[(\omega_1,0)],\dots, [(\omega_k,0)]$ as generators, written in terms of the model $\mathcal{K}(M_k)$.

We now compute $\mathfrak{M}_n$. 
Note that $\Mat_n(\R)\otimes \Mat_m(\R) \cong \Mat_{mn}(\R)$ and $\Mat_n(\K)\cong_{\R-\mathrm{alg}} \K \otimes_{\R} \Mat_n(\R)$. 
By Table \ref{Table:Clifford}, along with periodicity $\Cl(n+8) \cong \Cl(n) \otimes_{\R}\Cl(8)$, we find the structure of $\Cl(n)$ as an $\R$-algebra for all positive integers $n$, as is well-known. 
\begin{align}\label{CliffordStructure}
\begin{cases}
    \Cl(8k) &\cong \Mat_{16^k}(\R) \\
    \Cl(8k+1) &\cong \Mat_{16^k}(\C) \\
    \Cl(8k+2) &\cong \Mat_{16^k}(\Ha) \\
    \Cl(8k+3)  &\cong \Mat_{16^k}(\Ha)\oplus \Mat_{16^k}(\Ha)\\
    \Cl(8k+4)  &\cong \Mat_{2\cdot 16^k}(\Ha)\\
    \Cl(8k+5)  &\cong \Mat_{4\cdot 16^k}(\C)\\
    \Cl(8k+6)  &\cong \Mat_{8\cdot 16^k}(\R)\\
    \Cl(8k+7)  &\cong \Mat_{8\cdot 16^k}(\R)\oplus \Mat_{8\cdot 16^k}(\R)
\end{cases}
\end{align}
We see that, when $n \not \equiv 3 \bmod 4$, $\Mat_n(\K)$ has a unique irreducible $\K$-representation, and also a unique irreducible $\R$-representation (up to isomorphism), and conclude that $\mathfrak{M}_n \cong \Z$. 
When $n \equiv 3 \bmod 4$, since $\Mat_n(\K) \oplus \Mat_n(\K)$ admits two distinct irreducible real representations, one from each factor, the Grothendieck group is $\mathfrak{M}_n \cong \Z\oplus \Z$ for $n \equiv 3 \bmod 4$. 
\end{appendices}

\bibliographystyle{alphaurl}
\bibliography{biblio}

\begin{thebibliography}{GGKW17}

\bibitem[ABS64]{ABS64}
Michael Atiyah, Raul Bott, and Arnold Shapiro.
\newblock {Clifford} modules.
\newblock {\em Topology}, 3:3--38, 1964.

\bibitem[Ada60]{Ada60}
J.~F. Adams.
\newblock On the {Non-Existence} of {Elements} of {Hopf Invariant One}.
\newblock {\em Ann. of Math.}, 72(1):20--104, 1960.

\bibitem[Bae02]{Bae02}
John Baez.
\newblock The {Octonions}.
\newblock {\em Bull. Amer. Math. Soc.}, 39:145--205, 2002.

\bibitem[BH14]{BH14}
John~C. Baez and John Huerta.
\newblock {$G_2$} and the rolling ball.
\newblock {\em Trans. Amer. Math. Soc.}, 366(10):5257--5293, 2014.

\bibitem[BK25]{BK25}
Jonas Beyrer and Fanny Kassel.
\newblock $\mathbb{H}^{p,q}$-convex cocompactness and higher higher {Teichm\"uller} spaces.
\newblock {\em Geom. Func. Anal. (to appear)}, 2025.
\newblock \href {https://arxiv.org/abs/2305.15031} {\path{arXiv:2305.15031}}.

\bibitem[Bou05]{Bou05}
N.~Bourbaki.
\newblock {\em {Lie Groups and Lie Algebras-- Chapters 7-9}}.
\newblock {Springer Berlin, Heidelberg}, 2005.

\bibitem[Bre93]{Bre93}
Glen Bredon.
\newblock {\em Topology and Geometry}.
\newblock Springer New York, NY, 1993.

\bibitem[Bry20]{Bry20}
Robert~L. Bryant.
\newblock Notes on spinors in low dimension, 2020.
\newblock \href {https://arxiv.org/abs/2011.05568} {\path{arXiv:2011.05568}}.

\bibitem[BS10]{bonsante2010maximal}
Francesco Bonsante and Jean-Marc Schlenker.
\newblock Maximal surfaces and the universal {T}eichm{\"u}ller space.
\newblock {\em Invent. Math.}, 182(2):279--333, 2010.

\bibitem[Car14]{cartan1914}
Elie Cartan.
\newblock Les groupes r\'eels simples, finis et continus.
\newblock {\em Ann. Sci. \'Ec Norm. Sup. (3)}, 31:263--355, 1914.

\bibitem[CT20]{CT20}
Richard Canary and Konstantinos Tsouvalas.
\newblock Topological restrictions on {Anosov} representations.
\newblock {\em J. Topol.}, 13(4):1497--1520, 2020.

\bibitem[CT24]{CT24}
Brian Collier and J{\'e}r{\'e}my Toulisse.
\newblock Holomorphic curves in the 6-pseudosphere and cyclic surfaces.
\newblock {\em Trans. Amer. Math. Soc.}, pages 6465--6514, 2024.

\bibitem[DE25]{DavEva25}
Colin Davalo and Parker Evans.
\newblock Geometric structures for $\mathrm{SO}_0(p,p+1)$ and $\mathrm{G}_2'$-representations (to appear), 2025.

\bibitem[Dey25]{Dey25}
Subhadip Dey.
\newblock On {Borel} {Anosov} subgroups of $\mathrm{SL}(d, \mathbb{R})$.
\newblock {\em Geom. Topol.}, 29(1):171--192, 2025.

\bibitem[DGK18]{DGK18}
Jeffrey Danciger, Fran\c{c}ois Gu\'{e}ritaud, and Fanny Kassel.
\newblock Convex cocompactness in pseudo-{R}iemannian hyperbolic spaces.
\newblock {\em Geom. Dedicata}, 192:87--126, 2018.

\bibitem[DGK23]{DGK23}
Jeffrey Danciger, François Guéritaud, and Fanny Kassel.
\newblock Convex cocompact actions in real projective geometry.
\newblock {\em Ann. Sci. \'Ec Sup. (to appear)}, 2023.
\newblock \href {https://arxiv.org/abs/1704.08711} {\path{arXiv:1704.08711}}.

\bibitem[DGR24]{DGR24}
Subhadip Dey, Zachary Greenberg, and J.~Maxwell Riestenberg.
\newblock Restrictions on {A}nosov subgroups of {${\rm Sp}(2n,\Bbb R)$}.
\newblock {\em Trans. Amer. Math. Soc.}, 377(10):6863--6882, 2024.

\bibitem[DK23]{DK23}
Subhadip Dey and Michael Kapovich.
\newblock {Klein–Maskit} combination theorem for {Anosov} subgroups: free products.
\newblock {\em Math. Z.}, 305:35, 2023.

\bibitem[DMG25]{DMG25}
Cristina Draper and Cándido Martín-González.
\newblock A perspective on totally geodesic submanifolds of the symmetric space {$\mathrm{G}_2/\SO(4)$}, 2025.
\newblock URL: \url{https://arxiv.org/abs/2504.07586}, \href {https://arxiv.org/abs/2504.07586} {\path{arXiv:2504.07586}}.

\bibitem[{Dra}18]{Fon18}
Cristina {Draper Fontanals}.
\newblock Notes on ${G}_2$: The {Lie} algebra and the {Lie} group.
\newblock {\em Differential Geom. Appl.}, 57:23--74, 2018.
\newblock (Non)-existence of complex structures on $\mathbb{S}^6$.

\bibitem[DZ19]{DZ19}
Jeffrey Danciger and Tengren Zhang.
\newblock {Affine actions with Hitchin linear part}.
\newblock {\em Geom. Func. Anal.}, 29(5):1369--1439, 2019.

\bibitem[Ebe96]{eberlein96book}
Patrick~B. Eberlein.
\newblock {\em Geometry of nonpositively curved manifolds}.
\newblock Chicago Lectures in Mathematics. University of Chicago Press, Chicago, IL, 1996.

\bibitem[Esc18]{Esc18}
J.H. Eschenburg.
\newblock {Geometry} of {Octonions}, 2018.
\newblock \href {https://arxiv.org/abs/https://myweb.rz.uni-augsburg.de/~eschenbu/geomoct.pdf} {\path{arXiv:https://myweb.rz.uni-augsburg.de/~eschenbu/geomoct.pdf}}.

\bibitem[Eva24a]{Eva24Thesis}
Parker Evans.
\newblock {\em Geometry of ${G}_2'$-Harmonic Maps and Representations}.
\newblock PhD thesis, Rice University, 2024.

\bibitem[Eva24b]{Eva24}
Parker Evans.
\newblock Polynomial almost-complex curves in $\hat{S}^{2,4}$, 2024.
\newblock \href {https://arxiv.org/abs/2208.14409} {\path{arXiv:2208.14409}}.

\bibitem[FG06]{FG06}
V.~Fock and A.~Goncharov.
\newblock Moduli spaces of local systems and higher {Teichmüller} theory.
\newblock {\em Publ. Math. IHES}, 103:1–211, 2006.

\bibitem[Fro77]{Fro78}
G.~Frobenius.
\newblock Ueber lineare substitutionen und bilineare formen.
\newblock {\em J. Reine Angew. Math.}, 84:1--63, 1877.

\bibitem[GGKW17]{GGKW17a}
Fran\c{c}ois Gu\'{e}ritaud, Olivier Guichard, Fanny Kassel, and Anna Wienhard.
\newblock Anosov representations and proper actions.
\newblock {\em Geom. Topol.}, 21(1):485--584, 2017.

\bibitem[GSV03]{GSV03}
Michael Gekhtman, Michael Shapiro, and Alek Vainshtein.
\newblock The number of connected components in double {B}ruhat cells for nonsimply-laced groups.
\newblock {\em Proc. Amer. Math. Soc.}, 131(3):731--739, 2003.

\bibitem[GW12]{GW12}
Olivier Guichard and Anna Wienhard.
\newblock Anosov representations: domains of discontinuity and applications.
\newblock {\em Invent. Math.}, 190(2):357--438, 2012.

\bibitem[Hat17]{Hat17}
Alan Hatcher.
\newblock Vector {Bundles} and {$K$}-theory, 2017.
\newblock URL: \url{https://pi.math.cornell.edu/~hatcher/VBKT/VB.pdf}.

\bibitem[Hei05]{Hei05}
Juha Heinonen.
\newblock {\em Lectures on {L}ipschitz analysis}, volume 100 of {\em Report. University of Jyv\"askyl\"a{} Department of Mathematics and Statistics}.
\newblock University of Jyv\"askyl\"a, Jyv\"askyl\"a, 2005.

\bibitem[HL82]{HL82}
Reese Harvey and H.~Blaine Lawson.
\newblock {Calibrated geometries}.
\newblock {\em Acta Math.}, 148:47 -- 157, 1982.

\bibitem[Hur98]{Hur98}
Adolf Hurwitz.
\newblock {Ueber die Composition der quadratischen Formen von belibig vielen Variablen}.
\newblock {\em Nachr. Ges. Wiss. G\"ottingen}, pages 309--316, 1898.

\bibitem[Kar53]{Karpalevic53subgroup}
F.~I. Karpelevi\v{c}.
\newblock Surfaces of transitivity of a semisimple subgroup of the group of motions of a symmetric space.
\newblock {\em Doklady Akad. Nauk SSSR (N.S.)}, 93:401--404, 1953.

\bibitem[Kar20]{Kar20}
Spiro Karigiannis.
\newblock {\em Introduction to {$G_2$} geometry}, pages 3--50.
\newblock Springer US, 2020.

\bibitem[KLP17]{KLP17}
Michael Kapovich, Bernhard Leeb, and Joan Porti.
\newblock Anosov subgroups: dynamical and geometric characterizations.
\newblock {\em Eur. J. Math.}, 3(4):808--898, 2017.

\bibitem[KLP18]{KLP18b}
Michael Kapovich, Bernhard Leeb, and Joan Porti.
\newblock A {M}orse lemma for quasigeodesics in symmetric spaces and euclidean buildings.
\newblock {\em Geom. Topol.}, 22(7):3827--3923, 2018.

\bibitem[KP22]{KP22}
Fanny Kassel and Rafael Potrie.
\newblock Eigenvalue gaps for hyperbolic groups and semigroups.
\newblock {\em J. Mod. Dyn.}, 18:161--208, 2022.

\bibitem[KR71]{KR71}
Bertram Kostant and Stephen Rallis.
\newblock Orbits and representations associated with symmetric spaces.
\newblock {\em Amer. J. Math.}, 93(3):753--809, 1971.

\bibitem[KT24]{KT24}
Clarence Kineider and Roméo Troubat.
\newblock Connected components of the space of triples of transverse partial flags in $\mathrm{SO}_0(p,q)$ and {A}nosov representations, 2024.
\newblock \href {https://arxiv.org/abs/2411.08679} {\path{arXiv:2411.08679}}.

\bibitem[Lab06]{Lab06}
Fran\c{c}ois Labourie.
\newblock Anosov flows, surface groups and curves in projective space.
\newblock {\em Invent. Math.}, 165(1):51--114, 2006.

\bibitem[Leb87]{Leb87}
Claude Lebrun.
\newblock Orthogonal complex structures on $\mathbb{S}^6$.
\newblock {\em Proc. Amer. Math. Soc.}, 101, 1987.

\bibitem[LM89]{LM89}
H.~Blaine Lawson and Marie-Louise Michelsohn.
\newblock {\em Spin Geometry}.
\newblock Princeton University Press, 1989.

\bibitem[Mes07]{Mes07}
G~Mess.
\newblock {Lorentz} spacetimes of constant curvature.
\newblock {\em Geom. Dedicata}, pages 3--45, 2007.

\bibitem[MNS21]{MNS21}
Gianni Manno, Pawe{\l} Nurowski, and Katja Sagerschnig.
\newblock The {Geometry} of {Marked} {Contact} {Engel} {Structures}.
\newblock {\em J. Geom. Anal.}, 31(8):7686--7708, 2021.

\bibitem[Mos55]{Mostow55decomposition}
G.~D. Mostow.
\newblock Some new decomposition theorems for semi-simple groups.
\newblock {\em Mem. Amer. Math. Soc.}, 14:31--54, 1955.

\bibitem[MR02]{MR02}
R.~J. Marsh and K.~Rietsch.
\newblock The intersection of opposed big cells in the real flag variety of type ${G}_2$.
\newblock {\em Proc. Lond. Math. Soc.}, 85(1):22–42, 2002.

\bibitem[Nie24]{Nie24}
Xin Nie.
\newblock Cyclic {H}iggs bundles and minimal surfaces in pseudo-hyperbolic spaces.
\newblock {\em Adv. Math.}, 436:109402, 2024.

\bibitem[PSW21]{PSW21}
Maria~Beatrice Pozzetti, Andr\'es Sambarino, and Anna Wienhard.
\newblock Conformality for a robust class of non-conformal attractors.
\newblock {\em J. Reine Angew. Math.}, 774:1--51, 2021.

\bibitem[PT24]{PT24}
Maria~Beatrice Pozzetti and Konstantinos Tsouvalas.
\newblock On projective {A}nosov subgroups of symplectic groups.
\newblock {\em Bull. Lond. Math. Soc.}, 56(2):581--588, 2024.

\bibitem[San25]{San25}
Jesus Sanchez.
\newblock {Explicit Families of Spinor Representations}, 2025.
\newblock \href {https://arxiv.org/abs/2505.16203} {\path{arXiv:2505.16203}}.

\bibitem[Sch95]{Sch95}
R.D. Schafer.
\newblock {\em {Introduction} to {Non-Associative Algebras}}.
\newblock Dover, 1995.

\bibitem[Tit71]{Tits1971}
J.~Tits.
\newblock Représentations linéaires irréductibles d'un groupe réductif sur un corps quelconque.
\newblock {\em J. Reine Angew. Math.}, 247:196--220, 1971.

\bibitem[Tso20]{tsouvalas2020borel}
Konstantinos Tsouvalas.
\newblock On {Borel} {Anosov} representations in even dimensions.
\newblock {\em Comment. Math. Helv.}, 95(4):749--763, 2020.

\bibitem[TZ24]{TZ24}
Konstantinos Tsouvalas and Feng Zhu.
\newblock Topological restrictions on relatively {Anosov} representations.
\newblock {\em Trans. Amer. Math. Soc. (to appear)}, 2024.
\newblock \href {https://arxiv.org/abs/2401.03050} {\path{arXiv:2401.03050}}.

\bibitem[Yok77]{Yok77}
Ichiro Yokota.
\newblock Non-compact simple {Lie} group {$G_2'$} of type {$G_2$}.
\newblock {\em Jour. Fac. Sci., Shinsu University}, 12(1), 1977.

\bibitem[Zel00]{Zel00}
Andrei Zelevinsky.
\newblock {Connected Components of Double Bruhat Cells}.
\newblock {\em Int. Math. Res. Not.}, 2000(21):1131--1154, 2000.

\bibitem[Zor31]{Zor31}
Max Zorn.
\newblock Theorie der alternativen ringe.
\newblock {\em Abh. Math. Semin. Univ. Hambg.}, pages 123--147, 1931.

\end{thebibliography}

\end{document}